\magnification=\magstep1
\hfuzz=2pt
\hfuzz=2.2pt
\headline={\ifnum\pageno=1 \hfil\else\hss{\tenrm\folio}\hss\fi}
\footline={\hfil}
\hoffset=-.2cm
\font\titlefont= cmbx10 scaled \magstep1
\def\<{\langle}
\def\>{\rangle}
\def\I{I^{\!^{\!^{\!\circ}}}}
\def\J{J^{\!^{\!^{\!\circ}}}}
\def\disk{\rm I\hskip-1.9pt D}
\def\plex{{\rm I\kern-0.6em C}}
\def\integer{\mathchoice{\rm I\hskip-1.9pt N}{\rm I\hskip-1.9pt N}
{\rm I\hskip-1.4pt N}{\rm I\hskip-.5pt N}}
\def\real{\rm I\hskip-1.9pt R}

\def\Romannumeral#1{\uppercase\expandafter{\romannumeral#1}}
\def\date{\line{\number\day/\number\month/\number\year\hfil}}

\def\bz{{Z\kern-0.5em Z}}

\newskip\refskip\refskip=3em
\def\refsize{\advance\leftskip by \refskip}
\def\ref#1#2{\noindent\hskip -\refskip\hbox to
\refskip{[#1]\hfil}{\noindent #2\hfil}\smallskip}
\def\Bbb#1{{\bf Z}}


\def\epsilon{{\varepsilon}}
\def\sqr#1#2{{\vcenter{\vbox{\hrule height.#2pt 
	\hbox{\vrule width.#2pt height#1pt \kern#1pt 
	     \vrule width.#2pt} 
	\hrule height.#2pt}}}} 
\def\square{\mathchoice\sqr35\sqr35\sqr{2.1}4\sqr{1.5}4}
\def\endproof{{\quad$\square$
               \medskip}}
\def\emptyset{\hbox{\O}}
\def\centre#1{{\centerline{#1}}}
\def\ni{{\noindent}}
\font\tensmc=cmcsc10
\def\smc{\tensmc}
\def\proclaim#1.#2{{\medbreak\ni{\bf #1.\enspace}
    {\it #2\par} \ifdim \lastskip <\medskipamount
    \removelastskip \penalty 55\medskip \fi}}
\def\sproclaim#1.#2.#3{{\medbreak\ni{\bf #1.\enspace}
    {\it #2\par}{\it #3\par}\ifdim \lastskip <\medskipamount
    \removelastskip \penalty 55\medskip \fi}}
\def\proof{\noindent{\it Proof.}\ }
\def\bin#1.#2{{\left({{\scriptstyle #1}
    \atop{\scriptstyle #2}}\right)}}

\def\medtype{\font\eightrm=cmr8 
	\font\eightbf=cmbx8 
	\font\eightit=cmti8 
	\font\eightsl=cmsl8 
	\font\eightmus=cmmi8
	\font\eightsmc=cmcsc8
	\font\sixmus=cmmi6
	\font\sixrm=cmr6
	\let\rm=\eightrm \let\it=\eightit \let\sl=\eightsl 
	\let\bf=\eightbf \let\mus=\eightmus 
	\textfont0=\eightrm  \scriptfont0=\sixrm
	\textfont1=\eightmus  \scriptfont1=\sixmus
	\rm}

\def\d{\displaystyle}
\def\Nkn#1.#2{{N^{#1,#2}}}
\def\fkn#1.#2{{F^{#1}_{#2}}}
\def\fknt#1.#2{{F^{#1}_{#2}(t)}}

\input psbox
\psfordvips
\let\fillinggrid=\relax
\def\SBIMSMark#1#2#3{
 \font\SBF=cmss10 at 10 true pt
 \font\SBI=cmssi10 at 10 true pt
 \setbox0=\hbox{\SBF Stony Brook IMS Preprint \##1}
 \setbox2=\hbox to \wd0{\hfil \SBI #2}
 \setbox4=\hbox to \wd0{\hfil \SBI #3}
 \setbox6=\hbox to \wd0{\hss
             \vbox{\hsize=\wd0 \parskip=0pt \baselineskip=10 true pt
                   \copy0 \break%
                   \copy2 \break%
                   \copy4 \break}}
 \dimen0=\ht6   \advance\dimen0 by \vsize \advance\dimen0 by 8 true pt
                \advance\dimen0 by -\pagetotal
 \dimen2=\hsize \advance\dimen2 by .25 true in
%
%
  \openin2=publishd.tex
  \ifeof2\setbox0=\hbox to 0pt{}
  \else 
     \setbox0=\hbox to 3.1 true in{
                \vbox to \ht6{\hsize=3 true in \parskip=0pt  \noindent  
                \input publishd.tex 
                \vfill}}
  \fi
  \closein2
  \ht0=0pt \dp0=0pt
 \ht6=0pt \dp6=0pt
 \setbox8=\vbox to \dimen0{\vfill \hbox to \dimen2{\copy0 \hss \copy6}}
 \ht8=0pt \dp8=0pt \wd8=0pt
 \copy8
 \message{*** Stony Brook IMS Preprint #1, #2 ***}
}

\SBIMSMark{1996/13}{November 1996}{}
\vglue 1.5cm
\centerline{\titlefont ASYMPTOTIC RIGIDITY OF SCALING RATIOS }
\vskip .5cm
\centerline{\titlefont FOR CRITICAL CIRCLE MAPPINGS}

\vskip 1cm
\centerline{{\smc Edson de Faria}}
\vskip 1cm
{\medtype
\centerline{Instituto de Matem\'atica e Estat\'\i stica,}
\centerline{Universidade de S\~ao Paulo, Caixa Postal 66281,}
\centerline{05389-970 S\~ao Paulo SP - Brazil}}
\vfootnote*{{\font\sevenrm=cmr8
\sevenrm
This work has been supported by FAPESP Grant 95/3187-4.}}
\vskip 3.8cm
{
\noindent{\bf Abstract.} \medtype
Let $f$ be a smooth homeomorphism of the circle having one cubic-exponent
critical point and irrational rotation number of bounded combinatorial type.
Using certain pull-back and quasi-conformal surgery techniques, we prove that
the scaling ratios of $f$ about the critical point are asymptotically
independent of $f$. This settles in particular the {\it golden mean
universality conjecture}.  
We introduce the notion of {\it holomorphic commuting pair}, a complex
dynamical system that, in the analytic case, represents an extension
of $f$ to the complex plane and behaves somewhat as a quadratic-like
mapping. We define a suitable renormalization operator that acts on
such objects. Through careful analysis of the family of entire
mappings given by $z\mapsto z+\theta -{1\over {2\pi}}\sin{2\pi z}$,
$\theta$ real, we construct examples of holomorphic commuting pairs,
from which certain necessary limit set pre-rigidity results are
extracted. 
The rigidity problem for $f$ is thereby reduced to one of
renormalization convergence. We handle this last problem by means of
Teichm\"uller extremal methods made available through the recent work of
Sullivan on Riemann surface laminations and renormalization of unimodal
mappings. 
  
}
\vfill
{\medtype
\ni

\ni
 
\smallskip
\ni E.~de Faria e-mail: edson@ime.usp.br
}

\vfill\eject
\centerline{\smc Introduction}
\smallskip
\ni 
The problem of describing the fine scale geometric structure of
one-dimensional dynamical systems has been the subject of intense
investigation in recent years. A fairly complete theory has
emerged through the work of Herman [H$_1$] in the case of smooth
diffeomorphisms of the circle. Under reasonable smoothness
assumptions, Herman showed that any diffeomorphism $f:S^1\to
S^1$ with diophantine rotation number is differentiably conjugate to a
rotation. In particular, the scaling structure of the orbits of $f$
is asymptotically rigid, and completely determined by its rotation
number. Herman's results were subsequently sharpened by Yoccoz [Yo$_1$]
and Katznelson \& Ornstein [KO], and his proofs simplified
in some cases with the help of renormalization methods, as in the
works of Stark [St], Khanin \& Sinai [KS] and Rand [Ra$_3$].

No smooth classification theory as complete as this one exists yet for other
non-expanding one-dimensional dynamical systems. When
critical points are present, the classical Denjoy estimates used by
Herman are no longer sufficient to control the non-linearity of
iterates, and even simple bounds on the geometry of orbits seem to require
these techniques to be used in conjunction with the cross-ratio distortion
tools introduced by Yoccoz, de Melo \& van Strien, \'Swi\c{a}tek and
Sullivan, among others (see [MS] for a historical account). In a
recent {\it tour-de-force} by Sullivan, the asymptotic scaling
structure of the critical orbit of an infinitely renormalizable,
quadratic-like unimodal mapping of the interval was shown to be a
universal function of its kneading invariant, in the cases where such
invariant is of bounded type (cf.~[S$_1$],[MS]). 

In this work we study the scaling problem for the simplest smooth,
non-expanding dynamical systems on the circle besides
diffeomorphisms, namely smooth homeomorphisms with exactly one
critical point. These are called {\it critical circle mappings}. The
prototypical examples are the mappings in the Arnold family,
$$
x\mapsto x+\theta -{1\over 2\pi}\sin{2\pi x}\ \ ({\rm mod}\;1).
\eqno{(1)}
$$
The topological classification of such mappings is just as interesting
as that of diffeomorphisms. As Hall showed in [Ha], Denjoy
examples exist among critical circle mappings with flat critical
points. If the critical point is non-flat, however, then a
topological conjugacy to the corresponding rotation always exists
[Yo$_2$], provided the rotation number is irrational. This
conjugacy can of course never be smooth. Herman and \'Swi\c{a}tek have
shown that it is quasisymmetric if and only if the rotation number 
is an irrational of bounded combinatorial type (this is still
unpublished, but see [Sw$_1$], [H$_2$]). 
On the other hand, Khanin proved in [Kh] that in the unbounded type
case the conjugacy is always purely singular with respect to Lebesgue
measure. These facts reinforce the idea that critical circle mappings should
be compared to each other, not with rotations, and motivate the following
intrinsic rigidity question: 
{\it Are any two topologically conjugate smooth critical circle mappings
always smoothly conjugate?} To give a precise meaning to this question, let us
agree from this point on that a map is {\it smooth} if it is differentiable of
class at least $C^3$ away from critical points. Let us also say that the 
critical point of a critical circle mapping $f$ {\it has type} $s>1$
if $f$ is locally $C^3$-conjugate to $x\mapsto x|x|^{s-1}+a$, for
some $a$, in a neighborhood of the critical point. It is clear that having
critical points of the same type, if their types are defined at all, is a
necessary condition for two smooth critical circle mappings to be smoothly
conjugate. 
The following is supported by numerical observations and
analogy with the unimodal case.
\proclaim{Conjecture}.{Any two smooth critical circle mappings with
the same irrational rotation number and the same type of critical
point are $C^{1+\beta}$-conjugate for some $0\le\beta <1$.}
In this paper we take a step towards proving this conjecture for
rotation numbers of bounded combinatorial type and
critical points of cubic type ($s= 3$). Further steps are taken in
[dFM$_1$] and [dFM$_2$].
Our methods can be adapted to cover
all odd exponents $s=2k+1$, $k\ge 1$, as well. 

We need a few definitions before we can state our results.
Let $f:S^1\to S^1$ be a critical
circle mapping with critical point $c$ and let 
$$
\rho(f)\;=\;[r_0,r_1,\cdots\,]
\;=\;{1\over{\displaystyle r_0 +{\strut 1\over{
\displaystyle r_1 +{\strut 1\over\cdots}}}}}\;. 
$$
be the continued-fraction development of its rotation number. We say
that $\rho(f)$ is a number of {\it bounded combinatorial type} if
$\max\,r_n <\infty$. Let $\{q_n\}_{n\ge 0}$ be the sequence of
return times of the forward orbit of c to itself, and for each $n\ge
0$ let $J_n$ be the closed interval on the circle with endpoints
$f^{q_n}(c)$ and $f^{q_{n+1}}(c)$ that contains $c$. Then $c$ divides
$J_n$ into two intervals, $I_n$ with endpoint $f^{q_n}(c)$ and
$I_{n+1}$ with endpoint $f^{q_{n+1}}(c)$. The ratio of lengths
$s_n(f)\,=\,|I_{n+1}|/|I_n|$ is called the {\it $n$-th scaling ratio}
of $f$.

Our main theorem is the analogue for critical circle mappings of the
Coullet-Tresser rigidity of infinitely-renormalizable Cantor
attractors of unimodal mappings proved by Sullivan in [S$_1$]. We
give here two equivalent versions of this result.

\proclaim{Theorem A}.{If $f$ and $g$ are smooth critical circle
mappings with the same irrational rotation number of bounded
combinatorial type, then they have asymptotically the same scaling
ratios, i.e.
$$
\lim_{n\to\infty}{s_n(f)\over{s_n(g)}}\;=\;1\ .\eqno\square
$$
}
This theorem improves upon the so-called real a-priori bounds for critical
circle mappings, according to which the ratios $s_n(f)/s_n(g)$ are eventually
bounded by a constant depending only on the common rotation number of both
maps. These bounds were  proved by \'Swi\c{a}tek (cf.~[Sw$_1$], [Sw$_2$]),
Herman [H$_2$] and Yoccoz (unpublished); the proofs assume
$C^3$-smoothness and a negative Schwarzian property near critical
points. An important corollary to Theorem A is the so-called {\it
golden mean universality conjecture}.
\proclaim{Corollary}.{If $f$ is a smooth critical circle mapping and 
$\rho(f) = {{\sqrt{5}-1}\over{2}}=[1,1,1,\ldots]$ (the golden mean),
then the scaling ratios of $f$ converge to a universal
constant.\hfill\endproof}

Computer-assisted work by Shenker [Sh] (cf.~[Ra$_2$]) shows that the value of
this universal constant is $0.7760513\ldots$ . We emphasize the golden mean
case here owing to its historical significance, yet a more general
corollary is the universality of scaling ratios for mappings whose rotation
number is a quadratic algebraic number, i.e. has an eventually periodic
continued fraction development.

A proper formulation of the second version involves the notion of
quasi-symmetry. If $h$ is a homeomorphism of the line (or circle, with its
linear coordinate), we let the {\it quasi-symmetric distortion} of $h$ at
scale $t>0$ be the number 
$$
k(h,t)\;=\;\sup_{0<s\le t}\sup_{x}\; {{h(x+s)-h(x)}\over{h(x)-h(x-s)}}\ .
$$
If $k(h,t)\le k<\infty$ for all $t>0$ then $h$ is a {\it
$k$-quasisymmetric} mapping. If moreover $k(h,t)\to 1$ as $t\to 0$,
then $h$ is said to be {\it symmetric}.
\proclaim{Theorem B}.{Any two smooth critical circle mappings with
the same rotation number of bounded combinatorial type are
conjugate by a symmetric homeomorphism.\hfill\endproof}

Both theorems will follow from certain {\it
renormalization convergence results}. Just as in the unimodal case,
one can define a renormalization scheme for critical circle mappings,
thanks to the fundamental notion of {\it commuting pair} developed by
Lanford and Rand (cf.~[L$_1$], [Ra$_1$]). Commuting pairs represent
whole conjugacy classes of circle mappings, and so each of them has a
rotation number of its own. The first return map to
$J_n$, consisting of $f^{q_n}$ restricted to $I_{n+1}$ and $f^{q_{n+1}}$
restricted to $I_n$, is the principal example of a commuting
pair, in this case the {\it $n$-th renormalization} of $f$
(see section I). This renormalization scheme acts as the {\it Gauss
map} on rotation numbers. As first observed by Ostlund, Rand, Sethna
\& Siggia in [ORSS], and also by Feigenbaum, Kadanoff \& Shenker in
[FKS], renormalization can be viewed as an operator
acting on an infinite-dimensional space of commuting pairs. In both
works, the same claim was made that a hyperbolic fixed-point for this 
renormalization operator exists, corresponding to an analytic
critical circle mapping with golden-mean rotation number. This claim
can be generalized in an obvious way to cover all rotation numbers
that are periodic under the Gauss map.
A computer assisted proof of the existence and hyperbolicity of
a golden-mean fixed-point, along the lines of
Lanford's proof for the Feigenbaum case, was given by Mestel in [Me].
Later, Epstein and Eckmann proved the existence without essential
help from the computer [EE]. Their proof uses Schauder's
theorem, and therefore guarantees neither uniqueness nor hyperbolicity
of the fixed-point. 
Taking a broader perspective, and inspired by his own
computer-assisted work on unimodal mappings, Lanford conjectured that
the renormalization operator is {\it globally hyperbolic} and
possesses an infinite-dimensional {\it horseshoe-like attractor}.

Although in this paper we don't go so far as proving Lanford's conjectures in
full, we do prove the existence and {\it global} uniqueness of the
golden-mean fixed-point, as well as of all other fixed or
periodic points of the renormalization operator, and describe their
(codimension-one) stable sets. 
We prove that the successive renormalizations of any two
commuting pairs representing critical circle mappings with the same
rotation number of bounded type converge together in the
$C^0$-topology.  Indeed, a stronger form of convergence, implying
$C^k$-convergence for all $k<\infty$, takes place if both pairs are
real-analytic of a special kind (see section IX). Our methods don't give any
rate of convergence, however, which is unfortunate since an exponential rate
would yield the Conjecture in the cubic case.

Our approach is based on the deep
holomorphic and quasiconformal ideas of Sullivan presented in
[S$_1$], [S$_5$], and detailed in [MS,~Ch.~VI].  
Here is a brief outline of the paper. In section I we define a
special class of real-analytic commuting pairs, the {\it Epstein
class}, which contains all limits of renormalization due to the real
a-priori bounds.
In section II, we introduce certain complex-analytic dynamical systems
called {\it holomorphic commuting pairs}. These objects restrict to
real-analytic commuting pairs on the line, and resemble
quadratic-like mappings in many ways. For instance, they have
annular fundamental domains and Julia sets, just as quadratic-like
mappings do. 
A holomorphic commuting pair can be renormalized, and the result is
again an object of the same type.
In section III, we prove a pull-back theorem for holomorphic
commuting pairs. This permits us to assign a quasiconformal distance
between topologically equivalent objects of this type. In the
resulting metric spaces, any two points can be joined by special paths,
whose elements are quasiconformal deformations of the endpoints,
called {\it Beltrami paths} (cf.~section V). Renormalization carries
Beltrami paths to Beltrami paths. We say that a Beltrami path is {\it
efficient} if the distance between its endpoints is not much smaller
than its length.
In section VI, we show how to factor the long compositions
representing high renormalizations of commuting pairs in the Epstein
class so that the factors satisfy the hypotheses of Sullivan's {\it
sector theorem}. This is the point where we have to assume that the
rotation number is of bounded combinatorial type. The factoring
combined with Sullivan's {\it sector inequality} proves, as stated in
section VII, that any sufficiently high renormalization of a
commuting pair in the Epstein class can be extended to a holomorphic
commuting pair, whose fundamental domain is a definitely thick
annulus. In particular, renormalizing a very long but efficient
Beltrami path sufficiently many times, we see using the pull-back
theorem that its endpoints are brought within a fixed distance. 
It is a beautiful discovery by Sullivan that in this situation the
image Beltrami path necessarily {\it coils}, i.e. it cannot be
efficient. Therefore the distances between points along the path are
contracted, and this implies strong renormalization convergence, as
we show in section IX.

This outline overlooks several important points. Thus, since the
boundaries of domain and range of holomorphic commuting pairs are
quite arbitrary, it is necessary to work with the {\it germs} of such
objects around their Julia sets, with a germ version of the
qc-distance called the {\it Julia-Teichm\"uller} distance, and also
with the infinitesimal form of that distance. Using examples of
holomorphic commuting pairs built from the complexified Arnold
family, we show in section IV that the Julia sets of these objects
carry no invariant line fields. Therefore all quasiconformal
deformations of a fixed germ are supported in the {\it external
class} of that germ, which is the Cantor repeller constructed in section
VIII. The space of backward-orbits of this Cantor repeller is a compact
{\it Riemann-surface lamination} in the sense of Sullivan. This gives
us a space to which Sullivan's coiling idea can be applied. The
qc-structures that are invariant for the repeller can be lifted to
the lamination. By Sullivan's {\it almost geodesic principle}, if
a structure of this kind comes from an optimal qc-conjugacy between two
germs, then it can be used to generate a very long but efficient
Beltrami path of structures on the lamination. The coiling lemma used
to prove contraction is a partial converse to this fact.

The results in this paper have a number of interesting applications. 
The basic theory of holomorphic pairs introduced here has been used
recently by McMullen [McM] in his elegant study of self-similarity
properties of Siegel disks.
We mention one further application, connected with the scalings of
frequency-locking intervals of one-parameter families of circle
mappings. In the Arnold family (1), the values of $\theta$ for which
the corresponding map has irrational rotation number form a Cantor
set. The gaps of this Cantor set have been examined numerically by
Cvitanovic \& S\"oderberg in [CS], and its Hausdorff dimension
estimated at about 0.87. \'Swi\c{a}tek gave a rigorous proof that this
Cantor set has zero Lebesgue measure in [Sw$_1$]. Later, in [GrS], he
and Graczyk proved that the Hausdorff dimension is less than 1 but
not smaller than $1\over 3$. Our results can be combined with a
careful analysis of the unstable manifolds of the renormalization
operator to establish the universality of the Hausdorff dimension
among cubic families. The analysis will be carried out in a  
forthcoming paper.

\smallskip
\ni
{\bf Acknowledgements.\enspace} I am grateful to my thesis advisor,
D.~Sullivan, for many beautiful lectures and insights. I wish to express
my thanks to W.~de~Melo, C.~Tresser, M.~Lyubich, O.~Lanford and
C.~McMullen for various conversations about renormalization, to
F.~Gardiner for teaching me Teichm\"uller theory, and to J.~Milnor for
his comments and kind suggestions upon reading the manuscript. I am
indebted also to H.~Epstein for showing me a computer picture of the
golden mean fixed-point that inspired the definition of holomorphic
commuting pair.  
\bigskip

\medskip
\ni\centerline{\smc I.\enspace Renormalization of real commuting pairs
and the Epstein class} 
\smallskip
\ni
Let $f : S^1 \to S^1$ be a smooth, orientation-preserving
homeomorphism having exactly one critical point $c\in S^1$ of cubic
type. That is, let $f$ be such that we can represent it in the
form $h\circ f_\theta\circ H$, where $h$ and $H$ are smooth,
orientation-preserving diffeomorphisms, $\theta$ is a real number and
$f_\theta$ is the mapping whose lift $E_\theta$ to the real line is
given by
$$
E_\theta (x) = x+\theta -{1\over {2\pi}}\sin{2\pi x}
\ .
$$
We call $f$ a {\it critical circle mapping}. We also refer to $f=h\circ
f_\theta\circ H$ as an $hQ$-decomposition of $f$ (cf. [S$_1$]).
Our standing assumption in this paper is that $f$ has no periodic
points, i.e. that its rotation number $\rho (f)$ is irrational. Thus
$f$ is topologically conjugate to the corresponding irrational
rotation, after a well-known theorem of Yoccoz (cf.~[Yo$_2$ ]). We write
the rotation number of $f$ as an infinite continued fraction $\rho(f)
=  [r_0, r_1, \ldots, r_n, \ldots]$ and let $(q_n)_{n\ge 0}$ denote
the successive closest return times given recursively by $q_0 = 1$,
$q_1 = r_0$ and
$$
q_{n+1} = r_nq_n + q_{n-1} \ , \eqno{(1)}
$$
for all $n\ge 1$. Recall
that each of these numbers appears as denominator in the truncated 
expansion of order $n$ of $\rho (f)$ in its irreducible form 
$$
{p_n\over{q_n}}\;=\;[r_0,r_1,\ldots ,r_{n-1}] \ \ , 
$$
where $p_0=0$, $p_1=1$ and $p_{n+1}=r_np_n +p_{n-1}$ for all $n\ge 1$. It is
also convenient to set $q_{-1}=0$.
We denote by $I_n(c)$ the closed interval in $S^1$ with
endpoints $c$ and $f^{q_n}(c)$ containing $f^{q_{n+2}}(c)$. The
dynamical first return map to the interval $I_n(c) \cup I_{n+1}(c)$
is given by $f^{q_{n+1}}$ on $I_n(c)$ and by $f^{q_n}$ on
$I_{n+1}(c)$. Each pair $(f^{q_n}, f^{q_{n+1}})$ yields an example of
what one calls {\it weakly commuting pair} or simply {\it
commuting pair}, after Lanford and Rand (cf.~[L$_1$], [L$_2$], [Ra$_1$],
[Ra$_2$]). Here is the abstract definition. 

\ni {\it Definition 1.\enspace} A commuting pair $\zeta = (\xi, \eta)$
consists of two orientation preserving smooth homeomorphisms
$\xi : I_\xi \to \xi(I_\xi), \eta : I_\eta \to \eta(I_\eta)$ into the
reals where

\itemitem{(a)} $I_\xi = [\eta(0), 0] \subseteq \real $, $I_\eta = [0, \xi(0)] 
\subseteq \real $; 

\itemitem{(b)} Both $\xi$ and $\eta$ have homeomorphic extensions,
with the same degree of smoothness, to interval neighborhoods of their
corresponding domains, and such extensions commute, i.e. $\xi\circ\eta
= \eta\circ\xi$, wherever both sides are defined; 

\itemitem{(c)} $\xi\circ\eta(0)$ belongs to $I_\eta$; 

\itemitem{(d)} We have $\xi' (x) \ne 0 \ne
\eta'(y)$, for all $x$ in $I_\xi \setminus \{0\}$ and  all $y$ in
$I_\eta \setminus \{0\}$.

\ni
A {\it critical commuting pair} is a commuting pair that has
$hQ$-decompositions $\xi = h_\xi \circ Q\circ H_\xi$ and $\eta = h_\eta
\circ Q\circ H_\eta$ where $h_\xi$, $h_\eta$, $H_\xi$, $H_\eta$
are smooth diffeomorphisms and $Q$ is the map $z\mapsto z^3$. 
\smallskip
An object which is either a commuting pair or obtained from a
commuting pair by conjugating $\xi$ and $\eta$ by $x\mapsto -x$
(resp. by $x\mapsto \lambda x$, $\lambda \ne 0$) is called a
{\it commuting pair up to orientation} (resp. {\it up to linear
rescaling}). 
A critical circle mapping $f$ gives rise to a sequence of critical
commuting pairs in the following way. Let $\overline f$ be a lift of
$f$ to the real line satisfying ${\overline f}'(0)=0$ and $0<{\overline
f}(0)<1 $. For each $n\ge 0$, let $J_n\subseteq \real$ be the closed
interval adjacent to zero that projects down homeomorphically onto
$I_n(c)$ via the exponential mapping. Let $T:\real \to \real$ be the
translation $x\mapsto x+1$ and let $\xi_{f,n}:J_{n+1}\to \real$ be
given by $\xi_{f,n}(x)=T^{-p_n}\circ {\overline f}^{q_n}(x)$;
similarly, let $\eta_{f,n}:J_{n}\to \real$ be given by
$\eta_{f,n}(x)=T^{-p_{n+1}}\circ {\overline f}^{q_{n+1}}(x)$ . Then
$\zeta_{f,n} = (\xi_{f,n}, \eta_{f,n})$ is a critical commuting
pair up to orientation.

Conversely, regarding $I=[\eta(0),\xi(0)]$ as the circle (identifying
$\eta(0)$ and $\xi (0)$) and letting $f_\zeta :I \to I$ be given
by 
$$
f_\zeta (x) = \cases{\xi (x), &if $\ \eta(0)\le x < 0$ \cr
{} & {} \cr
\eta (x), &if $\ 0 \le x \le \xi(0)$ \cr} \ , \eqno{(3)}
$$
we recover a plethora of critical circle mappings from a critical commuting
pair $\zeta =(\xi ,\eta )$.
We perform the glueing of $\eta(0)$ to $\xi (0)$ via the mapping 
$\eta\xi^{-1}$, which by conditions (b) and (d) above maps a small
neighborhood of $\xi(0)$ diffeomorphically onto a small neighborhood
of $\eta(0)$. We obtain a smooth, closed one-manifold $M$ as the
quotient space, and $f_\zeta$ projects down to a smooth homeomorphism
$F_\zeta:M\to M$. Each identifying diffeomorphism $\varphi :M \to S^1$
gives rise to a critical circle mapping $f^\varphi=\varphi\circ
F_\zeta\circ \varphi^{-1}$.
Although there is no canonical
choice for $\varphi$, any two choices are such that the corresponding
$f^\varphi$ 's differ by a diffeomorphism. Therefore we recover a
whole {\it smooth conjugacy class} of critical circle mappings
(see [dFM$_1$] for a detailed exposition of the glueing procedure,
first introduced by Lanford).  
We will abuse language henceforth and call $f_\zeta$
{\it the} critical circle mapping of $\zeta$. 
We let $I_n\subseteq I$ be the closed interval that
corresponds to $I_n(c)$ for any representative $f^\varphi$, for each
$n\ge 0$. The endpoints of $I_n$ are $0$ and $f_\zeta^{q_n}(0)$, where
$\{q_n\}$ is the sequence of return times of any such representative. 

Letting $\rho(\zeta) = \rho(f_\zeta)$ be the rotation number of
$\zeta$, we are ready to define the renormalization operator for
commuting pairs. If $\rho(\zeta) = [r+1, r_1, r_2, \ldots, r_n, 
\ldots]$, then $\eta^{r+1}\xi (0)<0<\eta^r\xi (0)$ and one verifies that the
mappings $\eta |[0,\eta^r\xi (0)]$ and $\eta^r\circ \xi|I_\xi$ constitute a
commuting pair up to orientation.
\smallskip
\ni {\it Definition 2.\enspace} The commuting pair 
$$
{\cal R} \zeta = (\eta^{\#},
(\eta^r \circ \xi)^{\#}) \ ,
$$ 
where $^{\#}$ denotes linear rescaling by the factor $\lambda =
\xi(0)/\eta(0) < 0$, is called the {\it first renormalization of
$\zeta$.} We also refer to $(\eta ,\eta^r\circ\xi)$ as the first
renormalization of $\zeta$ without rescaling.
\smallskip
Thus, in the notation introduced above, we have ${\zeta_{f,n}}^{\#} = 
{\cal R}({\zeta_{f,n-1}}^{\#})$ for all $n \ge 1$. These may
therefore be regarded as the successive renormalizations of $f$.
It is easy to see that $\rho({\cal R}\zeta) = [r_1+1, r_2, \ldots, r_n,
\ldots]$. Thus, renormalization acts essentially as the Gauss map on
rotation numbers. 

Now we define a class of commuting pairs containing the {\it attractor} of
renormalization.
\smallskip 
\ni {\it Definition 3.\enspace} A real-analytic commuting pair $\zeta = (\xi, \eta)$
is said to be in the Epstein class ${\cal E}$ if, for $\gamma = 
\xi, \eta$, there exists a decomposition $\gamma = h_\gamma \circ
Q$, where as before $Q:z\mapsto z^3$, such that

\itemitem{(a)} $h_\gamma : Q(I_\gamma) \to \gamma(I_\gamma)$ is an
orientation preserving diffeomorphism;  

\itemitem{(b)} $h_\gamma^{-1}$ extends to a schlicht mapping ${\bf
C}(\widetilde I_\gamma) \to {\bf C}$, where $\widetilde I_\gamma
\supseteq \gamma(I_\gamma)$ is some open interval (here ${\bf C}(I)\,=\,{\bf
C}\setminus (\real\setminus I)$).
\smallskip 

For each $s>0$, let us write, taking into account condition (b) above,
$$
{\cal E}_s\,=\,\{ \zeta \in {\cal E}: {\widetilde I}_\gamma \supseteq
I_\gamma^s \ , \gamma = \xi ,\eta \} \ ,
$$
where $I^s$ denotes the interval centered at the midpoint of $I$ whose length
is $(1+s)$-times the length of $I$ (cf. section VI). We also refer to each
${\cal E}_s$ as an Epstein class.

Now we have the following fact.
\proclaim{Lemma I.1}.{The Epstein class $\cal E$ is invariant under
renormalization. \hfill\endproof }  

The following lemma describes how the successive renormalizations
without rescaling of a commuting pair $\zeta$ are nested inside $\zeta$. It
will be also extremely useful in section VI, in the breaking-up of long
renormalization compositions leading to the complex bounds.

\proclaim{Lemma I.2}.{Let $\zeta =(\xi ,\eta )$ be a critical
commuting pair and let $(\xi_n, \eta_n)$ be the sequence
of renormalizations of $\zeta$ without rescaling. Then $I_{\xi_n}=I_n$
and $I_{\eta_n}=I_{n-1}$ for all $n\ge 1$ and we have the 
following hybrid representations
$$
\cases{n\quad\hbox{even}\quad \Rightarrow \cases{\xi_n(x) =
f_\zeta^{q_{n-1}-1}\circ \xi (x)&for all $x$ in $I_{\xi_n}=I_n$ \cr
{} & {} \cr
\eta_n(x) = f_\zeta^{q_n -1}\circ \eta(x) &for all $x$ in
$I_{\eta_n}=I_{n-1}$ \cr} \cr 
{} & {} \cr 
{} & {}\cr 
n\quad\hbox{odd}\quad\ \Rightarrow \cases{\xi_n(x) = 
f_\zeta^{q_{n-1}-1} \circ \eta (x) &for all $x$ in $I_{\xi_n}=I_n$ \cr 
{} & {} \cr
\eta_n(x) = f_\zeta^{q_n - 1} \circ \xi (x) &for all $x$ in
$I_{\eta_n}=I_{n-1}$ \ .\cr}  
\cr}
\eqno{(4)}
$$
}
\proof The first assertion follows easily by induction on $n$, using
the recurrence relations (1). The hybrid expressions in (4) are clear
if we observe that $I_n$ is contained in the domain of $\xi$ when $n$
is even, and in the domain of $\eta$ when $n$ is odd.
\hfill\endproof

Because it relates the dynamics of $\zeta$ with that of $f_\zeta$,
Lemma I.2 can be used to transfer certain well-known a-priori bounds
for critical circle mappings to corresponding ones for critical
commuting pairs. Let $s_n(\zeta)=|I_{n+1}|/|I_n|$ be the $n$-th
scaling ratio of $\zeta$. Also, if $\zeta_1$ and $\zeta_2$ have the
same rotation number, let their {\it quasi-symmetric} distance be the
number 
$$
d_{QS}(\zeta_1,\zeta_2)=\log{k(h)} \ ,
$$
where $h:[\eta_1(0),\xi_1(0)]\to [\eta_2(0),\xi_2(0)]$ is the
conjugacy between both pairs and $k(h)$ is the quasi-symmetric
distortion of $h$.
Then, we can use Lemma I.2 to re-state the well-known results of
Herman [H$_2$], Yoccoz (unpublished), \'Swi\c{a}tek [Sw$_2$] and Graczyk \&
\'Swi\c{a}tek [GrS] in the following combined form.  

\proclaim{Theorem I.3}.{Given $0<\alpha <1$ irrational, there exist
constants $K_1>1$ and  $K_2,K_3>0$ depending only on 
$\alpha$ such that the following statements hold. 
\itemitem{(a)} If $\zeta$ is a critical commuting pair with rotation
number $\alpha$, then for all sufficiently large $n$ we have
$K_1^{-1}|I_n|\le |I_{n+1}| \le K_1|I_n|$; in other words the scaling
ratios of $\zeta$ are bounded away from zero and infinity;
\itemitem{(b)}If $\zeta_1$ and $\zeta_2$
are critical commuting pairs with rotation number $\alpha$, then for
all sufficiently large $n$ we have 
$$ 
|\d{s_n(\zeta_1)\over{s_n(\zeta_2)}}-1|\le K_2
\ ,
$$
and moreover $d_{QS}({\cal R}^{n}\zeta_1, {\cal R}^{n}\zeta_2)\le K_3$.
\hfill\endproof}
Notice in particular that any two critical commuting pairs with the
same irrational rotation number are quasi-symmetrically conjugate.
Using Theorem I.3 and the bounded geometry results and
techniques of Sullivan [S$_1$, \S 4], one obtains the following
fundamental compactness result, which is essentially the {\it pure
singularity property} of \'Swi\c{a}tek [Sw$_2$]. A complete, detailed proof
of this theorem (and much more) can be found in [dFM$_1$]. 

\sproclaim{Theorem I.4}.{Let $\zeta$ be a critical commuting pair
of class $C^r$ ($r\ge 3$) with irrational rotation number $\rho(\zeta)$,
and consider the $hQ$-decompositions of its successive renormalizations
$\xi^{\#}_n = h_{\xi,n}\circ Q\circ H_{\xi ,n}$ and
$\eta^{\#}_n = h_{\eta ,n} \circ Q\circ H_{\eta ,n}$.
Then the families $\{\xi^{\#}_n\}_{n\ge 0}$ and $\{\eta^{\#}_n\}_{n\ge
0}$ are precompact in the sense that, for $\gamma = \xi, \eta$, the
following conditions hold.
\itemitem{(a)} The critical values of $\gamma^{\#}_n$ are bounded
away from zero;
\itemitem{(b)} There exist $s>0$ depending only on the rotation
number of $\zeta$, fixed intervals ${\cal I}_\gamma$, ${\cal J}_\gamma$
and a positive integer $N$ such that, for all $n\ge N$,
$h_{\gamma ,n}^{-1}$ is well-defined on ${\cal I}_\gamma$ and ${\cal
I}_\gamma \supseteq (\gamma_n^{\#}I^{\#}_{\gamma_n})^s$, and $H_{\gamma ,n}$
is well-defined on ${\cal J}_\gamma$ and ${\cal J}_\gamma \supseteq
(\gamma_n^{\#})^{-1}({\cal I_\gamma})$;
\itemitem{(c)} The family $\{h_{\gamma ,n}^{-1}|{\cal I}_\gamma
\}_{n\ge N}$ has compact closure in the $C^r$-topology on
diffeomorphisms;
\itemitem{(d)} The sequence $(H_{\gamma ,n}|{\cal
J}_\gamma)_{n\ge N}$ converges to the identity in the
$C^r$-topology on diffeomorphisms.}.
{\ni Moreover, every $C^r$-limit of $\zeta^{\#}_n = (\xi^{\#}_n,
\eta^{\#}_n)$ is a critical commuting pair in the 
Epstein class ${\cal E}_s$. \hfill\endproof}

\ni {\it Remark 1.\enspace} The facts stated in Theorems I.3 and I.4 are
collectively known as the {\it real a-priori bounds} for critical
circle mappings.

\ni {\it Remark 2.\enspace} As it turns out, proving the asymptotic rigidity
statements of Theorems A and B in the Introduction is tantamount to showing
that $d_{QS}({\cal R}^{n}\zeta_1, {\cal R}^{n}\zeta_2) \to 0$ as $n\to 
\infty$. An exponential rate of convergence would yield the Conjecture in the
cubic case, see [dFM$_1$], [dFM$_2$].

\bigskip

\ni\centerline{\smc II. \enspace Holomorphic commuting pairs }
\smallskip
\ni
Now we introduce special complex-analytic extensions of critical commuting
pairs, akin to quadratic-like mappings. We need an auxiliary definition.  
Let us say that a configuration of simply connected domains $({\cal D}, {\cal
O}_\xi, {\cal O}_\eta, {\cal O}_\nu)$ in the plane is a {\it bow-tie} if the
following conditions are satisfied: (a) all four domains are symmetric about
the real axis; (b) each ${\cal O}_\gamma $ is a Jordan domain
whose closure is contained in ${\cal D}$; (c) $\overline {\cal O}_\xi \cap
\overline {\cal O}_\eta =    
\{0\} \subseteq {\cal O}_\nu$; (d) the differences ${\cal O}_\gamma
\setminus  {\cal O}_\nu$ and ${\cal O}_\nu \setminus  {\cal O}_\gamma
$ are non-empty connected sets for $\gamma = \xi, \eta$; (e) the interval
${\cal O}_\xi \cap \real$ lies to the left of zero.    
Let $J_\gamma $ denote the open intervals ${\cal O}_\gamma  \cap
\real$ for $\gamma= \xi, \eta, \nu$. Then $J_\xi$ and $J_\eta$ share
an endpoint at the origin, and $J_\xi$ lies in the negative real axis. Also,
$J_\nu$ contains the origin and is contained in $\overline J_\xi \cup J_\eta$.
Moreover, due to condition (d) we know  
that ${\cal O}_\xi \cup {\cal O}_\eta \cup {\cal O}_\nu$, as well as ${\cal
O}_\gamma \cap {\cal O}_\nu$ ($\gamma = \xi, \eta$) are Jordan domains.  
A sketch of the situation we have in mind is shown in Figure 1.  


\bigskip
$$
\psannotate{
\psboxto(11cm;0cm){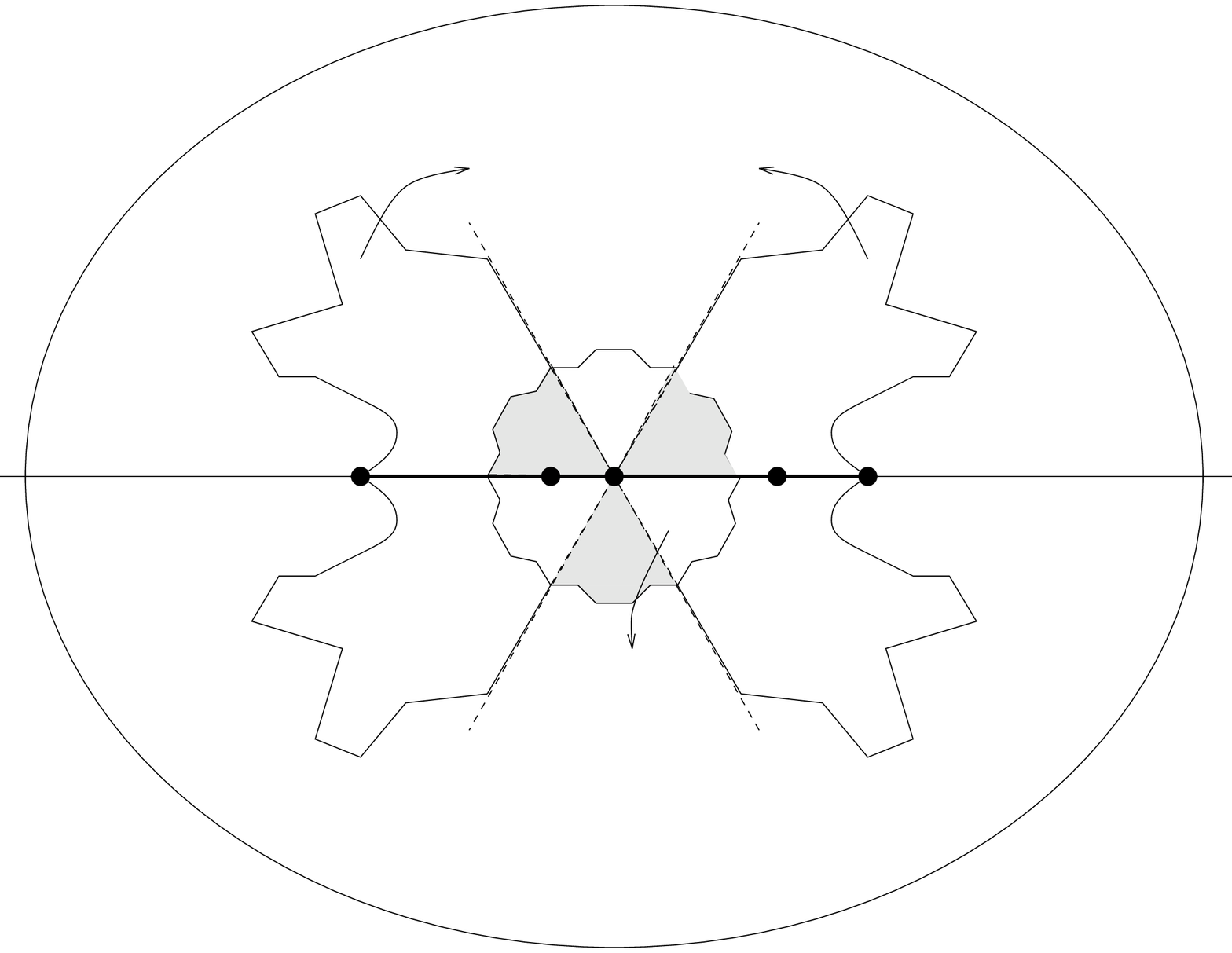}}
{\fillinggrid\at{7\pscm}{10\pscm}{${\cal O}_\xi$}
             \at{20.6\pscm}{10\pscm}{${\cal O}_\eta$}
             \at{10.5\pscm}{14\pscm}{$\xi$}
             \at{16.7\pscm}{14\pscm}{$\eta$}
             \at{20\pscm}{15\pscm}{${\cal V}$}
             \at{12.6\pscm}{11\pscm}{${\cal O}_\nu$}
             \at{12.7\pscm}{9\pscm}{$0$}
             \at{10.8\pscm}{7.5\pscm}{$\eta(\!0\!)$}
             \at{14.8\pscm}{7.5\pscm}{$\xi(\!0\!)$}
             \at{7.4\pscm}{7.5\pscm}{$a$}
             \at{16.2\pscm}{7.5\pscm}{$b$}
             \at{12.2\pscm}{5\pscm}{$\nu$}
}
$$
\centerline{\it Figure 1}
\bigskip

\ni {\it Definition 4.\enspace} A {\it holomorphic commuting pair} consists 
of a bow-tie $({\cal D}$, ${\cal O}_\xi$, ${\cal O}_\eta$, ${\cal O}_\nu$) together  
with complex analytic mappings $\xi, \eta, \nu$ having ${\cal O}_\xi, {\cal
O}_\eta, {\cal O}_\nu$ respectively as their domains and a positive integer
$m$ satisfying the following conditions: 
 
\itemitem{[H$_1$]} All three mappings commute with complex conjugation.

\itemitem{[H$_2$]} $\xi$ and $\eta$ are schlicht mappings onto ${\cal D} \cap {\bf
C}(\xi(J_\xi))$ and ${\cal D} \cap {\bf C}(\eta(J_\eta))$ respectively.

\itemitem{[H$_3$]} $\nu$ is a 3-fold branched covering onto ${\cal D} \cap {\bf
C}(\nu(J_\nu))$, with a unique critical point at zero. 

\itemitem{[H$_4$]} $\xi$ and $\eta$ have analytic extensions to a  
certain neighborhood of zero where both $\xi\circ\eta$ and $\eta 
\circ\xi$ are defined, and we have $\xi\circ\eta(z) = \eta \circ 
\xi(z) = \nu(z)$ for all $z$ in that neighborhood.

\itemitem{[H$_5$]} If $x \in J_\xi$ then $\xi(x) > x$, whereas if $x \in J_\eta$ 
then $\eta(x) < x$; moreover, $\xi(0), \nu(0) \in J_\eta$ and  
$\eta(0) \in J_\xi$.

\itemitem{[H$_6$]} If $a$ is the left endpoint of $J_\xi$ and $b$ is the right endpoint 
of $J_\eta$ then both $\xi^m(a)$ and $\eta(b)$ are well-defined as boundary
values, and we have $\xi^m(a)=\eta(0)$ and $\eta(b) = \xi(0)$. 
\smallskip 

Holomorphic commuting pairs will be denoted by $\Gamma$. Explicit examples
will be constructed in section IV. If $\Gamma$ is given, it is clear from
H$_2$, H$_3$, H$_4$ and H$_5$ that the restrictions 
$\xi|[\eta(0), 0]$ and $\eta|[0,\xi(0)]$ constitute a real analytic critical
commuting pair.  
Thus, we define the {\it rotation number} 
of $\Gamma$ to be the rotation number of its real commuting pair.  
Note that the intervals $J = \overline J_\xi \cup J_\eta$ and
$I=[\eta(0),\xi(0)]$ are both forward invariant
under the dynamical system generated by $\xi$ and $\eta$.
We call them the {\it large} and  {\it small} dynamical intervals of $\Gamma$,
respectively. We also say that the integer $m$ in condition H$_6$ is the {\it
height} of $\Gamma$. The following proposition is fundamental. 
 
\proclaim{Proposition II.1}.{In any holomorphic commuting pair, the  
mappings $\xi$ and $\eta$ have analytic extensions to ${\cal O}_\xi \cup {\cal O}_\nu$ 
and ${\cal O}_\eta \cup {\cal O}_\nu$ respectively. Moreover, the restrictions  
$\xi_* = \xi|{\cal O}_\nu$ and $\eta_* = \eta|{\cal O}_\nu$ are 3-fold branched  
covering maps onto ${\cal O}_\eta$ and ${\cal O}_\xi \cap {\bf C}([\xi^{-1}\circ
\eta(0), 0])$ respectively, and we have $\eta\circ\xi_* = \xi\circ\eta_* =
\nu$.  }
 
\proof
We use the 3-fold symmetry of ${\cal O}_\nu$ coming from $\nu$ in order to
extend $\xi$ and $\eta$ by Schwarz reflection in the following way.

By H$_2$, the composition $\eta\circ\xi$ is a well-defined schlicht mapping of
$V = \xi^{-1}({\cal O}_\eta)$ onto ${\cal D}\cap {\bf C}([\eta(0), \eta\xi(0)])$. 
Let $Y = \nu^{-1}([\nu(0), +\infty))$. Then, using conditions H$_1$ and H$_3$
we readily see that ${\cal O}_\nu \setminus  Y$ has exactly 3 connected
components, one of which, call it $W$, is symmetric about the real axis. 
We claim that $V = W$. Since $V$ is also symmetric about the real
axis, it is enough to show that $V^+ = W^+$. Now,  $\nu$ maps $W^+$
injectively onto ${\cal D}^+$; likewise, $\eta\circ   
\xi$ maps $V^+$ onto ${\cal D}^+$ injectively. Hence the composition $\phi 
= \nu^{-1}\circ(\eta\circ\xi)$ is well-defined in $V^+$ and maps it onto 
$W^+$. Since by H$_4$ we have $\eta \circ\xi \equiv \nu$ 
on some neighborhood ${\cal O}$ of zero, we deduce that $\phi(z) = z$ for  
all $z \in {\cal O} \cap V^+$. Therefore $\phi$ must be the identity
map, which settles the claim.   

From this, it follows that $\xi^{-1}(0)$ is the left endpoint of
$J_\nu$, and since $\nu$ agrees with $\eta\circ\xi$ over all of $W$,
we see that $\nu(\xi^{-1}(0)) = \eta(0)$.  
Switching the roles of $\xi$ and $\eta$ in this argument, we deduce that
$\eta^{-1}(0)$ is the right endpoint of $J_\nu$ and that $\nu(\eta^{-1}(0)) =
\xi(0)$. Therefore by H$_3$ the image of ${\cal O}_\nu$ under
$\nu$ is ${\cal D} \cap {\bf C}([\eta(0), \xi(0)])$, which by H$_2$ and
the last equality in H$_6$ is the image of ${\cal
O}_\eta$ under $\eta$.  
This shows that $\xi_* = \eta^{-1}\circ\nu : {\cal O}_\nu 
\to {\cal O}_\eta$ is well-defined. It is clearly  
a 3-fold branched covering onto ${\cal O}_\eta$, and since $\nu$ agrees with  
$\eta\circ\xi$ over $W$, we have $\xi_* \equiv \xi$ there. The proof for
$\eta$ is similar.
\hfill\endproof

Our next proposition introduces a renormalization operator for holomorphic 
commuting pairs which is compatible with the real 
renormalization operator of section I. In the proof we shall
use the following elementary set-theoretic remark.  
 
\proclaim{Lemma II.2}.{Let $\phi : A \to B$ be one-to-one and onto, and let 
$(B_n)_{n\ge 0}$ be the sequence of subsets of $B$ defined by  
$B_0 = B$ and $B_{n+1} = \phi(A \cap B_n)$. Then for each $n\ge
1$ the $n$-th iterate $\phi^n$ is well-defined over $A_{n-1} = A \setminus   
\bigcup_{i=0}^{n-1}\phi^{-i}(B\setminus A)$, and maps $A_{n-1}$ bijectively 
onto $B_{n-1}$. \hfill\endproof }

\proclaim{Proposition II.3}.{Let $\Gamma$ be a holomorphic commuting
pair. Then there exists a holomorphic commuting pair 
${\cal R}(\Gamma)$ whose underlying real commuting pair is the first  
renormalization of the real commuting pair of $\Gamma$.}
 
\proof
Let $\rho(\Gamma) = [r, r_1, r_2, \cdots, r_n, \cdots]$ be the
rotation number of $\Gamma$. Recall that the first renormalization 
of $(\xi, \eta)$ is the pair $(\eta, \eta^r \circ \xi)$ up to the linear 
rescaling given by $x \mapsto \lambda x$, where $\lambda = \xi(0)/\eta(0)  
< 0$.
We will obtain the desired ${\cal R}(\Gamma)$ up to such rescaling by
constructing domains ${\cal O}_{\widehat\xi}, {\cal
O}_{\widehat\eta}, {\cal O}_{\widehat\nu}$ and corresponding maps
$\widehat\xi, \widehat\eta, \widehat\nu$. From $\rho(\Gamma)=
(r,r_1,\cdots ,r_n,\cdots )$, we also know that
$$
\eta^r\xi(0) > 0 > \eta^{r+1}\xi(0) \ ,
$$
and therefore $\eta(0) < \eta^{k-1}\xi(0)$
provided $1 \le k \le r+1$.    
 
  
First we take ${\cal O}_{\widehat\xi} = {\cal O}_\eta$ and set $\widehat\xi = \eta$. 
Let $A= {\cal D} \cap {\bf C}([\eta(0), \xi(0)])$ be the image of 
$B= {\cal O}_\eta $ under $\eta$ and set $\phi = \eta^{-1} : A \to B$. For
each $k \ge 1$ we know by Lemma II.2 that $\phi^k$ is a well-defined one-to-one
map of
$$
A_{k-1} = A \setminus   
\bigcup_{i=0}^{k=1}\phi^{-i}(B\setminus A) = {\cal D} \cap {\bf C}([\eta(0), \eta^{k-1} 
\xi(0)])
$$
onto $B_{k-1} \subseteq 
{\cal O}_\eta$.
It follows that, for $1\le k\le r+1$, the domain $A_{k-1}$ is symmetric and
simply-connected, and since $\phi^k$ is schlicht, the same holds for
$B_{k-1}$.     
One can verify that $A \cap B_{k-1}$ equals $B_{k-1}$ minus the slit
$[\xi(0), b)$, which opens-up when $\phi = \eta^{-1}$ is applied.
Since $B_k = \phi(A \cap B_{k-1})$ and $B_0 = {\cal O}_\eta$ is a Jordan
domain, an inductive argument shows that $B_{k-1}$ is a Jordan domain for
$1 \le k \le r+1$. 
In particular, $B_{r-1} \subseteq {\cal O}_\eta$ is a Jordan domain. 
Let $\xi_* : {\cal O}_\nu \to {\cal O}_\eta$ be the 3-fold branched covering given 
by Proposition II.1, put $U = \xi_*^{-1}(B_{r-1})$ and define 
$\widehat\eta_* : U \to A_{r-1}$ by $\widehat\eta_*  = \eta^r \circ \xi_*$. 
Since the critical value of $\xi_*$ belongs to $B_{r-1}$, 
we see that $U$ is a Jordan domain and that $\widehat\eta_*$ is a 
3-fold branched covering onto its image. We then take  
${\cal O}_{\widehat\eta} = U \cap {\cal O}_\xi$ and set $\widehat\eta =  
\widehat\eta_*|{\cal O}_{\widehat\eta}$. As $U \cap {\cal O}_\xi =  
\xi_*^{-1}(B_{r-1}) \cap {\cal O}_\xi$ is  
mapped bijectively by $\widehat\eta_*$ onto $B_{r-1} \setminus  [\xi(0), b)$, 
we deduce that ${\cal O}_{\widehat\eta}$ is a Jordan domain and that
$\widehat\eta$ is schlicht over ${\cal O}_{\widehat\eta}$.
On the other hand, $B_r \subseteq {\cal O}_\eta$ is also a Jordan domain 
containing $\xi(0)$, the critical value of $\xi_*$. Hence, if we 
let ${\cal O}_{\widehat\nu} = \xi_*^{-1}(B_r) \subseteq {\cal O}_\nu$ and  put
$\widehat\nu = \eta^{r+1} \circ \xi_* : {\cal O}_{\widehat\nu} \to A_r$, we see at once  
that ${\cal O}_{\widehat\nu}$ is a Jordan domain and that $\widehat\nu$ is a  
3-fold branched covering onto its image. Moreover, $\widehat\nu =   
\eta \circ (\eta^r \circ \xi_*) = \widehat\xi\circ\widehat\eta_*$. 
 
Now, if we linearly rescale ${\cal D}, {\cal O}_{\widehat\xi}, {\cal O}_{\widehat\eta}, 
{\cal O}_{\widehat\nu}, \widehat\xi, \widehat\eta, \widehat\nu$ by the map  
$z \to \lambda z$, we get the desired  
holomorphic commuting pair ${\cal R}(\Gamma)$. Indeed, $({\cal D},
{\cal O}_{\widehat\xi}, {\cal O}_{\widehat\eta}, {\cal O}_{\widehat\nu})$ 
is a bow-tie up to linear rescaling. Moreover,
conditions H$_1$ through H$_4$ hold by construction. 
Condition H$_5$ is not satisfied until we do the rescaling (which 
reverses orientation on the line), when it becomes clear. Finally,
since we have
$$\widehat\xi^{r+1}(b) = \eta^{r+1}(b) = \eta^r\xi(0) =
\widehat\eta(0)
$$
and
$$\widehat\eta(c) = \eta^r\circ\xi(c) =
\eta^r(\eta^{-r+1}(0)) = \eta(0) = \widehat\xi(0) \ ,
$$
where $c$ is the left endpoint of ${\cal O}_{\widehat\eta} \cap \real$,
condition H$_6$ is satisfied if we take $m = r+1$. 
\hfill\endproof
 
We are interested in the dynamical system generated by 
the mappings $\xi|{\cal O}_\xi, \eta|{\cal O}_\eta$ and $\nu|{\cal O}_\nu$.
We shall identify $\Gamma$ itself with this dynamical system.    
It will be very important to know that the $\Gamma$-orbits 
can be encoded by a single (discontinuous, piecewise holomorphic) 
transformation. This will be made precise in the proposition below. Let $F : ({\cal
O}_\xi \cup {\cal O}_\eta \cup {\cal O}_\nu) \to {\cal D}$ be given by  
$$
F(z) = \cases{\xi(z) &if $z \in {\cal O}_\xi$ \cr 
\eta(z) &if $z \in {\cal O}_\eta$ \cr 
\nu(z) &if $z \in {\cal O}_\nu \setminus  ({\cal O}_\xi \cup {\cal O}_\eta)$ \cr }
$$ 
We call $F$ the {\it shadow} of the holomorphic 
commuting pair $\Gamma$. 
  
\proclaim{Proposition II.4}.{Given a holomorphic commuting pair  
$\Gamma$, consider its shadow $F$ and let  
${\cal U} = {\cal O}_\xi \cup {\cal O}_\eta \cup {\cal O}_\nu$ and $X = J \cup F^{-1}(J)$, 
where $J$ is the large dynamical interval of $\Gamma$. Then
\itemitem{(a)} The restriction of $F$ to ${\cal U} \setminus  X$ is a regular 3-fold
covering mapping onto ${\cal D}^+ \cup {\cal D}^-$;
\itemitem{(b)} $F$ and $\Gamma$ share the same orbits
as sets.
}  
 
\proof
Since ${\cal U} \setminus  X$ consists of six 
connected components, each of which is mapped  
bijectively onto either ${\cal D}^+$ or ${\cal D}^-$, part (a) follows. 
In order to prove (b), it suffices to show that the $\Gamma$-orbit of any point 
of ${\cal U}$ is contained in the corresponding $F$-orbit. 
Thus, let $z \in {\cal U}$ and let $\omega$ be 
any finite admissible word in the alphabet $\{\xi, \eta, \nu\}$. If the letter 
$\nu$ does not occur in $\omega$ then $\omega(z) = F^{|\omega|}(z)$,
where $|\omega| = $ length of $\omega$. Otherwise we write $\omega = 
\omega_L\nu\omega_R$ for some other words $\omega_L, \omega_R$ in the 
same alphabet (possibly empty); setting $x = \omega_R(z)$, we have 
three possibilities
 
\itemitem{(i)} $x \in {\cal O}_\nu \setminus  ({\cal O}_\xi \cup {\cal O}_\eta)$: in this case $\nu(x) = 
F(x)$ by definition so we may replace $\nu$ by $F$ in $\omega$. 

\itemitem{(ii)} $x \in {\cal O}_\xi \cap {\cal O}_\nu$: here we may write, using
Proposition II.1, $\nu(x) = \eta\xi(x) = \eta F(x)$; since $\xi({\cal O}_\xi \cap
{\cal O}_\nu) \subseteq {\cal O}_\eta$ by that same proposition, we deduce that
$F(x)$ is in ${\cal O}_\eta$, and so $\eta F(x) = F \circ F(x) = F^2(x)$.
Hence in this case we may replace $\nu$ by $F^2$ in $\omega$. 

\itemitem{(iii)} $x \in {\cal O}_\eta \cap {\cal O}_\nu$: same as (ii). 

\ni This substitution process applied to all occurences of $\nu$ in $\omega$ 
shows that $\omega(z) = F^n(z)$ for some $n \ge |\omega|$, and 
so part (b) is proved also. \hfill\endproof

\bigskip

\ni\centerline{\smc III. \enspace The pull-back theorem }
\smallskip
\ni
The principal reason why holomorphic commuting pairs are useful is
the fact that any quasi-symmetric conjugacy between the restrictions
of two such objects to the reals can be promoted to a global
quasi-conformal conjugacy between them. This is the contents of the
{\it pull-back theorem} below. 

Given a domain ${\cal O} \subseteq {\bf C}$  
symmetric about the real axis, we say that a homeomorphism $\psi : {\cal O} 
\to \psi({\cal O}) \subseteq {\bf C}$ is {\it symmetric} if it commutes with complex 
conjugation and satisfies $\psi({\cal O}^+) = \psi({\cal O})^+$
(where $A^+ =A\cap \{{\rm Im}\, z >0\}$). A given 
holomorphic commuting pair $\Gamma$ is said to have {\it geometric boundaries} 
if its bow-tie $({\cal D}, {\cal O}_\xi, {\cal O}_\eta,
{\cal O}_\nu)$ is such that $\partial{\cal D}$ and $\partial {\cal U}$ are
$K$-quasicircles for some $K \ge 1$, where ${\cal U} = {\cal O}_\xi \cup {\cal
O}_\eta \cup {\cal O}_\nu$. The smallest such $K$ together with the number
${\rm mod}\,({\cal D} \setminus D)$ are the {\it geometric parameters}
of $\Gamma$.   
 
\proclaim{Theorem III.1}.{Let $\Gamma_0, \Gamma_1$ be holomorphic commuting 
pairs with the same irrational rotation number and the same height, assume
they have geometric boundaries, and suppose $h : J_0 \to J_1$ is a
$k$-quasisymmetric conjugacy between the restrictions  
of $\Gamma_0, \Gamma_1$ to their respective large dynamical intervals. Then 
there exists a quasiconformal conjugacy $H : {\cal D}_0 \to {\cal D}_1$  
between $\Gamma_0$ and $\Gamma_1$ which extends $h$ and whose maximal  
dilatation depends only on $k$ and on the geometric parameters of both pairs.} 

One essential difference from Sullivan's original pull-back theorem must be
observed. The straightening theorem of Douady-Hubbard asserts that every
quadratic-like mapping is qc-conjugate to a quadratic polynomial [DH].
In particular, quadratic-like mappings have no wandering domains, after
another well-known theorem of Sullivan [S$_4$], and the pull-back
argument runs smoothly for them. By contrast, we are only able to rule out
wandering domains {\it a posteriori}, see Theorem IV.2. The
technical tool we use to deal with them in the proof of Theorem III.1
is the following {\it qc-sewing lemma} due to L. Bers (cf.~[B],
[Ric]). 
 
\proclaim{Lemma III.2}.{Let $\phi : {\cal O} \to \phi({\cal O}) \subseteq 
\widehat{\bf C}$ be a homeomorphism of an open set ${\cal O} \subseteq  
\widehat{\bf C}$ onto its image, let $\Lambda \subseteq {\cal O}$ be closed in  
$\widehat{\bf C}$ and assume that: 
(a) $\phi|\Lambda$ agrees with the restriction to $\Lambda$ of a  
$K_1$-quasiconformal homeo defined on some neighborhood of $\Lambda$; 
(b) $\phi | ({\cal O} \setminus \Lambda)$ is $K_2$-quasiconformal. Then $\phi$ 
is $K$-quasiconformal with $K \le \max\{K_1, K_2\}$. \hfill\endproof}
 
We need a few other geometric facts.
Let $A \subseteq \disk$ be a ring domain having $\partial\disk$
as its outer boundary. Set $\delta = \inf\{d(z,\partial {\disk }) : z
\in {\disk } \setminus A \}$ and suppose $A $ does not
contain the origin. Then we have the following inequalities due to
Teichm\"uller 
$$
{1 \over 2\pi}\log\, ({1 \over 1 -
\delta}) \le
{\rm mod}\,A  \le {1 \over 2\pi}\log\,\Psi({\delta \over 1 -
\delta}) \ , \eqno{(9)}
$$ 
where $\Psi$ is a universal monotone increasing function
(cf.~[A$_1$]). These inequalities may be combined with K\"obe's
distortion theorem and the Ahlfors-Beurling extension to yield proofs
of the following three lemmas.
Recall that a {\it $K$-quasidisk} is the image of a round disk in the 
extended complex plane under a global $K$-quasiconformal mapping. 

\proclaim{Lemma III.3}.{Let $Q_0, Q_1 \subseteq {\bf C}$ be $K$-quasidisks, 
symmetric with respect to the real axis and satisfying $\overline Q_0  
\subseteq Q_1$. Then the Jordan regions $(Q_1 \setminus Q_0)^{\pm}$, $Q_0 \cup 
Q_1^{\pm}$ are all $K'$-quasidisks with $K'$ depending only on $K$ and 
${\rm mod}\,(Q_1 \setminus Q_0)$. \hfill\endproof}
 
\proclaim{Lemma III.4}.{Let $I_0, I_1 \subseteq {\disk } \cap \real$ be closed  
intervals and let $\phi : I_0 \to I_1$ be a $k$-quasisymmetric homeomorphism. 
Then $\phi$ has a $K$-quasiconformal extension to a self-mapping of ${\disk }$ 
which is symmetric about the real axis and whose maximal dilatation $K$ 
depends only on $k$, ${\rm mod}\,({\disk } \setminus I_0)$ and ${\rm
mod}\,({\disk } \setminus I_1)$.  
\hfill\endproof} 
 
\proclaim{Lemma III.5}.{Let $A_0, A_1$ be disjoint closed arcs in  
$\partial {\disk }$ and let $Q$ be the oriented conformal quadrilateral  
determined by ${\disk }, A_0$ and $A_1$. If $h : \partial {\disk } \to
\partial {\disk }$ is a homeomorphism such that $h|(\partial {\disk } \setminus A_i)$
is $k$-quasisymmetric for $i = 0,1$ then $h$ is $k'$-quasisymmetric with $k'$
depending only on $k$ and ${\rm mod}\, Q$. \hfill\endproof } 

One more result
is needed before we move on to the proof of Theorem III.1.
Let $\Gamma_0$ and $\Gamma_1$ be holomorphic commuting pairs and $F_0, F_1$ 
be their corresponding shadows (section II), and suppose
$h : J_0 \to J_1$ is a conjugacy between the restrictions $F_i|J_i$. 
 
\proclaim{Lemma III.6}.{Let $\psi : {\cal D}_0 \to {\cal D}_1$ be any  
symmetric homeomorphic extension of $h$. Then there exists a symmetric 
homeomorphism $\widetilde\psi : {\cal U}_0 \to {\cal U}_1$ such that $F_1
\circ \widetilde\psi = \psi \circ F_0$ which is still an extension of $h$.} 

\proof
Writing $X_i = J_i \cup F_i^{-1}(J_i)$, 
we know by Proposition II.4 that $F_i|({\cal U}_i \setminus X_i) $ is a
regular 3-fold covering map onto $ {\cal D}_i^+ \cup {\cal D}_i^-$.
Thus we can lift the restriction $\psi|{\cal D}_0^+ \cup {\cal D}_0^-$ through
the $F_i$'s to get a homeomorphism $\widehat\psi : {\cal U}_0 \setminus X_0
\to {\cal U}_1 \setminus X_1$; such lift is uniquely determined if we require
in addition that it be symmetric. Since $F_i(X_i) \subseteq J_i$, and using
the fact that $F_i$ is schlicht when restricted to each of the six components
of ${\cal U}_i \setminus X_i$, we deduce that $\widehat\psi$
extends to a symmetric homeomorphism $\widetilde\psi : {\cal U}_0 \to {\cal
U}_1$, satisfying $F_1 \circ \widetilde\psi = \psi_0 \circ F_0$
everywhere by continuity. As $\psi|J_0 \equiv h$ and since $J_i$ is $F_i$-forward invariant,
we conclude that $\widetilde\psi|J_0 \equiv h$ also. 
\hfill\endproof

\ni{\it Proof of Theorem III.1.\enspace} 
By Lemma III.4 and the Riemann mapping theorem, there exists a
symmetric quasiconformal homeomorphism $G : {\cal D}_0 \to {\cal
D}_1$ extending $h$ and whose maximal dilatation depends only on $k$
and the geometric parameters. Applying Lemma III.6 to $\psi = G$
yields a symmetric quasiconformal lift $\widetilde G : {\cal U}_0 \to
{\cal U}_1$ with $K(\widetilde G) = K(G)$, still satisfying
$\widetilde G|J_0 \equiv h$. By Lemma III.3, the Jordan domains 
$({\cal D}_i  \setminus {\cal U}_i)^+$ and ${\cal U}_i \cup {\cal
D}_i^{-}$ are $K'$-quasidisks with $K'$ depending only on the geometric
parameters. 
Let $\tau_0$, $\tau_1$ be $K'$-quasiconformal mappings of the plane such that
$\tau_i(\disk)= ({\cal D}_i\setminus{\cal U})^+$. Let
$\beta :\partial\disk \to \partial\disk$ be given by 
$$
\beta(z)=\cases{\tau_1^{-1}\circ {\widetilde G}\circ\tau_0 (z) & if $z\in
\tau_0^{-1}(\partial({\cal U}_0 \cup {\cal D}_0))$ \cr
{} & {}\cr
\tau_1^{-1}\circ G\circ\tau_0 (z) & if $z\in
\tau_0^{-1}((\partial D)^+)$ \cr
} \ .
$$
By Lemma III.5, $\beta$ is quasisymmetric (with $k(\beta)$ depending only on
the geometric parameters). Let $B:\disk\to\disk$ be the
Ahlfors-Beurling extension of $\beta$ and let $\widehat
G=\tau_1\circ B\circ\tau_0^{-1} : ({\cal D}_0 \setminus {\cal U}_0)^+ \to
({\cal D}_1 \setminus {\cal U}_1)^+$.
We have $\widehat G \equiv \widetilde G$ over $(\partial {\cal
U}_0)^+$ and $\widehat G \equiv G$ over the remaining part of the boundary of
$({\cal D}_0 \setminus {\cal U}_0)^+$. Then, let $H_1 : {\cal D}_0 \to {\cal
D}_1$ be given by
$$
H_1(z) = \cases{\widetilde G(z) &if $z \in {\cal U}_0$ \cr 
{} & {}\cr
\widehat G(z) &if $z \in \overline{({\cal D}_0 \setminus {\cal
U}_0)^+}$ \cr  
{} & {}\cr
\sigma(\widehat G(\sigma z)) &if $z \in ({\cal D}_0 \setminus {\cal
U}_0)^-$ \cr} \ \ ,
$$  
where $\sigma$ denotes complex conjugation. This map is a  
quasiconformal homeomorphism with $K(H_1) = \max\{K(\widehat G),
K(\widetilde G)\}$, and $H_1|J_0 \equiv h$. 

Now we may start the {\it pull-back routine}. From $H_1$, we define 
inductively a sequence of homeomorphisms $H_n : {\cal D}_0 \to {\cal D}_1$
by 
$$
H_n(z) = \cases{H_1(z) &if $z \in ({\cal D}_0 \setminus {\cal U}_0)$
\cr 
{} & {}\cr
\widetilde H_{n-1}(z) &if $z \in {\cal U}_0$ \cr}
$$ 
where $\widetilde H_{n-1} : {\cal U}_0 \to {\cal U}_1$ is the lift
we obtain applying Lemma III.6 to  
$\psi = H_{n-1}$. The map $H_1$ has been constructed so that
$\widetilde H_1(z) = H_1(z)$ for each $z \in \partial {\cal U}_0$; it
follows inductively from Lemma III.6 that each $H_n$ is a
symmetric quasiconformal homeomorphism with $K(H_n) = K(H_1)$, and
that $H_n|J_0 \equiv h$ for all $n$. 
By the compactness principle for qc-mappings (cf.~[A$_1$]), this
sequence has a limit $H_\infty : {\cal D}_0 \to {\cal D}_1$. We have
$K(H_\infty) \le K(H_1)$ and $H_\infty | J_0 \equiv h$ as well. 

Notice that $\{H_n\}$ has the following {\it stabilization property} : 
if $z \in {\cal U}_0$ then $H_n \circ F_0(z) = F_1 \circ H_n(z)$ if and only if 
$H_{n+1}(z) = H_n(z)$. Let $E$ be the set of all $z \in {\cal U}_0$ which
iterated finitely many times by $F_0$ either land outside ${\cal U}_0$, where
$H_n \equiv H_1$ for all $n$, or land on $J_0$, which is forward invariant and
where $H_n \equiv h$ for all $n$. Then for every $z \in E$ the sequence
$\{H_n(z)\}$ is eventually constant, hence eventually equal to
$H_\infty(z)$. The stabilization property gives us
$H_\infty \circ F_0(z) = F_1 \circ H_\infty(z)$ for all  $z \in E$. Since $X_0
\subseteq E$, we have $({\cal U}_0 \cap \overline E) \setminus E \subseteq
{\cal U}_0 \setminus X_0$, where $F_0$ is continuous, and so for all $z \in
{\cal U}_0 \cap \overline E$ we have $H_\infty \circ F_0(z) = F_1 \circ
H_\infty(z)$ also.   

If $\Omega \subseteq {\cal U}_0 \setminus \overline E$ is a connected
component then the restriction $F_0 | \Omega$ is schlicht. Since $E$
is backward invariant, it follows by induction that $F_0^n(\Omega)
\subseteq {\cal U}_0 \setminus \overline E$ is a connected component
also, for all $n > 0$.  
Observe that $\partial\Omega\cap X_0$ consists 
of at most one point. For if $a, b \in \partial\Omega\cap X_0$ are two 
distinct points, then mapping $\Omega$ forward if necessary we may assume 
that $a, b \in J_0$ and choose $n > 0$ so that the points $F_0^n(a)$,
$F_0^n(b)$ lie in opposite sides of zero inside $J_0$. Since $\Omega$
is connected, this forces $F_0^n(\Omega) \cap X_0 \ne \emptyset$, a
contradiction.

Next, suppose $F_0^n(\Omega) = \Omega$ for some $n > 0$; there is no loss 
of generality in assuming that $\Omega \subseteq {\cal D}_0^+$. Then 
there exists an inverse branch $\Phi : {\cal D}_0^+ 
\to {\cal D}_0^+$ of $F_0^n$ such that $\Phi(\Omega) = \Omega$. Since  
${\cal D}_0^+$ is a Jordan domain, we known by the Denjoy-Wolff theorem 
(cf. [Mil], [S$_3$]) that either there exists 
$z \in {\cal D}_0^+$ such that $\Phi(z) = z$, necessarily attracting 
because $\Phi({\cal D}_0^+) \subseteq {\cal U}_0^+ \ne {\cal D}_0^+$,  
or there exists $z \in \partial{\cal D}_0^+$ such that  
$\Phi(z) = z$ ($\Phi$ extends continuously to $\partial{\cal D}_0^+$). The 
first possibility is incompatible with $\Phi(\Omega) = \Omega$, while
the second implies that $z$ is in the large dynamical interval of
$\Gamma_0$, which is impossible because $F_0$ has no periodic points
there. We deduce that each connected component of ${\cal U}_0
\setminus \overline E$ is a {\it wandering domain}, i.e. its forward
images are pairwise disjoint. Therefore $H_\infty$
conjugates $F_0$ and $F_1$ everywhere except along the grand-orbits 
of wandering domains.   

Now we perform a sequence of quasiconformal sewings 
in order to change $H_\infty$ into a global conjugacy between both 
pairs. Partitioning the connected components of ${\cal U}_0 \setminus
\overline E$ into grand-orbit equivalence classes and selecting one
representative from each 
class yields countably many domains $\{\Omega_n\}_{n \ge 1}$. First
we change  $H_\infty$ along the forward orbit of $\Omega_1$. Let
$\varphi_0 : \overline\Omega_1 \to \overline{F_0(\Omega_1)}$ denote
the homeomorphic extension of $F_0$ to the closure of $\Omega_1$. Then
$\varphi_0(z) = F_0(z)$ for all $z \in \overline\Omega_1$ with at most 
one exception $z_0 \in \partial\Omega_1$. Similarly, define $\varphi_1 : 
\overline{H_\infty(\Omega_1)} \to \overline{F_1H_\infty(\Omega_1)}$ as the 
homeomorphic extension of $F_1$ to the closure of $H_\infty(\Omega_1)$.  
Once again $\varphi_1 \equiv F_1$ with at most one exception $z_1 \in
\partial H_\infty(\Omega_1)$. 
Let $\phi : \overline{F_0(\Omega_1)} \to \overline{F_1H_\infty(\Omega_1)}$  
be given by $\phi = \varphi_1 \circ H_\infty \circ \varphi_0^{-1}$. Then 
$\phi$ is a $K(H_\infty)$-quasiconformal homeomorphism and {\it a priori}  
agrees with $H_\infty$ over $\partial F_0(\Omega_1)$ except possibly  
at one point, 
so by continuity of both maps $\phi \equiv H_\infty$ everywhere along  
$\partial F_0(\Omega_1)$. Hence, if we set $\psi^{(1)} \equiv H_\infty$  
over ${\cal D}_0 \setminus F_0(\Omega_1)$ and $\psi^{(1)} \equiv \phi$ over  
$F_0(\Omega_1)$ we get by the qc-sewing lemma a  
$K(H_\infty)$-quasiconformal homeo $\psi^{(1)} : {\cal D}_0 \to {\cal D}_1$ 
which satisfies the conjugacy equation $\psi^{(1)} \circ F_0(z) =  
F_1\circ \psi^{(1)}(z)$ for all $z \in ({\cal U}_0 \cap \overline E)
\cup \Omega_1$. 
Repeating this argument with $\psi^{(1)}$ replacing $H_\infty$ and 
$F_0(\Omega_1)$ replacing $\Omega_1$ and so on inductively, we get a
sequence $\psi^{(n)} : {\cal D}_0 \to {\cal D}_1$ of uniformly ($\le
K(H_\infty)$) quasiconformal mappings. Again by the compactness
principle we extract a limit $\psi_1$ of $\{\psi^{(n)}\}$. Feeding
$\psi_1$ into our pull-back routine in place of $H_1$ and once again
going to a limit yields a quasiconformal homeo $H_{1,\infty} : {\cal
D}_0 \to {\cal D}_1$ with $K(H_{1,\infty}) \le K(H_\infty)$ that, by
the stabilization property, conjugates $F_0$ and $F_1$ not only on
${\cal U}_0 \cap \overline E$ but also along the full grand-orbit of
$\Omega_1$. 
Proceeding inductively, we take care of the full 
grand-orbits of $\Omega_2, \Omega_3, \ldots$ through partial quasiconformal 
conjugacies $H_{2,\infty}, H_{3,\infty}, \ldots$ satisfying $K(H_{n,\infty}) 
\le K(H_\infty)$ as well as $H_{n, \infty} | J_0 \equiv h$ for all $n$. 
Going to a limit one final time we get $H : {\cal D}_0 \to {\cal D}_1$, a  
global quasiconformal conjugacy with $K(H) \le K(H_\infty)$ which is still 
an extension of $h$.
\hfill\endproof

\ni{\it Remark.\enspace} Let $\Gamma$ be a holomorphic commuting pair with
irrational rotation number. By analogy with the case of quadratic-like
mappings, we define the {\it filled-in limit set} of $\Gamma$ to be
$$
{\cal K}_\Gamma\;=\; \overline{\bigcap_{n\ge 0}\, F^{-n}({\cal D})} 
\ . 
$$
We also let the {\it limit} or {\it Julia set} $J_\Gamma$ be the set
${\cal K}_\Gamma$ minus the union of all grand-orbits of wandering
domains. As the proof of the pull-back theorem shows, ${\cal K}_\Gamma$ has no
other stable interior components. In the next section we will show in fact
that ${\cal K}_\Gamma = J_\Gamma$. It is not difficult to see that ${\cal
K}_\Gamma\cap \partial {\cal U}$ is the disjoint union of six arcs.

\bigskip

\ni\centerline{\smc IV. \enspace Existence and limit set qc-rigidity of
holomorphic commuting pairs} 
\smallskip
\ni 
For each $0\le \theta <1$, let $E_\theta : {\bf C} \to {\bf C}$
be the entire mapping given by $E_\theta(z) = z + \theta - {1 \over
2\pi}\sin(2\pi z)$. Since $E_\theta \circ T = T\circ E_\theta$, where
$T$ is the translation $z \mapsto z + 1$, $E_\theta$ is the lift to the
complex plane of a holomorphic self-mapping of the 
cylinder, $f_\theta : {\bf C}/\bz \cong {\bf C}^*\hookleftarrow$. Moreover, the restriction $E_\theta|\real$ maps the
real axis onto itself and satisfies $E_\theta'(x) \ge 0$ for all $x
\in \real$, and equality holds iff $x \in \bz$ (these constitute all the 
critical points of $E_\theta $). Therefore the restriction $f_\theta |{S}^1 $
is a critical circle homeomorphism with rotation number, say, $\rho(\theta)$.
It is well-known that $\theta\mapsto \rho(\theta)$ is a continuous,
non-decreasing map of $[0,1)$ onto itself such that the interval $\rho^{-1}(t)
\subseteq [0,1)$ degenerates to a point whenever $t$ is irrational
(see [H$_1$]).  
 
With the family $\{E_\theta\}$ at hand we shall construct in this section 
examples of holomorphic commuting pairs with arbitrary rotation number and
arbitrary height. More precisely, we shall prove the following theorem. 

\proclaim{Theorem IV.1}.{For each $n \ge 0$ and each $\theta$ 
such that $\rho(\theta)$ has a continued fraction expansion of 
length at least $n + 1$, the real commuting pair determined by  
$(f_\theta^{q_n}, f_\theta^{q_{n+1}})$ extends to a holomorphic commuting 
pair $\Gamma_{n,\theta}$ with geometric boundaries. The family  
$\{\Gamma_{n,\theta}\}$ runs through all possible pairs of combinatorial 
invariants at least once, and for each $(m, \rho) \in {\integer} \times 
{[0,1)}$ with $m \ge 2$ there exist countably many $(n,\theta) \in
{\integer} \times  {[0,1)}$ such that $\Gamma_{n,\theta}$ has height
$m$ and rotation number $\rho$. }

When combined with the results of section III, this construction 
yields two crucial properties of holomorphic commuting pairs, which we express
as follows. 

\proclaim{Theorem IV.2}.{Let $\Gamma$ be a holomorphic commuting pair 
with geometric boundaries and irrational rotation number. Then $\Gamma$ 
has no wandering domains and admits no non-trivial, symmetric, invariant 
Beltrami differentials entirely supported in its limit set. }
 
This theorem allows holomorphic commuting pairs to be parametrized by
conformal structures supported on the outer annulus of a fixed model, cf. next
section. 

The main analytic tool to be used in the proof of Theorem IV.1 is the following
growth estimate.  
 
\proclaim{Lemma IV.3}.{There exist a positive constant $C_0$ and a positive 
monotone non-decreasing function $\varphi(s)$ defined for $s \ge 0$ such 
that if $|y| \ge \varphi(|x|)$ then $|E_\theta(x + iy)| \ge C_0\exp(\pi 
|y|)$.}
 
\proof 
When $\theta = 0$, a straightforward computation  
yields 
$$
{\big |}E_0(x+iy){\big |}^2 \;=\;
{1 \over 4\pi^2}\cosh^2{(2\pi y)} + \big [x^2 + y^2 - {1 \over
4\pi^2}\cos^2{(2\pi x)}\big ] 
$$
$$
-  {1 \over \pi}\big [x\sin{(2\pi x)}\cosh{(2\pi y)} +  
y\cos{(2\pi x)}\sinh{(2\pi y)}\big ] {\rm .}  
$$
The first expression between brackets is positive as soon as, say, $|y| \ge 1$, 
while the second is dominated by $(|x| + |y|)\cosh{(2\pi y)}$. Thus, if 
$|y| \ge 1$ we have
$$
{\big |}E_0(x + iy){\big |}^2 \ge {1 \over 4\pi^2}\big [\cosh{(2\pi y)} -  
4\pi(|x| + |y|)\big ]\cosh{(2\pi y)}  \ .
\eqno (10)
$$ 
Now, let 
$$
\epsilon(t) = {1 \over 4\pi}\cosh{(2\pi t)} - t - 1 \ .
$$
This strictly convex function has a minimum at a certain $t_0> 0$ such
that $\epsilon(t_0) < 0$. Hence 
for each $s \ge 0$ there exists a unique $\overline\varphi(s) > t_0$
such that $\epsilon(\overline\varphi(s)) = s$. Since
$\epsilon(t)$ is strictly increasing for $t \ge t_0$, so is
$\overline\varphi(s)$ for $s \ge 0$, and $t \ge \overline\varphi(s)$
implies $\epsilon(t) \ge s$. Setting $\varphi(s) = \max\{1,
\overline\varphi(s)\}$ and observing that the expression between
brackets in $(10)$ is equal to $4\pi[\epsilon(|y|) + 1 - |x|]$, we
deduce that if $|y| \ge \varphi(|x|)$ then
$$
{\big |}E_0(x + iy){\big |}^2 \ge {1 \over \pi}\cosh(2\pi|y|) \ge {1 \over 2\pi}\exp 
(2\pi|y|) {\rm .} \eqno (11)
$$ 
On the other hand, when $0 < \theta < 1$ we have $E_\theta(z) =  
E_0(z) + \theta$, so that $|E_\theta(z)| \ge |1 - |E_0(z)|^{-1}|.|E_0(z)|$. 
Therefore, if $|y| \ge \varphi(|x|)$, we have by $(11)$
$$
{\big |}E_\theta(x+ 
iy){\big |} \ge {1 \over \sqrt{2\pi}}[1 - e^{-\pi}\sqrt{2\pi}]\exp(\pi|y|)
\ . \eqno\square
$$

 
We divide the work required to prove Theorem IV.1 into several steps. Let us
fix $\theta$ for the time being and write $\rho(\theta) = [r_0, r_1, \ldots,
r_n, \ldots]$. We conform with the notation established in the first section,
so that, in its irreducible form, ${\d {p_n \over q_n}} = [r_0, r_1, \ldots,
r_{n-1}]$ satisfies $p_0 = 0$, $q_0 = 1$; $p_1 = 1$, $q_1 = r_0$ and for $n
\ge 1$, $p_{n+1} = r_np_n + p_{n-1}$, $q_{n+1} = r_nq_n + q_{n-1}$. 

We need a brief geometric description of the map $E_\theta$.
The pre-image of the real axis under $E_\theta$ consists of $\real$ itself 
together with the family of analytic curves 
$$
{\cal S}_{\pm}^{(k)} :\  x=k\pm {1\over{2\pi}}\arccos{\Big[
{{-2\pi|y|}\over {\sinh{(2\pi y)}}}\Big]} \ ,
$$ 
where $k \in \bz$, arising as solutions to ${\rm Im}\,E_\theta(x+iy) 
= 0$. For each $k \in \bz$, the curves ${\cal S}_+^{(k)}$ and ${\cal S}_-^{(k)}$
meet at the critical point $c_k = k$, and are both asymptotic to the vertical
lines $x = k \pm {1 \over 4}$. Notice that each $c_k$ is a critical point of
cubic type.   
In the upper half-plane ${\bf C}^+$, let $V_k$ be the simply-connected region 
bounded by the arcs ${\cal S}_+^{(k-1)} \cap {\bf C}^+$ and  
${\cal S}_-^{(k)} \cap {\bf C}^+$ and 
the interval $[k-1, k] \subseteq \real$. Then $E_\theta | V_k$ is schlicht  
onto ${\bf C}^+$; we let $\phi_k : {\bf C}^+ \to V_k$ denote the  
corresponding inverse. Similarly, let $W_k 
\subseteq {\bf C}^+$ be the simply-connected region bounded by  
${\cal S}_-^{(k)} \cap \overline{{\bf C}}^+$ and ${\cal S}_+^{(k)} \cap {\bf C}^+$,  
observe that $E_\theta |  
W_k$ is schlicht onto ${\bf C}^-$ and let $\psi_k : {\bf C}^- \to W_k$ be the 
corresponding inverse.
 
Now let $A_n \subseteq {\bf C}^+$ be the unique connected component of  
$(E_\theta^{q_n})^{-1}({\bf C}^+)$ whose closure contains the point  
$T^{-p_{n+1}}\circ E_\theta^{q_{n+1}}(0) \in \real$. Similarly, let $B_n  
\subseteq {\bf C}^+$ be the unique connected component of  
$(E_{\theta}^{q_{n+1}})^{-1} 
({\bf C}^+)$ such that $T^{-p_n}\circ E_\theta^{q_n}(0) \in \overline B_n$. We 
have either $A_n \subseteq V_0$ and $B_n \subseteq V_1$ or $A_n \subseteq V_1$ and 
$B_n \subseteq V_0$, depending on whether $n$ is even or odd, respectively 
(Figure 2 illustrates the former case). 
 
\proclaim{Lemma IV.4}.{For each $n \ge 0$ there exists a unique $q_n$-tuple  
$(k_1, k_2, \ldots, k_{q_n})$ with $0 = k_1 \le k_2 \le \cdots \le k_{q_n}  
\le p_n + 1$ such that $A_n = \phi_{k_1}\circ \phi_{k_2}\circ \cdots \circ 
\phi_{k_{q_n}}({\bf C}^+)$. A similar statement holds for $B_n$.} 
 
\proof 
This is an easy consequence of the fact that $0 \le E_\theta^j(0) 
< p_n + 1$ for $j = 0, 1, \ldots, q_n$, for all $n \ge 0$, which in turn 
follows from the very definitions of $p_n, q_n$.
\hfill\endproof

\proclaim{Lemma IV.5}.{Let $f$ be a circle homeomorphism with
$\rho(f) = $ $[r_0$, $r_1$, $\ldots$, $r_n$, $\ldots]$, let $c \in
{S}^1$, and for each $n \ge 1$ let $J_n \subseteq {S}^1$ be the 
closed interval of endpoints $c$ and $f^{q_{n-1}-q_n}(c)$ containing
$f^{q_{n-1}}(c)$. If $j < q_n$ is such that $f^{-j}(c)$ belongs to
$J_n$, then $j \le 0$.\hfill\endproof }  
 
Let us use the notation $\<\alpha, \beta\>$ to represent a closed interval on
the line with endpoints $\alpha$ and $\beta$, irrespective of order.
 
\proclaim{Lemma IV.6}.{For each $n \ge 0$ we have ${\overline A_n
\cap \real }= \< \alpha_n, 0\>$ and $\overline B_n \cap \real =$   
$\< 0, \beta_n\>$, where $\alpha_0 = -1$, $\beta_0 = \alpha_1$ and for $n\ge
1$ the points $\alpha_n, \beta_n \in \real$ are 
uniquely determined by the requirements: 
$T^{-p_n}\circ E_\theta^{q_n}(\alpha_n) = T^{-p_{n-1}}\circ  
E_\theta^{q_{n-1}}(0)$ 
and $T^{-p_{n+1}}\circ E_\theta^{q_{n+1}}(\beta_n) = T^{-p_n}\circ  
E_\theta^{q_n}(0)$.}
 
\proof 
Consider $f = f_\theta$ and take $c$ to be the critical point of 
$f_\theta$. Then Lemma IV.5 says that there can be no critical points for
$f_\theta^{q_n}$ in the interior of $J_n$, for by the chain rule these are
precisely the pre-images $f_\theta^{-j}(c)$ with $0 \le j < q_n$. The result
follows. 
\hfill\endproof

Given $R > 0$, let ${\cal D}_R = \{z : |z| < R\}$ and let $A_{n,R}$ be  
the unique connected component of $(T^{-p_n}\circ E_\theta^{q_n})^{-1} 
({\cal D}_R^+)$ contained in $A_n$. Let $B_{n,R}$ be 
similarly defined. If $R$ is sufficiently large 
($R > p_n + 1$ is good enough) we see that ${\overline A_{n,R} \cap \real }=  
\overline A_{n} \cap \real$ and ${\overline B_{n,R} \cap \real }= \overline B_{n}  
\cap \real$ for $n \ge 0$. It is clear that both $A_{n,R}$ and  
$B_{n,R}$ are Jordan domains, in fact quasidisks, and that they are mapped 
respectively by $T^{-p_n}\circ E_\theta^{q_n}$ and $T^{-p_{n+1}}\circ  
E_\theta^{q_{n+1}}$ bijectively onto ${\cal D}_R^+$.

\proclaim{Lemma IV.7}.{For every sufficiently large $R$ we have $\overline A_{n,R} 
\subseteq {\cal D}_R \cap \overline{{\bf C}^+}$ and $\overline B_{n,R}  
\subseteq {\cal D}_R \cap \overline{{\bf C}^+}$.}
 
\proof
For $s, R$ positive numbers, 
let
$$
\delta(s, R) = \varphi(s) + {1 \over \pi}\log^+(C_0^{-1}R) \ ,
$$
where $\varphi$ and $C_0$ are given by Lemma IV.3. Then $|y| \ge \delta(|x|,
R)$ implies $|E_\theta(x + iy)| \ge R$, which in turn means that $E_\theta(x +
iy) \in {\bf C} \setminus {\cal D}_R$. Therefore, for each $k \in \bz$ we have 
$$
\phi_k(\overline{{\cal D}_R^+}) \subseteq \overline V_k \cap \Big\{x + iy : y
\le \delta(|x|, R)\Big\} \ .
$$  
Let $V_{k, R}$ denote this last intersection. Since $\delta(s,R)$ has logarithmic growth 
in $R$, every sufficiently large $R$ satisfies the inequality $R > p_n + 1 +  
\delta(p_n + 1, 2R)$; for a given $R$ as such, if $0 \le k \le p_n +1$ and 
$z$ is any point in $V_{k,2R}$ with $z = x + iy$, then
$$
|z| \le |x| +  
\delta(|x|, 2R) \le p_n + 1 + \delta(p_n + 1, 2R) < R  \ ,
$$ 
and so it follows that $z \in {\cal D}_R \cap \overline{{\bf C}^+}$. Thus,
if $0 \le k \le p_n + 1$ then $\phi_k(\overline{{\cal D}_{2R}^+}) 
\subseteq {\cal D}_R \cap  \overline{{\bf C}^+} \subseteq \overline{{\cal
D}_{2R}^+}$. Since $T^{p_n}(\overline{{\cal D}_R^+}) \subseteq \overline{{\cal
D}_{2R}^+}$, if we take $(k_1, k_2, \ldots, k_{q_n})$ as
in Lemma IV.4 we deduce that
$$
\overline A_{n,R} = \phi_{k_1} \circ
\phi_{k_2}\circ \cdots \circ  
\phi_{k_{q_n}} (T^{p_n}\overline{{\cal D}_R^+}) \subseteq \phi_{k_1} \circ  
\phi_{k_2} \circ \cdots \circ \phi_{k_{q_n}}(\overline{{\cal D}_{2R}^+}) 
\subseteq {\cal D}_R \cap \overline{{\bf C}^+} \ .
$$
This proves the first inclusion; the second is proved in similar fashion.
\hfill\endproof  

{\ni
{\it Remark.\enspace} Observe that if we define ${\cal U}_{n,R} =
\phi_{k_2}\circ \phi_{k_3} \circ \cdots \circ \phi_{k_{q_n}}({\cal D}_R^+)$
and set $A_{n,R}' = \phi_1({\cal U}_{n,R})$ and $A_{n,R}'' =
\psi_0\sigma({\cal U}_{n,R})$, where $\sigma : {\bf C} \to {\bf C}$ is complex
conjugation, then the above argument applies {\it mutatis mutandis} to yield
$\overline{A_{n,R}'} \subseteq {\cal D}_R \cap \overline{{\bf C}^+}$,
$\overline{A_{n,R}''} \subseteq {\cal D}_R \cap \overline{{\bf C}^+}$  
as well, for every sufficiently large $R$ and all $n \ge 0$. }
\smallskip 
 
\ni {\it Proof of Theorem IV.1.}\enspace
Given $n \ge 0$, let $R_n > 0$ be large enough for the conclusion of Lemma IV.7
to hold. Let $\xi_n = T^{-p_n}\circ E_\theta^{q_n}$ and $\eta_n =  
T^{-p_{n+1}}\circ E_\theta^{q_{n+1}}$ and let ${\cal O}_{\xi_n}, {\cal O}_{\eta_n}  
\subseteq {\bf C}$ be the symmetric Jordan domains (quasidisks) such that 
${\cal O}_{\xi_n}^+ = A_{n, R_n}, {\cal O}_{\eta_n}^+ = B_{n, R_n}$.
Then $\xi_n$ and $\eta_n$ commute, and $\overline {\cal O}_{\xi_n}, 
\overline {\cal O}_{\eta_n} \subseteq {\cal D}_{R_n}$, by Lemma IV.7. The  
restrictions $\xi_n|{\cal O}_{\xi_n}$ and $\eta_n|{\cal O}_{\eta_n}$ are schlicht and 
onto their images, which by Lemma IV.5 are ${\cal D}_{R_n} \cap {\bf
C}(\<\xi_n(\alpha_n), \xi_n(0)\>)$ and ${\cal D}_{R_n} \cap {\bf
C}(\<\eta_n(0), \eta_n(\beta_n)\>)$, respectively. 
Also, let ${\cal O}_{\nu_n} \subseteq {\bf C}$ be the connected component of $\xi_n^{-1} 
({\cal O}_{\eta_n})$ containing the origin and let $\nu_n = \xi_n \circ \eta_n$. 
Then the restriction $\nu_n|{\cal O}_{\nu_n}$ is a holomorphic 3-fold branched 
covering map onto its image, $\nu_n({\cal O}_{\nu_n}) = {\cal D}_{R_n} \cap 
{\bf C}(\<\eta_n(0), \xi_n(0)\>)$. Moreover, by the remark following Lemma IV.7, we
have
$$
\overline{{\cal O}_{\nu_n}^+} \subseteq \overline A_{n,R} \cup  
\overline{A_{n,R}'} \cup \overline{A_{n,R}''} \subseteq {\cal D}_{R_n} \cap 
\overline{{\bf C}^+} \ ,
$$
and so $\overline {\cal O}_{\nu_n} \subseteq {\cal D}_{R_n}$. 
It follows at once that $({\cal D}_{R_n}, {\cal O}_{\xi_n}, {\cal O}_{\eta_n},
{\cal O}_{\nu_n})$ is a {\it bow-tie}.

Now we claim that this bow-tie together with the maps $\xi_n, \eta_n, \nu_n$
determine a holomorphic commuting pair $\Gamma_{n,\theta}$ with geometric
boundaries, up to orientation, with rotation number $\rho(\Gamma_{n,\theta}) = [r_{n +1}+1,
r_{n+2}, \ldots]$ and height given by $m(\Gamma_{0,\theta}) = r_0$ when $n =
0$, and by $m(\Gamma_{n,\theta}) = r_n + 1$ when $n > 0$.
We have indirectly checked all conditions in Definition 4, except perhaps condition
H$_6$. We check it for $n > 0$; the case $n = 0$ is just as easy. Using the
commutativity of $T$ with $E_\theta$, Lemma IV.4 and the recurrence relations
defining $p_{n+1}$ and $q_{n+1}$,  
we get
$$
\xi_n^{r_n + 1}(\alpha_n) \;=\; (T^{-p_n}\circ E_{\theta}^{q_n})^{r_n} 
(T^{-p_n} \circ E_{\theta}^{q_n}(\alpha_n))  
\;=\; T^{-p_{n+1}} \circ E_{\theta}^{q_{n+1}}(0) = \eta_n(0) {\rm .} 
$$
Similarly, we have $\eta_n(\beta_n) = \xi_n(0)$. Thus  
condition H$_6$ is satisfied too, and  $m = r_n + 1$ is the height of
$\Gamma_{n,\theta}$. The statement on rotation  numbers is clear.
\hfill\endproof  


$$
\psannotate{
\psboxto(\hsize;0cm){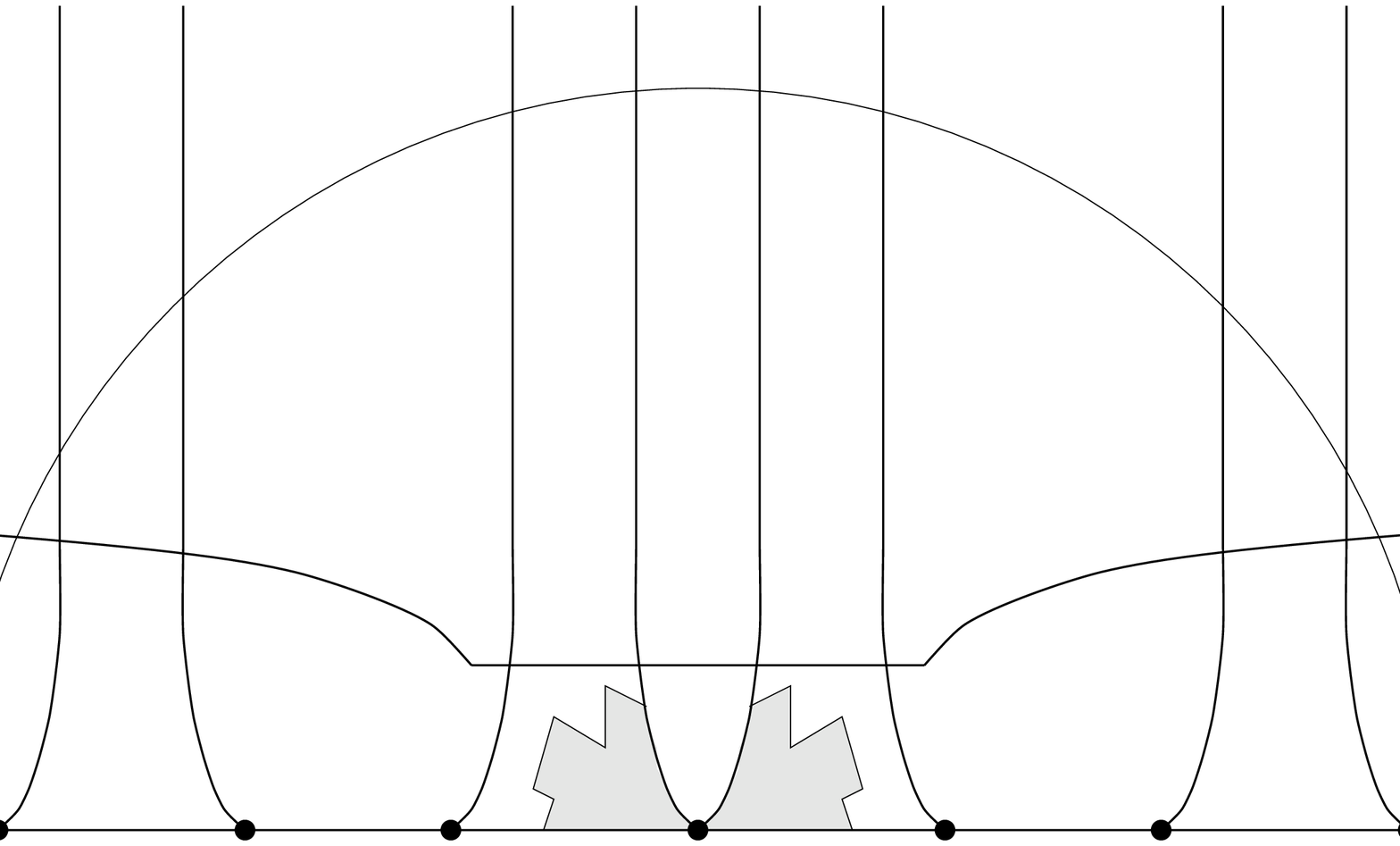}}
{\fillinggrid\at{14.8\pscm}{1.8\pscm}{$0$}
             \at{10.4\pscm}{1.8\pscm}{$-1$}
             \at{18.2\pscm}{1.8\pscm}{$1$}
             \at{21.3\pscm}{1.8\pscm}{$p_n$}
             \at{25.3\pscm}{1.8\pscm}{$Re\,z$}
             \at{12\pscm}{3\pscm}{$A_{n,R}$}
             \at{14.4\pscm}{3\pscm}{$B_{n,R}$}
             \at{6.5\pscm}{7\pscm}{$y=\delta(|x|,2R)$}
             \at{7\pscm}{13\pscm}{${\cal D}_R^+$}
             \at{11\pscm}{15\pscm}{$V_0$}
             \at{14.6\pscm}{15\pscm}{$V_1$}
             \at{21\pscm}{15\pscm}{$V_{p_n+1}$}
             \at{9.9\pscm}{1.85\pscm}{$\alpha_n$}
             \at{14.4\pscm}{1.8\pscm}{$\beta_n$}
}
$$
\centerline{\it Figure 2}
\bigskip

\ni {\it Remark.\enspace } Because of Proposition II.3, once $R_0$ is chosen so 
that the above construction works for $n = 0$, we may take  
$R_n = R_0$ thereafter. If this is done then, for each $n \ge 0$, 
$\Gamma_{n+1, \theta}$ becomes the first renormalization of $\Gamma_{n,\theta}$ 
up to linear rescaling.   
\smallskip
 
Now we turn our attention to Theorem IV.2. We shall extract the stated rigidity
properties of holomorphic commuting pairs from corresponding ones found naturally in
the family $\{f_\theta\}$ introduced above.  
 
If $f : {\bf C}^*\to {\bf C}^*$ is holomorphic, we denote by
$S_f$ the set of singular values of $f$, i.e. points in ${\bf
C}^*$ all neighborhoods $U$ of which are such that $f^{-1}(U)
\buildrel f \over \longrightarrow U$ fails to be a covering map. We
also write $X_f = {\bf C}^* \setminus S_f$, so that $f^{-1}(X_f)
\buildrel f \over \longrightarrow X_f$ is always a covering map. For
example, since $1 \in {\partial }{\disk }$ is the unique critical
point of $f_\theta$, it is easy to see that $S_{f_\theta} =
\{f_\theta(1)\}$; in this case $f_{\theta}^{-1}(X_{f_\theta})$ 
has an infinite discrete complement in ${\bf C}^*$. 
We let $J_f$ be the Julia set of $f$.
 
\proclaim{Lemma IV.8}.{The family $\{f_\theta\}$ is topologically
complete, i.e. every symmetric, normalized holomorphic
self-map of ${\bf C}^*$ which is topologically conjugate to a member of
the family is a member also.}
\proof
Let $f : {\bf C}^* \to {\bf C}^*$ be holomorphic and suppose  
$h : \widehat{{\bf C}} \to \widehat{{\bf C}}$ is an orientation preserving
homeo fixing $0$ and $\infty$ and satisfying $h \circ f_\theta = f
\circ h$. 
Let $A \in {\rm Aut}\,(\widehat{{\bf C}})$ be given by $A(z) = \lambda z$, 
where $\lambda = h\circ f_\theta(1)/f_\theta(1)$. This $A$ is homotopic to $h$
relative to $S_{f_\theta} \cup \{0, \infty\}$, so the covering homotopy theorem 
yields a holomorphic lift $\widehat A : f_\theta^{-1}(X_{f_\theta}) \to
f^{-1}(X_f)$, which is then homotopic to $h$ relative to
$f_\theta^{-1}(S_{f_\theta}) \cup \{0, \infty\}$. 
Some easy topology and the removable singularity theorem show that $\widehat
A$ is M\"obius and fixes $0$ and $\infty$. In
particular, if $f$ is symmetric about
${\partial }{\disk }$ and is normalized so that its critical point
lies at $1 \in {\partial }{\disk }$, then $\widehat A$ is the identity
and $|\lambda| = 1$, say $\lambda = e^{2\pi i\alpha}$. Therefore $f =
A\circ f_\theta\circ \widehat A^{-1} = f_{\theta +\alpha}$.
\hfill\endproof
 
\proclaim{Theorem IV.9}.{The mapping $f_\theta$ has no wandering
domains. Moreover, if $\rho(\theta)$ is irrational then $f_\theta$
admits no non-trivial, symmetric, invariant Beltrami differentials
entirely supported in its Julia set. }
 
\proof
Since $S_{f_\theta}$ is a finite set, the first  
assertion follows from a theorem due to L. Keen [K].
Now suppose $\mu$ is an $f_\theta$-invariant Beltrami differential 
in $\widehat{{\bf C}}$ with support in $J_{f_\theta}$; assume also
that $\mu$ is symmetric about ${\partial }{\disk }$. For all
sufficiently small real $t$, let $h_t : \widehat{{\bf C}}\to
\widehat{{\bf C}}$ be the unique solution to $\overline{\partial }h_t
= (t\mu).{\partial }h_t$ fixing $\{0, 1,\infty\}$ pointwise, and let $f_t
= h_t \circ f_\theta \circ h_t^{-1}$.
Since $t\mu$ is symmetric and $f_\theta$-invariant, each $f_t$ is
symmetric and holomorphic, and has a single critical point at $1 \in
{\partial }{\disk }$. 
Using Lemma IV.8, we have $f_t = f_{\theta_t}$ for some $\theta_t$.
But then $\rho(\theta_t) = \rho(\theta)$ is irrational, so 
$\theta_t = \theta$ for all $t$, by remark in the first paragraph of
this section. Therefore, $h_t$ commutes with $f_\theta$
for all $t$; in particular $h_t$ must permute the elements of $Y_n =
f_\theta^{-n}(1)$, which is discrete in ${\bf C}^*$, for each $n \ge
0$. Since $h_0 = id_{\widehat{{\bf C}}}$ and for each $z \in
\widehat{{\bf C}}$ the path $t \to h_t(z)$ is continuous by the
Ahlfors-Bers theorem [A$_1$], we deduce that $h_t$ fixes $Y_n$
pointwise for all $n \ge 0$, for all $t$. 
But by Montel's theorem, 
$
J_{f_\theta} \subseteq \overline{\bigcup_{n \ge 0} Y_n} \ ,
$ 
so $h_t$ agrees with the identity over $J_{f_\theta}$
for all $t$. Since $h_t$ is conformal off $J_{f_\theta}$, it follows
that $h_t \equiv id_{\widehat{{\bf C }}}$ for all $t$, and so 
$\mu \equiv 0$ a.e.
\hfill\endproof
 
\ni{\it Proof of Theorem IV.2.\enspace}
Combining Theorem III.1 with Theorem IV.1, we know that $\Gamma$ is
conjugate to $\Gamma_{0,\theta}$ for some $\theta$ by a
qc-homeomorphism $H$. Let $\mu$ be a $\Gamma$-invariant Beltrami
differential with support in $J_{\Gamma}$. Then $\mu ' =
H^{*}\mu$ is $\Gamma_{0,\theta}$-invariant. Spreading $\mu '$
through the entire complex plane via the mappings defining
$\Gamma_{0,\theta}$ we get a Beltrami differential $\nu$ invariant
under both $E_{\theta}$ and $T^{-1}\circ E_{\theta}^{r_0}$, and
therefore invariant under $T$ also. Thus $\nu$ 
projects down to a Beltrami differential on the cylinder which is
$f_{\theta}$-invariant and supported in $J_{f_\theta}$. By Theorem
IV.9, this Beltrami differential must vanish a.e., and so $\mu\equiv
0$ a.e. also. A similar argument rules out wandering domains. 
\hfill\endproof

\bigskip

\ni\centerline{\smc V. \enspace Beltrami paths and the Julia-Teichm\"uller
metric} 
\smallskip
\ni 
If $\Gamma_1$ and $\Gamma_2$ are holomorphic commuting pairs with the
same rotation number and the same height, we define
their {\it quasiconformal distance} to be $d_{QC}(\Gamma_1,
\Gamma_2)={1\over 2}\inf_H\log{K(H)}$, where $H$ ranges over all possible
symmetric quasiconformal conjugacies between both pairs. This number is
obviously zero if and only if $\Gamma_1$ is holomorphically
equivalent to $\Gamma_2$. In order to have a space with an actual
metric on it, we proceed just as in the case of Riemann surfaces or
Fuchsian groups.
Thus, let $\Gamma$ be a fixed holomorphic commuting pair and let ${\cal D}$
be the outer disk of its bow-tie. Let ${\rm Def}(\Gamma)$ be the class of 
all holomorphic commuting pairs which are conjugate to $\Gamma$ via a 
symmetric qc-homeomorphism. In ${\rm Def}(\Gamma)$, declare $\Gamma_0$ to be 
equivalent to $\Gamma_1$ iff there exists a symmetric conformal mapping 
${\cal D}_0 \to {\cal D}_1$ conjugating $\Gamma_0$ to $\Gamma_1$. 
Then define the {\it Teichm\"uller space} of $\Gamma$,  
say ${\rm Teich}(\Gamma)$, to be the quotient of  
${\rm Def}(\Gamma)$ by the above equivalence relation.
 
If $[\Gamma_0], [\Gamma_1] \in {\rm Teich}(\Gamma)$, set $d_T([\Gamma_0], 
[\Gamma_1]) = {1\over 2}\inf_H\log K(H)$ where $H$ ranges over all possible symmetric 
qc-conjugacies between any two representatives $\widetilde\Gamma_0 \in [\Gamma], 
\widetilde\Gamma_1 \in [\Gamma_1]$, and where as before $K(H)$ denotes the  
maximal dilatation of $H$. This defines the {\it Teichm\"uller metric}
on ${\rm Teich}(\Gamma)$. 
As in the case of Fuchsian groups, an alternative description of  
${\rm Teich}(\Gamma)$ as an orbit space is available. Observe that if
$G$ is a group of qc-selfhomeomorphisms of ${\cal D}$ and $B^\infty$
is the unit ball of $L^\infty({\cal D}, {\bf C})$ then there is a
natural action $G \times B^\infty \to B^\infty$,
$$
(h,\mu) \mapsto h^*\mu = {\mu_h + (\mu\circ h).{\d\overline h_z \over\d h_z} 
\over 1 + \overline\mu_h.(\mu\circ h).{\d\overline h_z \over\d h_z}}\ ,
$$ 
which consists of taking the pull-back under $h$ of $\mu \in B^\infty$ 
viewed as a Beltrami differential on 
${\cal D}$. Each $h^* : B^\infty \to B^\infty$ is biholomorphic and
$(h^*)^{-1} = (h^{-1})^*$.  
Call $\mu \in B^\infty$ (a) symmetric, if $\mu$ 
commutes with complex conjugation, and (b) $\Gamma$-invariant, if  
$\gamma^*\mu = \mu$ for $\gamma = \xi, \eta, \nu$. 
If we take $G$ to be the group of all symmetric qc-selfhomeos of  
${\cal D}$ which commute with $\Gamma$ and let $M(\Gamma) = \{\mu  
\in B^\infty : \mu$ is symmetric and $\Gamma$-invariant 
$\}$, then the above $G$-action on $B^\infty$ restricts to an action  
$G \times M(\Gamma) \to M(\Gamma)$. Let ${\rm Orb}_G(\Gamma)$ be the
corresponding orbit space. This space can be given the following
metric 
$$
d([\mu_0], [\mu_1]) = 
{1\over 2}\inf\log K(h^{\widetilde\mu_0} \circ (h^{\widetilde\mu_1})^{-1}) \ ,
$$
where the infimum is taken over all $\widetilde\mu_0, \widetilde\mu_1 \in  
M(\Gamma)$ in the $G$-orbits of $\mu_0$ and $\mu_1$, respectively. Here $h^\mu$ denotes 
the unique symmetric qc-homeo ${\cal D} \to {\cal D}$ with $h^\mu(0) = 0$  
and such that $\mu_{h^\mu} = (h^\mu)^*(0) = \mu$; existence and uniqueness 
are guaranteed by the measurable Riemann mapping theorem.

Given $\mu \in M(\Gamma)$, let us consider the Jordan domains ${\cal D}$ and
${\cal O}_{\gamma^\mu} = h^\mu({\cal O}_\gamma)$, and the maps $\gamma^\mu =
h^\mu\circ\gamma\circ (h^\mu)^{-1} : {\cal O}_{\gamma^\mu} \to \gamma({\cal
O}_{\gamma^\mu})$, for $\gamma = \xi, \eta, \nu$, which are holomorphic because
$\mu$ is $\Gamma$-invariant.  
 
\proclaim{Lemma V.1}.{These objects determine a holomorphic commuting pair.}
\proof 
Of all conditions in the definition given in section II, the only one that is
not immediate is H$_4$. We must check that both $\xi^\mu$ and $\eta^\mu$ extend   
holomorphically across some neighborhood of zero, where they ought to   
commute. The point is that, by Proposition II.1, the map $\eta^{-1}\circ\nu$ 
is well-defined over ${\cal O}_\nu$ and agrees with $\xi$ over ${\cal
O}_\xi\cap {\cal O}_\gamma$.
Therefore $(\eta^\mu)^{-1} \circ \nu^\mu = h^\mu\circ(\eta^{-1}\circ\nu) 
\circ(h^\mu)^{-1}$ is well-defined over ${\cal O}_{\nu^\mu}$ 
and agrees with $\xi^\mu$ over ${\cal O}_{\xi^\mu} \cap 
{\cal O}_{\nu^\mu}$. Similarly, $(\xi^\mu)^{-1}\circ \nu$ 
is well-defined over ${\cal O}_{\nu^\mu}$ and extends $\eta^\mu$. It follows 
at once that both extensions commute on that common part of their domains. 
\hfill\endproof
 
If we denote by $\Gamma^\mu$ the resulting holomorphic commuting pair then 
this lemma gives us the right to write formally $\Gamma^\mu = h^\mu \circ 
\Gamma\circ (h^\mu)^{-1}$. Therefore, just as with Fuchsian groups
(cf.~[A$_1$]), we have the following statement.

\proclaim{Proposition V.2}.{The orbit space ${\rm Orb}_G(\Gamma)$ with the  
metric $d$ is naturally isomorphic to ${\rm Teich}(\Gamma)$ with its  
Teichm\"uller metric $d_T$. }
\proof
Let $\Phi : {\rm Orb}_G(\Gamma) \to  
{\rm Teich}(\Gamma)$ be given by $\Phi([\mu]) = [\Gamma^\mu]$. The proof that
$\Phi$ is an isometry is standard.
\hfill\endproof

We see at once that ${\rm Teich}\,(\Gamma)$ is a path-connected space. 
A {\it Beltrami path} in $M(\Gamma)$ is a path $t\mapsto \mu_t$ such
that for almost every $z\in{\cal D}$, the path $t\mapsto \mu_t(z)$ is
a geodesic in $\disk$. This definition is equivariant with respect to
the action of the group $G$, so we have Beltrami paths in ${\rm
Orb}_G(\Gamma)$, and therefore Beltrami pahts in ${\rm
Teich}(\Gamma)$ also, joining any two points in the space.

We will need a germ version of the qc-distance called the {\it
Julia-Teichm\"uller distance}.
Let $\Gamma |{\cal O}$ denote the restriction of (all arrows of)
$\Gamma$ to the open set ${\cal O}\subseteq {\cal U}$. Consider the
class of pairs $(\Gamma , {\cal O})$ where ${\cal O}$ is an open
neighborhood  of the Julia set of $\Gamma$. Define an
equivalence relation on such pairs as follows: $(\Gamma_1,{\cal
O}_1)\sim (\Gamma_2,{\cal O}_2)$ iff $\Gamma_1|({\cal O}_1\cap {\cal
O}_2)\equiv \Gamma_2|({\cal O}_1\cap {\cal O}_2)$. The resulting
equivalence classes are the germs of holomorphic commuting pairs
around their limit sets. The germ of $\Gamma$ {\it up
to conformal equivalence} will be denoted by $\<\Gamma\>$. Now let 
$$
d_{JT}(\<\Gamma_1\>,\<\Gamma_2\>)={1\over 2}\inf_{H} K(H)
$$
where $H:(\Gamma_1',{\cal O}_1')\simeq (\Gamma_2',{\cal O}_2')$
ranges over all possible quasiconformal conjugacies between all
representatives of both germs.
\smallskip
\ni {\it Definition 5.\enspace} $d_{JT}(\<\Gamma_1\>,\<\Gamma_2\>)$ is
the {\it Julia-Teichm\"uller distance} between $\<\Gamma_1\>$ and
$\<\Gamma_2\>$.
\smallskip
The Julia-Teichm\"uller distance is clearly a pseudo-metric. However,
it is not clear yet that it is a {\it metric}, cf. section VIII. More
importantly, it is weakly contracted by the renormalization operator.
We also verify without difficulty that the map from ${\rm
Teich}(\Gamma)$ to the space of germs up to analytic conjugacy given
by $[\Gamma_\mu] \mapsto \<\Gamma_\mu\>$ is distance-nonincreasing.

\bigskip
\ni\centerline{\smc \enspace VI. The factoring of renormalization compositions }
\smallskip
\ni
In order to develop complex bounds for renormalization of holomorphic
commuting pairs, we shall need a generalization of the so-called {\it sector
theorem} of Sullivan [S$_3$, \S 5], which we proceed to state.

Given $a, b \in {\real}$ with $a < b$, let $S(a,b)$ be the class of all  
schlicht mappings $\phi$ defined on ${\bf C}(I_\phi) = {\bf C}\setminus
({\real}\setminus I_\phi)$, where $I_\phi \supseteq (a,b)$ is some open
interval, which preserve both half-planes ${\bf C}^+, {\bf C}^-$ and are such
that $\phi((a,b)) = (a,b)$. We refer to $I_\phi$ as the {\it base} of $\phi
\in S(a,b)$: it is the largest interval containing $(a,b)$ restricted to which
$\phi$ is a homeomorphism into the reals. 
An element $A \in S(a,b)$ is a {\it left $\alpha$-root} (where  
$0 < \alpha < 1$) if there exists $a_0 \le a$ such that $A(z) = u.(z -  
a_0)^\alpha + v$, where $u,v \in {\real}$ and the branch of $z \mapsto  
(z - a_0)^\alpha$ are uniquely determined by the requirements $A(a) = a$, 
$A(b) = b$. The point $a_0 \in {\real}$ is the {\it pole} of $A$.  
Right roots are defined similarly. 
Given a bounded interval $J \subseteq {\real}$ and some $\lambda > 0$, we 
denote by $J^\lambda$ the closed interval centered at the midpoint 
of $J$ whose length is $(1 + \lambda)$-times the length of $J$. 
 
\sproclaim{Theorem VI.1}.{Let there be given $A_i, B_i \in S(a,b)$, for  
$i = 1,2, \cdots, m$, and constants $\lambda, K, s > 0$ and $0 < \alpha 
< 1$ satisfying
\itemitem{(a)} Each $A_i$ is a left $\alpha_i$-root with $\alpha_i 
\le \alpha$ and pole at $a_i$, where $a_1 = a$ and $a_i < a$ for all  
$i \ge 2$;
\itemitem{(b)} There exists a finite sequence of stopping times $1 =  
i_0 < i_1 < \cdots < i_q = m$ with $i_{n+1} - i_n \le s$ such that, 
setting $d_n = \min\{|a_i - a| : i_n \le i < i_{n+1}\}$, the inequality 
$\sum_{j \ge n} d_j^{-1} \le Kd_n^{-1}$ holds for all $n$;
\itemitem{(c)} The following holds for all $i \ge 2$: if $I_i$ is the base of
$B_i$ then $B_i(I_i) \supseteq [a_i, b]^\lambda$ and setting $J_i =
B_i^{-1}([a_i,b])$ then $J_i^\lambda \subseteq I_i$.}.
{\ni Under these assumptions, there exists 
a positive angle $\theta = \theta(\alpha, s, K, \lambda)$ such that 
the image of the upper half-plane by the composition $A_mB_m  
\cdots A_iB_i \cdots A_1B_1$ is contained in the sector $0 \le {\rm arg} 
(z - a) \le \pi - \theta$. }

A complete proof of this powerful tool is given in [dF$_2$].
Now, considering the long renormalization compositions of 
a critical commuting pair in the Epstein class, we 
would like to break them up into factors that will satisfy the hypotheses 
of Theorem VI.1 after affine rescaling. This is accomplished 
at the end of this section. 
 
Let $f : S^1 \to S^1$ be an orientation preserving homeomorphism with 
irrational rotation number $\rho(f) = [r_0, r_1, \cdots, r_n, \cdots]$.  
As in section I, given $x \in S^1$ and $k \ge 0$, let 
$I_k(x) \subseteq S^1$ be the unique closed interval with endpoints 
$x$ and $f^{q_k}(x)$ containing $f^{q_{k+2}}(x)$. For a distinguished 
point $c \in S^1$, we shall write $I_k$ instead of $I_k(c)$. 
Fix some large $n$, and consider the ordered collection 
of intervals ${\cal B} = \{f^i(I_n) : 1 \le i \le q_{n+1} - 1\}$. 
These intervals have pairwise disjoint interiors. For $k = 0, 1,  
\cdots, n+1$, let $j_k$ be the largest $j \ge 1$ such that 
$f^i(I_n) \cap \I_k = \emptyset$ for $1 \le i < j$. Observe that 
$j_0 = 1$. 
 
\proclaim{Lemma VI.2}.{For $1 \le k \le n+1$, we have $j_k = q_k$ 
if $n-k$ is odd, while $j_k = q_k + q_{k+1}$ if $n-k$ is even.}
 
\proof
For all $k$ in that range, either $I_n \subseteq I_k$ 
or $I_n$ is adjacent to $I_k$, depending on whether $|n - k|$ is even or 
odd. The lemma follows, then, from the dynamical interpretation 
of $\{q_i\}_{i \ge 0}$ as a sequence of return times.  
\hfill\endproof
 
Let us consider the {\it blocks} ${\cal B}_k = \{f^i(I_n) \in {\cal B} :  
j_{k-1} \le i < j_{k+1}\}$, for $k = 1, 2, \cdots, n$. Notice that 
${\cal B} = \bigcup_{k=1}^n {\cal B}_k$ and that ${\cal B}_k \cap 
{\cal B}_{k+2} = \emptyset$ for $1 \le k \le n-2$. These blocks correspond 
roughly to what Sullivan calls {\it epochs} in [S$_1$]. 
Let the {\it scale} of $f^i(I_n) \in  
{\cal B}$ be the largest $k \ge 1$, if any, such that $f^i(I_n) \subseteq 
I_{k-1} \setminus \I_{k+1}$, and let it be equal to zero otherwise. 
An element of ${\cal B}_k$ is a  
{\it $k$-marked interval} if its scale is equal to $k$. We also call  
an element of ${\cal B}_1$ a {\it 0-marked interval} if its scale is  
zero and it precedes all 1-marked intervals in the forward dynamical 
order of ${\cal B}$. 
We denote by ${\cal M}_k$ the collection of all $k$-marked intervals, 
for $k = 0, 1, \cdots, n$. It is not difficult to see that ${\cal M}_0$ 
has either 0 or $r_0 - 1$ elements, depending on whether $n$ is odd or 
even, respectively. 
 
\proclaim{Lemma VI.3}.{If $1 \le k \le n$, then $r_k \le  
{\rm card}({\cal M}_k) \le r_k(r_{k+1} + 1)$, and moreover
${\rm card}({\cal M}_k) = r_k$ whenever $k \equiv n$ (mod~2). }

\proof
Observe that the intervals $J_0 = f^{q_{k-1}} 
(I_k)$, $J_s = f^{sq_k}(J_0)$, $s = 1, 2, \cdots, r_k - 1$, constitute 
a partition of $I_{k-1} \setminus I_{k+1}$ modulo endpoints. 
Suppose $i < j$ are such that the intervals $f^i(I_n), f^j(I_n)$ belong 
to ${\cal B}_k$ and are both in the same $J_s$. Then they are $k$-marked 
by definition, and by Lemma VI.2 either (a) $j - i \le q_{k+1} - q_{k-1}$ 
or (b) $j - i \le r_{k+1}q_{k+1} + r_kq_k$, depending on whether $n-k$ is 
even or odd. Since $f^{j-i}(J_s) \cap \J_s \ne \emptyset$, we have  
$f^{j-i}(I_k) \cap \I_k \ne \emptyset$ as well, which implies $j - i \ge 
q_{k+1}$. In case (a) this yields a contradiction, and so each  
$J_s$ contains at most one element of ${\cal M}_k$, i.e. ${\rm card} 
({\cal M}_k) \le r_k$. In case (b) it shows that the number of 
elements of ${\cal M}_k$ in each $J_s$ is at most the smallest integer 
greater than $(r_{k+1}q_{k+1} + r_kq_k)/q_{k+1}$, which is $r_{k+1} +1$, 
and so ${\rm card}({\cal M}_k) \le r_k(r_{k+1} + 1)$. 
In either case we have $J_s \supseteq f^{q_{k-1} + (s+\epsilon_k) 
q_k}(I_n)$ for $s = 0, 1, \cdots, r_k - 1$, where $\epsilon_k$ is the  
remainder of $n-k$ modulo 2. Since such images of $I_n$ are in ${\cal B}_k$ 
by Lemma VI.2, they are in ${\cal M}_k$ too, hence ${\rm card}({\cal M}_k) \ge
r_k$. \hfill\endproof 

Now let $f$ be a critical circle mapping,  $c$ its critical point, so
that the bounded geometry results of section I are valid for $f$. More
precisely, we assume the following axioms.

\smallskip 
\ni{\it Axiom 1.} There exists $0<\sigma<1$ such that the inequality 
$|I_{n+1}(f^ic)| \le \sigma |I_{n-1}(f^ic)|$ holds for all $n \ge 1$ and all 
$i \in \bz$. 
\smallskip 
\ni{\it Axiom 2.} There exists $\lambda > 0$ such that the following 
holds for all $n \ge 1$: if $0 < i < i+j \le q_{n+1} - 1$ and $J \supseteq 
f^i(I_n)$ is the largest interval restricted to which $f^j$ is a diffeo onto 
its image then $[f^i(I_n)]^\lambda \subseteq J$. 
\smallskip
 
Both axioms are straightforward consequences of Theorem I.3
(interpreted directly for circle mappings). The second axiom is in
fact obtained from the Koebe principle for distortion of cross-ratios,
cf. [MS, Ch.~VI].
 
Let the {\it polar-ratio} of a non-degenerate 
interval $J$ with respect to a point $x$ be the number $P(x,J) = {\rm dist} 
(x, J)/|J|$. 
Observe that, under a map with bounded cross-ratio distortion, polar-ratios 
do not decrease by more than a multiplicative factor depending only on the 
cross-ratio distortion of the map. 
 
\proclaim{Lemma VI.4}.{There exist constants $C > 0$ and $\mu > 1$, 
depending only on constant $\sigma$ of Axiom 1, such that $P(c,
f^i(I_n)) \ge C\mu^{n-k}$ for each interval $f^i(I_n)$ whose scale is
equal to $k$. } 
 
\proof
Let $x = f^ic$; all intervals written $[a,b]$ in  
this proof will be contained in $S^1 \setminus \{x\}$. 
We assume that $n-k$ is even; the odd case is similar. Since $f$ is
topologically conjugate to the corresponding rotation, we have

\itemitem{(a)} $f^{-q_j}(x) \in [f^{q_{j-1}}(x), f^{q_{j+1}}(x)]$ for all 
$j \ge 1$;

\itemitem{(b)} if $f^i(x) \in [f^{q_{k-1}}(x), f^{q_{k+1}}(x)]$ then  
$f^{-i}(x) \in [f^{-q_{k-1}}(x), f^{-q_{k+1}}(x)]$. 
 
\ni Putting these two facts together we get $c = f^{-i}(x) \in
I_{k-2}(x)\setminus I_{k+2}(x)$, and in particular $d(c,f^i(I_n))\ge 
|I_{k+2}(x)|$. On the other hand, since $n-k$ is even, we have
$$
I_n(x)\subseteq I_{n-2}(x)\subseteq \cdots \subseteq I_{k+2}(x)
\ ,
$$
Applying Axiom 1 for $x = f^ic$ and using a
{\it telescoping} trick, we deduce that 
$$
P(c, f^i(I_n)) \ge \big({1\over{\sqrt\sigma}}\big)^{n-k-2}
\ ,
$$
and this proves the Lemma.
\hfill\endproof
 
Now we are ready to exhibit the promised 
factoring of the $n$-th renormalization of a critical commuting pair 
$\zeta =(\xi, \eta)$ in the Epstein class, for all sufficiently large $n$.
We know that there exist open intervals $\widetilde I_\xi \supseteq  
\xi(I_\xi)$ and $\widetilde I_\eta \supseteq \eta(I_\eta)$, as well as 
symmetric schlicht mappings $h_\xi^{-1} : {\bf C}(\widetilde I_\xi) \to 
{\bf C}$, $h_\eta^{-1} : {\bf C}(\widetilde I_\eta) \to {\bf C}$ such that $\xi \equiv 
h_\xi \circ Q$ and $\eta \equiv h_\eta \circ Q$, where $Q$ denotes the  
cubic polynomial $z \mapsto z^3$. Renormalizing $\zeta$ enough times if
necessary, we may assume also, by Theorem I.4, that $\gamma(I_\gamma)$  
sits inside $\widetilde I_\gamma$ with universal space around it ($\gamma =
\xi, \eta$). In particular, each restriction $h_\gamma^{-1}  
|\gamma(I_\gamma)$ has universally bounded cross-ratio distortion (cf. 
observation preceding Lemma VI.4). 
 
Consider the successive renormalizations of $(\xi, \eta)$ without 
rescaling: $(\xi_0, \eta_0) = (\xi, \eta)$ and $(\xi_{n+1}, \eta_{n+1}) 
= (\eta_n, \eta_n^{r_n} \circ \xi_n)$, for all 
$n \ge 0$.
We concentrate on the case $n$ even (and large), 
and within it we show how to achieve the desired factoring only for $\xi_n$; 
the other cases are similarly handled.
By Lemma I.2, we have the hybrid representation
$$
\xi_n (x)\,=\, f^{q_{n-1}-1}\circ \xi (x)
$$ 
for all $x$ in $I_n$, where $f=f_\zeta : I\to I$ is the circle
mapping associated to $\zeta$, and $I_j = I_j(c)$ for $c = 0$,
the critical point.  
We examine here the diffeomorphic part of $\xi_n^{-1}$, namely the
composition
$$ \widehat\xi_n = (f^{q_{n-1}-1})^{-1} : f^{q_{n-1}}(I_n) \to
f(I_n) \ . 
\eqno (13)
$$ 
Let $(13)$ be written as a word in $\xi^{-1}$, $\eta^{-1}$. A factor 
$\gamma^{-1}$ in this composition is called a {\it left} or a {\it right} 
factor, according to whether $\gamma = \eta$ or $\gamma = \xi$.
Each such factor, remember, has a further decomposition  
$\gamma^{-1} \equiv Q^{-1}\circ h_\gamma^{-1}$. A {\it left root} is  
the part of a left factor corresponding to $Q^{-1}$; {\it right roots} 
are similarly defined. The $h_\gamma^{-1}$ are called $h$-{\it factors}. 
A left root is said to be $k$-marked 
($k$ necessarily even) if the domain of its left factor 
is some $J \in {\cal M}_k$. Similarly, a right root is said to be  
$k$-marked ($k$ necessarily odd) if the domain of its  
right factor is some $J \in {\cal M}_k$. 
Here $k$ ranges from 0 up to $n-1$. Lemma VI.3 bounds the number of  
$k$-marked roots for any such $k$ in terms of the combinatorics of the 
rotation number $\rho(f) = \rho(\zeta)$. 
Organize all marked left roots in the composition giving $\widehat\xi_n$ 
by their order of appearance from right to left in that composition, 
and call them successively $\widehat A_1, \widehat A_2, \cdots, \widehat A_m$. 
In this order, first come the $(n-2)$-marked roots, then come the  
$(n-4)$-marked roots, and so on, and $\widehat A_m$ is the very  
last factor in that composition. We have 
$m = (r_0 - 1) + r_2 + \cdots + r_{n-2}$, after Lemma VI.3. 
In $(13)$, let $\widehat B_1$ be the sub-composition going from the first 
factor on the right up to and including the left-most factor before 
$\widehat A_1$, which is precisely the $h$-factor associated to the 
left-root $\widehat A_1$. Also, for $j = 2, 3, \cdots, m$, 
let $\widehat B_j$ be the sub-composition running strictly between  
$\widehat A_{j-1}$ and $\widehat A_j$. In this new notation, $(13)$ 
becomes
$$
\widehat\xi_n = \widehat A_m \circ \widehat B_m \circ \cdots \circ 
\widehat A_j \circ \widehat B_j \circ \cdots \circ \widehat A_1  
\circ \widehat B_1\ . \eqno (14)
$$ 
Let $T_2 \subset I$ be the largest open interval containing $f^2(I_n)$ 
restricted to which $f^{q_{n-1}-2}$ is a diffeo onto its image. 
Remembering that $\eta$ is defined well to the right of $\xi(0)$
(which is the right endpoint of $f(I_n)$), let
$(a, b) = T_1 = \eta^{-1}(T_2)$ and put also $T_i =    
f^{i-2}(T_2)$ for $i = 3, \cdots, q_{n-1}$. Then each $T_i$ contains 
the corresponding $f^i(I_n)$ plus definite, {\it beau} space on both 
sides, after Axiom 2. 
Let all factors $\widehat A_j, \widehat B_j$ be rescaled 
via the affine, orientation preserving maps taking the relevant $T_i$'s 
back onto $(a,b)$. Call the rescaled mappings $A_j, B_j$, respectively: 
these are now elements in the class $S(a,b)$.
We are ready to state and prove the promised factoring of renormalization
compositions.  
\bigskip
$$
\psannotate{
\psboxto(11.5cm;0cm){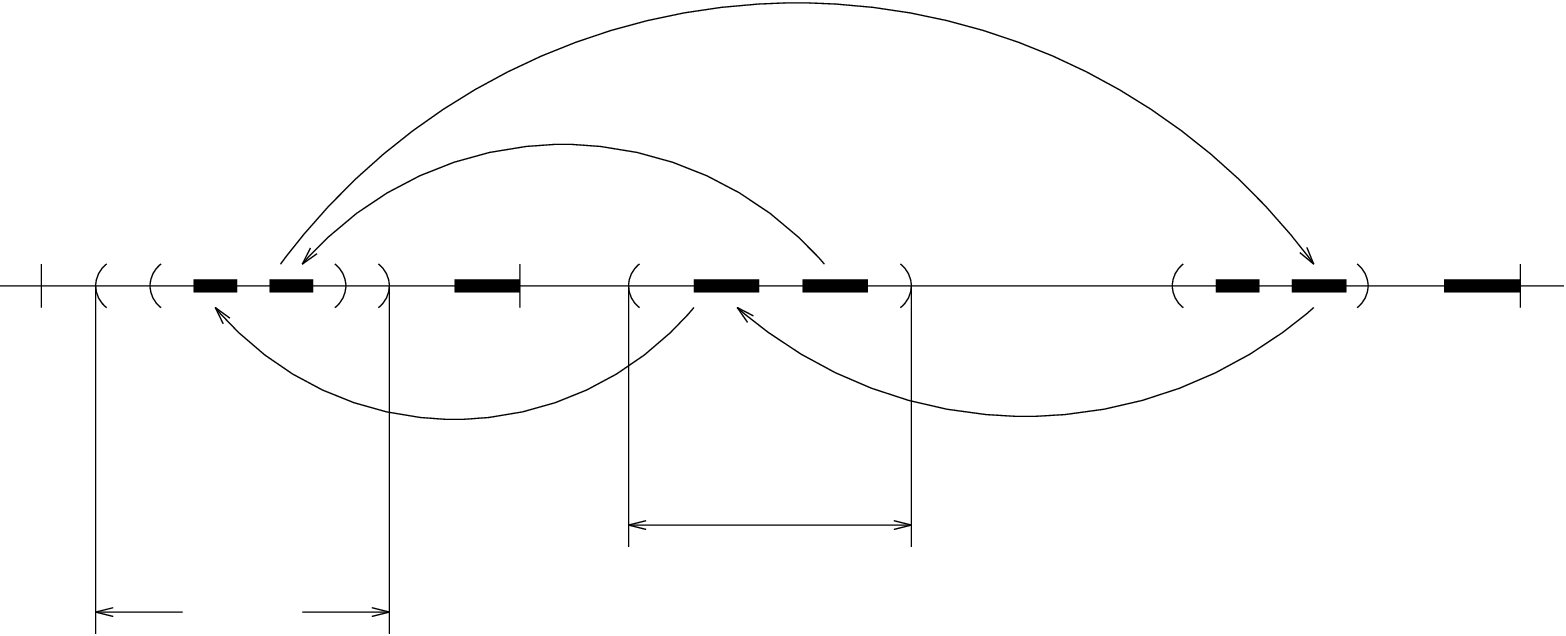}}
{\fillinggrid\at{2.5\pscm}{4\pscm}{$\eta(0)$}
             \at{9.75\pscm}{4\pscm}{$0$}
             \at{23.9\pscm}{4\pscm}{$\xi(0)$}
             \at{8.8\pscm}{5.7\pscm}{$I_n$}
             \at{22.7\pscm}{5.7\pscm}{$f\,I_n$}
             \at{16\pscm}{9\pscm}{$\xi$}
             \at{9.6\pscm}{7.5\pscm}{$\eta$}
             \at{17\pscm}{2.5\pscm}{$f^{q_k-2}$}
             \at{3.95\pscm}{0.12\pscm}{$J$}
             \at{10.4\pscm}{0.5\pscm}{$I_{k-1}\!\setminus\! I_{k+1}$}
             \at{6.5\pscm}{2.5\pscm}{$\eta$}
}
$$

\centerline{\it Figure 3}
\medskip
 
\proclaim{Theorem VI.5}.{If the rotation number $\rho(\zeta)$ is of 
bounded combinatorial type and $n$ is sufficiently large, then the rescaled
composition $A_m \circ B_m \circ \cdots \circ A_j \circ B_j \circ \cdots \circ
A_1 \circ B_1$ satisfies all hypotheses of the generalized sector theorem, and
the bounds involved are beau. }
 
\proof
Since each $A_j$ is a left $\alpha_j$-root with $\alpha_j = {1 \over 3}$ 
(i.e., a cubic root), and the pole of $A_1$ is in dynamical correspondence 
with $a$, assumption (a) of Theorem VI.1 is satisfied.
Since the number of marked left-roots at each scale is uniformly bounded 
by the hypothesis on $\rho(\zeta)$ and Lemma VI.3, the bounded 
gap condition of assumption (b) is fulfilled if we group the roots  
together by scales. 
Now, let $a_j$ be the pole of $A_j$, and let $i$ be such that  
$f^i(I_n)$ is the marked left interval corresponding to $A_j$. Combining 
the observation following Axiom 2 with the remarks preceding the statement of
this theorem and the definition of polar ratio, we obtain
$$
|a_j - a| \ge C_0|a-b|P(c,f^i(I_n)) \ ,
$$
for a certain {\it beau} constant $C_0$. Therefore, by Lemma VI.4, $|a_j - a|$  
grows exponentially with $n - k$, where $k$ is the scale of $A_j$, and 
this takes care of the series condition of assumption (b). 
It remains to check whether assumption (c) holds.
For all $j \ge 2$ we may write
$$
\widehat B_j = h_j \circ (f^{q_k - 1})^{-1} \ ,
$$
provided $k$ is the scale of $\widehat A_j$, 
where $h_j$ is the $h$-factor associated to $\widehat A_j$. Thus, we have 
the situation depicted in Figure 3 (there are two cases,  
depending on whether the scale of $\widehat A_{j-1}$, 
the preceding marked left-root, is equal to $k$ or $k+2$; 
Figure 3 illustrates the former). 
By Axiom 2, $f^{q_k - 1}$ is a diffeomorphism on an interval $J \supseteq  
[f(I_{k-1})]^\lambda$. Hence assumption (c) 
is indeed verified if we take into account that: (i) $J$ is in dynamical 
correspondence with the base of $B_i$; (ii) $f(0)=\eta(0)$, the left  
endpoint of $f(I_{k-1})$, is in dynamical correspondence with the pole 
of $A_j$, and (iii) $h_j$ is a map of {\it beau} bounded cross-ratio
distortion. The case $j = 1$ is similarly proved. 
Since all bounds involved are {\it beau}, we are done. \hfill\endproof

\bigskip
\ni\centerline{\smc VII.\enspace The complex bounds} 
\smallskip
\ni
We recall Sullivan's {\it sector inequality}. Given real numbers $u<u_*$
and $0<\vartheta\le {\pi\over 3}$, consider in the upper half-plane the
truncated sector 
$$
{\cal S}(u,u_*,\vartheta)\;=\;{\bf C}^+\cap \big\{arg\,(z-u)<\pi -\vartheta\big\}\cap
\big\{arg\,(z-u_*)> {{2\pi}\over 3}\big\}
\ .\eqno{(15)}
$$
Also, let $w'<u<u'<u_*<u_*'<w_*'$ and let $w$, $w_*$ be points in ${\cal
S}(a,b,\vartheta)$. Assume $C>1$ is a constant such that all non-zero distances
between points in $\{u,u_*,w,w_*,u',u_*',w',w_*'\}$ lie between $C^{-1}$ and $C$.
Moreover, denote by ${\cal N}_R$ the portion in the upper half-plane
of the Poincar\'e neighborhood of radius $R$ of the geodesic $(w',w_*')$
in ${\bf C}((w',w_*'))$.

\proclaim{Lemma VII.1}.{There exists $r>0$ depending only on $\vartheta$
and $C$ such that the following holds. For each $R>r$ there exists
$R_0>R$ such that, if $\psi :{\cal N}_{R_0}\to {\cal S}(u,v,\vartheta)$ is
univalent and maps $u$, $u_*$, $w$, $w_*$ respectively to $u'$, $u_*'$, $w'$,
$w_*'$, then $\psi({\cal N}_R)\subseteq {\cal N}_{R/2}$. \hfill\endproof}

Now, let us define the {\it conformal type} of a holomorphic
commuting pair to be the modulus of the annulus determined by the
inner and outer domains of its bow-tie. Let us also say that an
irrational number has {\it combinatorial type bounded by $N>0$} if
the convergents of its continued-fraction development are bounded by
$N$. 

\proclaim{Theorem VII.2}.{Given a positive integer $N$, there exists
$\tau=\tau (N)>0$ with the following property. If $\zeta$ is a critical 
commuting pair with rotation number of combinatorial type bounded by
$N$ and if $\zeta$ either belongs to some Epstein class or extends to a
holomorphic commuting pair, then for all $n$ sufficiently large,
${\cal R}^n\zeta$ extends to a holomorphic commuting pair
$\Gamma_n(\zeta)$ with geometric boundaries and conformal type
bounded from below by $\tau$.} 

\proof If $\zeta$ extends to a holomorphic commuting pair, then $\zeta$ is
analytically equivalent to a critical commuting pair in the Epstein class.
Thus one may assume that $\zeta$ is in the Epstein class already. Let $a$,
$\theta$ and $n$ be as in Theorem VI.5. To define $\Gamma_n(\zeta)$, one first
constructs its bow-tie $({\cal D}_n, {\cal O}_{\xi_n}, {\cal O}_{\eta_n},
{\cal O}_{\nu_n})$. For ${\cal D}_n^+$ one takes a large Poincar\'e
neighborhood ${\cal N}_R$ as above, with $w'$ and $w_*'$ to satisfy the
conditions below and $R$ given by Lemma VII.1. Then one takes
$$
{\cal O}_{\xi_n}^+\;=\; \xi_n^{-1}({\cal
D}_n^+)\;=\;(f^{q_{n-1}-1}\circ\xi)^{-1}({\cal D}_n^+)
\ .
$$
where $f$ is as in section VI (and where $n$ is assumed to be even). By Theorem
VI.5, $\xi({\cal O}_{\xi_n}^+)$ is 
contained in a sector in the upper half-plane with angle $\theta$ on the left.
Since $\xi^{-1}$ is a right cubic-root factor, it follows that ${\cal
O}_{\xi_n}^+$ lies within a truncated sector of type (15), where $\vartheta
=\theta /3$, $u_*$ is the origin and $u$ is a point in dynamical
correspondence with $a$. Therefore by the sector inequality the closure of
${\cal O}_{\xi_n}$ lies well within ${\cal D}_n$. One constructs ${\cal
O}_{\eta_n}$ by a similar procedure, so that its closure is also contained in
${\cal D}_n$. One then completes the bow-tie following the method in section IV.
Now let ${\cal U}_n$ be the inner domain of this bow-tie, and let $J_n= {\cal
U}_n\cap \real$. The points $w'$ and $w_*'$ had to be chosen from the start so
that $[w',w_*']$ contains the interval $J_n$ with {\it beau} space on both
sides. The boundary of ${\cal U}_n$ consists of finitely many analytic arcs
meeting at finitely many corners with internal angles $\ge \pi /3$. The total
number of corners is bounded in terms of the number $k$ of critical values of
$\xi_n$ found in $[w',w_*']$. If the size of this interval is commensurable
with the size of $J_n$, then by the real a-priori bounds (cf.~Theorems I.3 and
I.4) one has $k$ bounded in terms of $N$ only. Therefore, $\partial {\cal
U}_n$ is a $K$-quasicircle with $K=K(N)$, by Ahlfors' characterization of
quasicircles in [A$_1$]. Finally, the minimum distance between the boundaries
of the annulus ${\cal D}_n\setminus {\cal U}_n$ depends only on the space of
$J_n$ inside $[w',w_*']$. Therefore its modulus is controlled, by the 
Teichm\"uller inequalities (9). \hfill\endproof
A major consequence of this theorem is the following compactness
property enjoyed by renormalization. 

\proclaim{Corollary VII.3}.{Given a positive integer $N$, there exists $B_N>0$
with the following property. If $\zeta_1$ and $\zeta_2$ are
critical commuting pairs with the same irrational rotation
number of combinatorial type bounded by $N$, and if each of them either
belongs to an Epstein class or extends to a holomorphic commuting pair,
then for all sufficiently large $n$ we have, in the notation of
Theorem VII.1,
$$
d_{JT}(\<\Gamma_n(\zeta_1)\>,\<\Gamma_n(\zeta_1)\>)\,\le\,
d_{QC}\big(\Gamma_n(\zeta_1), \Gamma_n(\zeta_2)\big)\,\le\, B_N
\ . 
$$
}
\proof Combine Theorem VII.2 with the pull-back argument of Theorem III.1.
\hfill\endproof 

One of the keys to renormalization contraction, this corollary states that if
we renormalize a finite Beltrami arc sufficiently many times, then its endpoints,
no matter how far apart in the Julia-Teichm\"uller sense, come within a fixed
distance from each other.
\smallskip
\ni
{\it Remark.\enspace } In a very recent work, Yampolsky [Ya] proved
without using the sector theorem that the modulus of ${\cal
D}_n\setminus {\cal U}_n$ is always bounded from below, {\it
independently of the combinatorics}, thereby obtaining complex bounds
for critical circle maps with arbitrary rotation numbers.

\bigskip
\ni\centerline{\smc VIII.\enspace A Cantor repeller and its Riemmann surface
lamination}  
\smallskip
\ni
A {\it Cantor repeller} consists of two collections $\{ D_i\}_{0\le i\le n}$
and $\{\Delta_j\}_{0\le j\le m}$ of topological disks in the plane and a
surjective holomorphic map $\phi :\bigcup_{i=0}^{n}D_i \to
\bigcup_{j=0}^{m}\Delta_j$ such that 
(a) in each collection, any two disks have pairwise disjoint
closures; (b) each disk of the first collection is compactly contained in some
disk of the second collection; (c) each $\phi |D_i$ is schlicht onto
$\Delta_j$ for some $j$. The invariant limit set $K_\phi =
\bigcap_{n=0}^{\infty} \phi^{-n}(\bigcup_{j} 
\Delta_j)$ is a Cantor set, hence the name, and the restriction $\phi |K_\phi$
is topologically conjugate to a certain shift of finite type (cf.~[Bo]). We
call this shift the {\it topological type} of our Cantor repeller.
Writing $U$ and $V$ for domain and range of $\phi$, we represent the
Cantor repeller by $(U,\phi ,V)$. We also say that a Cantor repeller
is {\it linear} if $U$, $V$ and $\phi$ are symmetric about the real line.

Now let $\Gamma$ be a holomorphic commuting pair with irrational rotation
number, and  let ${\cal K}_\Gamma$ be its (connected) limit set. We want to
show how to extract from within $\Gamma$ a certain Cantor repeller
that, off of its limit set, turns out to be conformally conjugate to
$\Gamma$ in the vicinity of its small dynamical interval. More
precisely, we have the following result. 

\proclaim{Theorem VIII.1}.{Let $\Gamma$ be a holomorphic commuting
pair having a connected filled-in limit set ${\cal K}_\Gamma$. Then
there exist an open set ${\cal O}$ and a linear Cantor repeller
$(U,\phi ,V)$ of fixed topological type such that 
\itemitem{(a)} $\phi |U^+$ is conjugate to $\Gamma |
({\cal O}\setminus {\cal K}_\Gamma )$ by a holomorphic map $H$;
\itemitem{(b)} For each $n\ge 0$ there exists an open neighborhood
${\cal V}_n$ of the small dynamical interval of $\Gamma$ contained in
${\cal O}$ such that $H(\phi^{-n}(U^+)) = {\cal
V}_n\setminus {\cal K}_\Gamma$.}
\proof
Since $\widehat{\bf C} \setminus {\cal K}_\Gamma$ is simply-connected,
let $\Phi : \widehat{\bf C} \setminus \overline{\disk} \to \widehat{{\bf
C}} \setminus {\cal K}_\Gamma$ be the 
Riemann-mapping, normalized to be symmetric about the real axis
and fixing $\infty$.
Consider the simply-connected regions $V_0 = {\cal O}_\xi^+ \setminus 
{\cal K}_\Gamma$, $V_1 = {\cal O}_\nu^+ \setminus ({\cal O}_\xi \cup {\cal
O}_\eta \cup {\cal K}_\Gamma)$, 
$V_2 = {\cal O}_\eta^+ \setminus {\cal K}_\Gamma$, $V_3 = {\cal O}_\xi^-
\setminus {\cal K}_\Gamma$,
$V_4 = {\cal O}_\nu^- \setminus ({\cal O}_\xi \cup {\cal O}_\eta \cup
{\cal K}_\Gamma)$,
$V_5 = {\cal O}_\eta^- \setminus {\cal K}_\Gamma$, and also $W_0 =
\Delta^+ \setminus {\cal K}_\Gamma$, $W_1 = \Delta^- \setminus {\cal
K}_\Gamma$.  
Let ${\cal O}_i =\Phi^{-1}(V_i)$ and $\Omega_j = \Phi^{-1}(W_j)$. 
Now, let $\psi : \bigcup_{i=0}^5 {\cal O}_i \to \Omega_0 \cup \Omega_1$ be
the mapping $\Phi^{-1} \circ F \circ \Phi$, where $F$ is the shadow
of $\Gamma$ (cf. section II). Then each 
$\psi | {\cal O}_i$ is schlicht onto either $\Omega_0$ or $\Omega_1$, 
depending on whether $i$ is even or odd.
The intervals $I_i = \partial{\cal O}_i \cap \partial\disk$  
are pairwise disjoint, and each $\psi | {\cal O}_i$ carries the corresponding
$I_i$ onto $(\partial\disk)^+$ if $i$ is even, and onto
$(\partial\disk)^-$ if $i$ is odd. 
Next, let $\Delta_0 \subseteq {\bf C}^+$ be a Jordan domain such that (a)
$\Delta_0$ is symmetric about $\partial\disk$ under geometric inversion;
(b) $\overline\Delta_0 \setminus \overline{\disk} \subseteq \Omega_0$; 
(c) $\overline{{\cal O}}_0 \cup \overline{{\cal O}}_1 
\cup \overline{{\cal O}}_2 \subseteq \Delta_0$. Let $\Delta_1
\subseteq {\bf C}^-$ be similarly defined, and
consider the inverse mappings $\psi_i = f^{-1} : \Omega_0 
\to {\cal O}_i$ for $i$ even. Restrict each $\psi_i$ to $\Delta_0 
\setminus \overline{\disk}$ and then extend
the corresponding restriction to $\Delta_0$ by Schwarz's reflection.
Continue to denote these extensions by the
same names, and let $D_i = \psi_i(\Delta_0)$, $i$ even. Define 
$D_i$ for $i$ odd in similar fashion, using $\Delta_1$.
Then $\psi : D_0 \cup \cdots \cup D_5 \to \Delta_0 \cup \Delta_1$ is
a Cantor repeller. 
Moreover, each of the open sets
$$
{\cal V}_n= int\Big(\Phi\,\Big(\psi^{-n}\big(\bigcup_{i=0}^{5}
D_i\setminus \overline{\disk} \big)\Big)\cup {\cal K}_\Gamma\Big) \ ,
$$
contains the small dynamical interval of $\Gamma$ (cf. remark at the
end of section III). Note that ${\cal
V}_{n+1}\subseteq {\cal V}_n$, and take ${\cal O}={\cal V}_0$.
Finally, let $M$ be a fractional linear transformation taking
$\widehat{\bf C}\setminus \disk$ onto ${\bf C}^+$ and, say, the point
$-1$ to $\infty$. Then let $H=\Phi\circ M^{-1}$, $\phi=M\circ\psi
\circ M^{-1}$ and let $U$ and $V$ be domain and range of $\phi$. This puts
our Cantor repeller in linear form and proves (a) and (b).
\hfill\endproof

The Cantor repeller given by this theorem is determined only up to
holomorphic conjugacy. It is a weak analogue of the
Douady-Hubbard external class for holomorphic commuting pairs, even
though a straightening theorem is missing.

We define the {\it germ} of a Cantor repeller around its
limit set just as in the case of holomorphic commuting pairs (cf.
section V), and write $\<\phi\>$ for the germ of $(U,\phi ,V)$ up to
holomorphic equivalence. We
also define the corresponding Julia-Teichm\"uller 
distance between such germs, but call it instead the {\it germ
distance} to avoid confusion, and denote it by $d_G$. Theorem VIII.1
implies that, if $\Gamma_1$ and $\Gamma_2$ are qc-conjugate and
$\phi_1$ and $\phi_2$ are the corresponding Cantor repeller maps,
then $d_{JT}(\<\Gamma_1\>,\<\Gamma_2\>)\ge d_G(\<\phi_1\>,\<\phi_2\>)$.

\sproclaim{Theorem VIII.2}.{Let $\Gamma_1$ and $\Gamma_2$ be
holomorphic commuting pairs and let $\phi_1$ and $\phi_2$ be the
corresponding Cantor repeller maps. Consider the following five statements.
\itemitem{(a)} $\Gamma_1$ and $\Gamma_2$ have the same germ up to holomorphic
equivalence;
\itemitem{(b)} $d_{JT}(\<\Gamma_1\>,\<\Gamma_2\>)\,=\,0$;
\itemitem{(c)} $d_G(\<\phi_1\>,\<\phi_2\>)\,=\,0$;
\itemitem{(d)} $\phi_1$ and $\phi_2$ are analytically conjugate on the line.
\itemitem{(e)} There exists $k\ge 0$ such that $d_{JT}(\<{\cal
R}^k\Gamma_1\>,\<{\cal R}^k\Gamma_2\>)\,=\,0$.
}.
{\ni Then we have (a) $\Rightarrow$ (b) $\Rightarrow$ (c) $\Leftrightarrow$
(d) $\Rightarrow$ (e).} 
\proof 
We only prove that (c) $\Rightarrow$ (d) $\Rightarrow$ (e) and refer to [MS,
Ch.VI, \S 4, Corollary 1] for the other implications. Assume (c) holds.
For each $\varepsilon >0$ we find neighborhoods
$U_1^{\varepsilon}\supseteq K_{\phi_1}$ and
$U_2^{\varepsilon}\supseteq K_{\phi_2}$ and a
$(1+\varepsilon)$-quasiconformal map
$h_\varepsilon:U_1^{\varepsilon}\to U_2^{\varepsilon}$ conjugating
$\phi_1$ and $\phi_2$. This conjugacy restricts to a quasi-symmetric
map on the line that is H\"older with exponent $1-O(\varepsilon)$.
Hence, the {\it scaling functions} of $K_{\phi_1}$ and $K_{\phi_2}$ differ
on corresponding points of their dual Cantor sets by $O(\varepsilon)$.
Letting $\varepsilon\to 0$ we deduce that $K_{\phi_1}$ and
$K_{\phi_2}$ have the same scaling function. Therefore $\phi_1$ and
$\phi_2$ are analytically conjugate on some real-line neighborhoods
of both Cantor sets, by [S$_3$]. Now, if (d) holds, then both repellers are
analytically conjugate on neighborhoods of their limit sets in the complex
plane. Therefore by Theorem VIII.1(b), $\Gamma_1$ and $\Gamma_2$ are
analytically conjugate on neighborhoods of their small dynamical intervals. By
the complex bounds, these neighborhoods contain the inner domains of the
bow-ties of all sufficiently large renormalizations of both pairs, so (e)
follows.\hfill\endproof 

We note a further relationship between the Julia-Teichm\"uller
distance on holomorphic commuting pairs and the germ distance on
Cantor repellers. 
\proclaim{Proposition VIII.3}.{Let $\Gamma$ and $\Gamma'$ be
topologically equivalent holomorphic commuting pairs. For each
$\varepsilon >0$ there exists $N>0$ such that for all $n\ge N$
$$
d_{JT}(\<\Gamma_n\>,\<\Gamma'_n\>)\;\le\;
d_G(\<\phi\>,\<\phi'\>)\,+\,\varepsilon  \ ,
$$
where $\Gamma_n ={\cal R}^n\Gamma$ and $\Gamma'_n ={\cal R}^n\Gamma'$ (cf.
Proposition II.3).} 
\proof This follows at once from Theorem VIII.1(b) together with the
complex bounds given by Corollary VII.2. \hfill\endproof

Now we associate to the germ of a Cantor repeller
around its limit set a compact Riemann surface lamination
in the sense of Sullivan (cf.~[S$_1$]).
Roughly, this will be the space of backward branch orbits (or 
{\it threads}) of points in any deleted neighborhood of $K_\phi$ 
factored by the equivalence relation determined by the dynamics of $\phi$
itself.
Recall that a {\it Riemann surface lamination (or RSL -) 
structure} on a Hausdorff topological space $X$ consists of a maximal atlas 
$\{(U_\alpha, \varphi_\alpha)\}$ covering $X$ such that (a) each 
$\varphi_\alpha$ maps the corresponding $U_\alpha$ homeomorphically
onto $D_\alpha \times T_\alpha$, where $D_\alpha \subseteq {\bf C}$ 
is a disk and $T_\alpha$ is a Hausdorff space, and
(b) each ovelapping homeo $\varphi_\beta \circ \varphi_\alpha^{-1} : 
\varphi_\alpha(U_\alpha \cap U_\beta) \to \varphi_\beta(U_\alpha \cap
U_\beta)$ is of the form $(z,t) \mapsto (\psi_t(z), \phi(t))$ with 
$\psi_t$ holomorphic for each $t$.
Provided with such a structure, the space
$X$ is called a {\it Riemann surface lamination}.
One defines the leaves of a Riemann surface lamination just as in the
case of foliations; leaves come with obvious intrinsic
structures making them into Riemann surfaces in a natural way.
A lamination (qc-) morphism $X \to Y$ between two Riemann surface laminations
is a continuous map that sends leaves of $X$ into leaves of $Y$ and is 
holomorphic (resp. quasiconformal) on leaves.
Every locally trivial bundle
over a Riemann surface with totally disconnected fiber above each 
element of the base has a unique RSL-structure on the total
space that makes the projection map into a lamination morphism. This
fact yields the following lemma on inverse limits.

\proclaim{Lemma VIII.4}.{Let $\ldots \to
S_n \buildrel g_n \over \longrightarrow S_{n-1} \to \ldots 
\buildrel g_1 \over \to S_0$ be an inverse system where each $S_n$
is a free union of Riemann surfaces and each $g_n$ is a proper
holomorphic covering map, let $S_\infty$ be its topological inverse
limit and let $\pi_n : S_\infty \to S_n$, $n \ge 0$ be the canonical
projections. Then $S_\infty$ is a locally compact,
2$^{\hbox{nd}}$-countable space and has a unique RSL-structure making
each $\pi_n$ into a lamination morphism. \hfill\endproof}

We have also the following facts.

\proclaim{Lemma VIII.5}.{Let $X_0 \buildrel g_0 \over 
\longrightarrow X_1 \buildrel g_1 \over \longrightarrow \cdots
\to X_n \buildrel g_n \over \longrightarrow \cdots$ be a direct
system where each $X_n$ is a Riemann surface lamination and each 
$g_n$ is an open, injective lamination morphism. Then the direct limit
space $X^\infty$ has a unique RSL-structure making the canonical maps
$\rho_n : X_n \to X^\infty$, $n \ge 0$, as well as the direct limit
map $g^\infty : X^\infty \to X^\infty$, into open, injective
lamination morphisms.\hfill\endproof}

\proclaim{Lemma VIII.6}.{Let $X_n = X$ and $g_n = g$ for each
$n \ge 0$ in Lemma VIII.5, where $X$ is locally compact and first
countable and $g$ acts discontinuously on $X$.
Then the orbit space $X^\infty / \<g^\infty\>$ has a unique
RSL-structure for which the canonical projection $X^\infty \to X^\infty/
\<g^\infty\>$ is a lamination morphism. Moreover, if $\phi$ has a relatively
compact fundamental domain in $X$ then $X^\infty/\<g^\infty\>$ is a 
compact space.\hfill\endproof}

Following Sullivan, we say that a Riemann surface lamination is {\it
hyperbolic} if each of its leaves is covered by the disk.

\proclaim{Theorem VIII.7}.{For each Cantor repeller $(U,\phi , V)$
there exists a compact, hyperbolic Riemann surface lamination
$L(U,\phi ,V)$ such that
\itemitem{(a)} If $(U, \phi ,V)$ and $(\widetilde U, \widetilde\phi
,\widetilde V)$ represent the same germ,
then we have a lamination isomorphism $L(U, \phi ,V) \cong L(\widetilde U, 
\widetilde\phi ,\widetilde V)$; 
\itemitem{(b)} Every qc-conjugacy
$(U_1, \phi_1, V_1) \sim (U_2, \phi_2, V_2)$ induces a 
qc-isomorphism of laminations $L(U_1, \phi_1, V_1) \cong L(U_2,
\phi_2, V_2)$.} 
\proof Let $V_0 = V\setminus K_\phi$. Consider the 
inverse system
$$ 
\cdots \buildrel \phi \over \longrightarrow \phi^{-n}V_0 \buildrel
\phi \over \longrightarrow \phi^{-(n-1)}V_0 \to \cdots \to \phi^{-1}V_0 
\buildrel \phi \over \longrightarrow V_0 \eqno (16)
$$
and its sub-system
$$ 
\cdots \buildrel \phi \over \longrightarrow \phi^{-(n+1)}V_0 \buildrel
\phi \over \longrightarrow \phi^{-n}V_0 \to \cdots \to \phi^{-2}V_0 
\buildrel \phi \over \longrightarrow \phi^{-1}V_0 \eqno (17)
$$
and let $V_\infty$, $V_\infty'$ be their
inverse limit spaces. Both are Riemann surface laminations
by Lemma VIII.4, and since $(17)$ is cofinal in $(16)$, we have a
lamination isomorphism $\Phi : V_\infty \to V_\infty'$. 
The inclusions $\phi^{-(n+1)}V_0 \subseteq \phi^{-n}V_0$ yield an open,
injective lamination morphism $\Psi : V_\infty' \hookrightarrow
V_\infty$, so $g = \Psi\circ\Phi$ is a map with these properties also. 
Since the multivalued map $\phi^{-1}$ acts discontinuously on $V_0$,
$g$ acts discontinuously on $V_\infty$.
Applying Lemmas VIII.5 and VIII.6 to the direct system 
$V_\infty \buildrel g \over \longrightarrow V_\infty \buildrel g \over 
\longrightarrow \cdots \to V_\infty \buildrel g \over \longrightarrow 
\cdots$, we get the desired $L(U,\phi ,V)$ as the orbit space of the direct
limit map $g^\infty$ acting on the direct limit space $V^\infty$.
This lamination is compact because $V_0\setminus \phi^{-1}V_0$ is a
relatively compact fundamental domain for the action of $\phi^{-1}$
on $V_0$. All leaves are hyperbolic Riemann surfaces, cf. geometric
description of the dyadic pair-of-pants lamination in [S$_5$, Example~3]. Parts
(a) and (b) are straightforward. \hfill\endproof 

\bigskip

\ni\centerline{\smc \enspace IX. Renormalization convergence }
\smallskip
\ni
We are now in a position to use the Teichm\"uller theory of Riemann surface
laminations, introduced by Sullivan in the appendix to [S$_1$], in order to
prove that renormalization contracts the Julia-Teichm\"uller distance. For
details on the unproved assertions in this section, see the book by de Melo
and van Strien [MS, Ch.~VI].

In [S$_1$], Sullivan defined Beltrami vectors and
quadratic differentials on a compact Riemann surface lamination $X$
as cross-sections of suitable tensor bundles over $X$. Thus,
a Beltrami vector $\mu$ locally on each flow-box chart
$(D_\alpha \times T_\alpha, \psi_\alpha)$ is a Borel measurable
function $\mu_\alpha : D_\alpha \times T_\alpha \to {\bf C}$ satisfying (a) $\mu_\alpha
(\cdot, t) \in L^\infty(D_\alpha)$ for each $t \in T_\alpha$, and the
map $t \mapsto \mu_\alpha(\cdot, t)$ is continuous if we provide $L^\infty
(D_\alpha)$ with the weak topology; (b) if $\psi_{\alpha\beta}$ denotes the
chart transition $\psi_\beta \circ \psi_\alpha^{-1}$ and we write $\psi_{\alpha\beta}
= (\psi_{\alpha\beta}^z, \psi_{\alpha\beta}^t)$ then we have
$$
\mu_\alpha = {\overline{\partial\psi_{\alpha\beta}^z} \over
\partial\psi_{\alpha\beta}^z} 
\mu_\beta \circ \psi_{\alpha\beta} \ . 
$$
A Beltrami coefficient on $X$ is an essentially bounded Beltrami
vector with essential norm less than one.
Sullivan also defined quadratic differentials on $X$ as the
corresponding dual objects. More precisely, a quadratic differential
$\varphi$ on $X$ is an assignment of a $\sigma$-finite measure class $[m_\alpha]$ 
to the transversal $T_\alpha$ of each flow-box chart satisfying (a) the
transversal components $\psi_{\alpha\beta}^t$ of chart transitions are 
absolutely continuous as maps $(T_\alpha, [m_\alpha]) \to (T_\beta, [m_\beta])$;
(b) for each choice of representative $m_\alpha \in [m_\alpha]$ there exists
a measurable function $\varphi_\alpha : D_\alpha \times T_\alpha \to
{\bf C}$ such that, on overlappings
$$
\varphi_\alpha = \varphi_\beta \circ
\psi_{\alpha\beta}\,\left[{\partial\psi_{\alpha\beta}^z  
\over \partial z}\right]^2\,{\rm Jac}(\psi_{\alpha\beta}^t) \ ,  \eqno (18)
$$
where the Jacobian is measured with
respect to the measures $m_\alpha$ and $m_\beta$; (c) each $\varphi_\alpha$ is 
integrable with respect to the product measure $dz\,d\overline z\,dm_\alpha$ on
$D_\alpha \times T_\alpha$. It follows from this definition that there exists
a well-defined measure $d|\varphi|$ associated to a quadratic differential on
$X$. Its expression on a given chart is
$|\varphi_\alpha|\,dz\,d\overline z\,dm_\alpha$  
for each choice of measure $m_\alpha$, and from $(18)$ any two choices
differ by the Jacobian of the identity with respect to both transversal
measures, i.e. by their Radon-Nikodym derivative. If the
total mass 
$$
|\varphi| = \int_X d|\varphi| 
$$
is finite, we say 
that $\varphi$ is an {\it integrable} quadratic differential, and
$|\varphi|$ is the {\it norm} of $\varphi$.
A quadratic differential is said to be {\it holomorphic} if it is 
holomorphic on almost all leaves with respect to the transversal measure
class that it defines.
The {\it Teichm\"uller norm} of a Beltrami vector $\omega$ is 
$$
|\omega|_T = \sup|\int_X\omega\varphi| \ ,
$$
where the supremum 
is taken over all integrable holomorphic quadratic differentials of
norm $|\varphi| = 1$.
These definitions are set-up so that the natural pairing 
in the second member
is well-defined.
Given $\varepsilon \ge 0$, a Beltrami vector
$\omega$ on $X$ is called {\it $\varepsilon$-extremal} if $|\omega|_\infty \le
(1 + \varepsilon) |\omega|_T$.

Sullivan used an elegant ergodic argument to prove a {\it
generalized Gr\"otzsch inequality} relating the dilatation of
qc-isomorphisms on $X$ which are leafwise isotopic to the identity
and the {\it metric up to a multiple} given by an integrable
quadratic differential of norm one on $X$ (cf.~[MS, Ch.VI, \S 7] for
details).
He then used this inequality to prove the almost geodesic principle
below. Let us say that a Beltrami coefficient $\mu$ on $X$ is
{\it dynamical} if, integrating $\mu$ via the MRMT along the leaves of $X$,
we get a transversally continuous map $X\to X$ that is qc on leaves,
i.e. a lamination qc-morphism. If $\mu$ is dynamical, let $c(\mu)$ be
the RSL-structure on $X$ given by $\mu$.
A {\it dynamical} Beltrami vector $\omega$ on $(X,c(\mu))$
is one for which there exists a (unique) path of dynamical Beltrami
coefficients $\mu_t$, $t\ge 0$, with $\mu_0=\mu$ and tangent to $\omega$ at
$t=0$, such that for all $t$ the smallest maximal dilatation of a qc-morphism
between $c(\mu)$ and $c(\mu_t)$ is equal to $e^{2t}$.
We write $c_t(\omega)=c(\mu_t)$ and call the path $\mu_t$ the {\it Beltrami
ray} of $\omega$ at $\mu$.
\smallskip
\ni {\it Example.\enspace} Let $(U,\phi ,V)$ be a Cantor repeller and
let ${\cal L}_\phi$ be the lamination of Theorem VIII.7. Given any
$\phi$-invariant Beltrami vector $\widetilde\omega$ on the Riemann
surface $V_0=V\setminus K_\phi$, we pull it back via the natural
projection to the inverse limit space $V_\infty$ and then project it
down to a $g^{\infty}$-invariant Beltrami vector on the direct limit
space $V^{\infty}$, thus getting a dynamical Beltrami vector $\omega$
on ${\cal L}_\phi=V^{\infty}/\<g^{\infty}\>$. By [S$_5$], all
dynamical Beltrami vectors on ${\cal L}_\phi$ arise in this way (they are
precisely the transversally locally constant ones, in Sullivan's
terminology). In particular, the Teichm\"uller norm of $\omega$ can
be computed by pairing $\widetilde\omega$ on a fundamental domain for
$\phi$, such as $V\setminus U$, with holomorphic quadratic
differentials there. 
\smallskip

Now the {\it almost geodesic principle} can be stated as follows.

\proclaim{Theorem IX.1}.{Given $\varepsilon$, $L>0$, there exists
$\delta = \delta (\varepsilon , L)>0$ such that the following holds.
If $\mu$ is a dynamical Beltrami coefficient on a compact hyperbolic Riemann
surface lamination $X$, $\omega$ is a $\delta$-extremal dynamical
Beltrami vector on $X$ at $\mu$ and $\{\psi_t\}_{0 \le t \le 1}$
is a leafwise qc-isotopy between $(X, c(\mu))$ and $(X,
c_{\ell}(\omega))$, then we have $L\le K(1+\varepsilon)$, where $K$
is the maximal dilatation of $\psi_1$ and $\ell={1\over
2}\log{L}$. Moreover, $\delta\to 0$ as $\varepsilon\to
0$.\hfill\endproof}  

We need a converse to this theorem within the realm of the example above.
Consider the germ $\<\Gamma\>$ of a holomorphic commuting pair
and let $\<\phi\>$ be the germ of the corresponding Cantor repeller
constructed in Theorem VIII.1. Let $(\Gamma ,{\cal O})$ be
a representative of $\<\Gamma\>$ and let $\mu$ be a
Beltrami differential with domain ${\cal O}$. We call
$\mu$ {\it admissible} for $\<\Gamma\>$ if $\mu$
is $\Gamma$-invariant and vanishes a.e. on the limit set ${\cal
K}_\Gamma$. If $\mu$ defined on ${\cal O}$ is an admissible Beltrami
coefficient for $\<\Gamma\>$, let $h_{\mu}$ be a suitably normalized
qc-mapping with dilatation $\mu$ and let $\Gamma({\cal O},\mu)$ be the
holomorphic dynamical system generated by the mappings
$$
h_\mu\circ (\gamma |({\cal O}_\gamma\cap {\cal O}))\circ h_\mu^{-1}
\ ,
$$
where $\gamma=\xi ,\eta ,\nu$.
If $\sigma$ is another admissible Beltrami coefficient for $\<\Gamma\>$
defined on ${\cal O}'$, we say that $\mu$ and $\sigma$ are {\it equivalent} if
$\Gamma({\cal O}\cap {\cal O}',\mu)$ is analytically conjugate to
$\Gamma({\cal O}\cap {\cal O}',\sigma)$. 
We then let $|\mu |_{JT}=\inf{|\sigma|_\infty}$, where $\sigma$ runs through
all admissible Beltrami coefficients for $\<\Gamma\>$ that are equivalent to
$\mu$ (caution: this is a non-linear norm). We also say that $\mu$ is
{\it $\varepsilon$-efficient} if $|\mu|_\infty\le (1+\varepsilon)|\mu|_{JT}$.
Similarly, admissible Beltrami coefficients or vectors
for $\<\phi\>$ are those defined in the domain of a representative
of the germ of $\phi$ which are $\phi$-invariant
and vanish a.e. on the limit set $K_\phi$. The definitions we have just given
can be repeated here. We denote by $|\cdot |_G$ the
non-linear norm of admissible coefficients for $\<\phi\>$ that corresponds to
$|\cdot |_{JT}$. Observe that Theorem VIII.1 sets up a correspondence between
admissible Beltrami differentials for $\<\Gamma\>$ and admissible Beltrami
coefficients for $\<\phi\>$. Admissible objects for $\<\phi\>$ are precisely
those that lift to dynamical objects in the lamination ${\cal L}_\phi$.
Observe that a globally $\Gamma$-invariant $\widetilde\mu$ in the sense of
section V is admissible, and the definitions have been arranged so that
$$
d_{JT}(\<\Gamma\>,\<\Gamma^{\widetilde\mu}\>)\;=\;{1\over
2}\,\log{{1+|\widetilde\mu|_{JT} }\over 
{1-|\widetilde\mu|_{JT}}} \ .
$$
Likewise, if $\mu$ is the admissible coefficient for $\<\phi\>$ corresponding
to $\widetilde\mu$ and $\phi^\mu$ is the corresponding Cantor repeller, we
have (cf. section VIII)
$$
d_G(\<\phi\>,\<\phi^\mu\>)\;=\;{1\over
2}\,\log{{1+|\mu|_G }\over 
{1-|\mu|_G}} \ .
\eqno{(19)}
$$

Now Sullivan's {\it coiling lemma} can be stated as follows.
\proclaim{Theorem IX.2}.{Given $\varepsilon'>0$ and $0<d\le 1$, there
exists $\theta=\theta(\varepsilon',d)>0$ such that, if $\omega$ is a
an admissible Beltrami vector for $\<\phi\>$ and the admissible
Beltrami coefficient $\mu_s=s\omega$ is $\theta$-efficient for some
$0<s< d|\omega|_\infty^{-1}$, then $\omega$ is
$\varepsilon'$-extremal.\hfill\endproof }
Renormalization without rescaling acts on admissible Beltrami vectors
in a natural way. Thus, if $\widetilde\mu$ is admissible for $\<\Gamma\>$ and
defined on ${\cal O}$, let its $n$-th renormalization $\widetilde\mu_n$ be the
restriction of $\widetilde\mu$ to ${\cal U}_n\cap {\cal O}$, where ${\cal U}_n$
is the inner domain of the bowtie of $\Gamma_n\,=\, {\cal R}^n(\Gamma)$.
Then let $\mu_n$ be an admissible Beltrami coefficient for $\<\phi_n\>$ that
corresponds to $\widetilde\mu_n$, where $\<\phi\>$ is the Cantor repeller germ
associated to $\Gamma_n$. By Corollary VII.3, if $\Gamma$ is of bounded
combinatorial type then for every $n$ sufficiently large ${\cal U}_n$ is
contained in ${\cal O}$. Therefore the holomorphic commuting pair
$\Gamma_n^{\widetilde\mu_n}$ is well-defined (cf. section V), and we have
$\Gamma_{n+1}^{\widetilde\mu_{n+1}}\,=\,{\cal R}\,\Gamma_n^{\widetilde\mu_n}$,
for all sufficiently large $n$.

We have at last the main {\it renormalization contraction theorem} that follows.

\proclaim{Theorem IX.3}.{Let $\<\Gamma\>$ and $\<\Gamma'\>$ be
germs of holomorphic commuting pairs with the same rotation number of
bounded combinatorial type and the same height. Then the distance 
$d_{JT}({\cal R}^n\<\Gamma\>,{\cal R}^n\<\Gamma'\>)$ converges to
zero as $n\to \infty$.}
\proof
We argue as in the proof of [MS, Ch.~VI, Thm. 8.3].
By Proposition VIII.3, it suffices to show that renormalization contracts the
germ distance $d_G$ between the germs of corresponding repellers $\<\phi\>$,
$\<\phi'\>$.
Let $\widetilde\mu$ be the
Beltrami coefficient of a qc-conjugacy $(\Gamma ,{\cal O}) \to (\Gamma' ,{\cal
O}')$ which is $\varepsilon_1$-efficient for some $\varepsilon_1$ to be
specified below. Let $\mu$ be the corresponding admissible Beltrami
coefficient for $\<\phi\>$. Then $\mu$ is $\varepsilon_1$-efficient also. Let
$\omega = \mu /|\mu |_\infty$ and lift $\omega$ to a dynamical Beltrami vector
$\widehat\omega$ on the lamination ${\cal L}_\phi$ with its standard
structure. Note that $\widehat\omega$ is $\varepsilon_1$-extremal. For each
$t>0$, let $\mu (t)=(\tanh t)\omega$ and note that the lifted path
$\widehat\mu (t)$ is the Beltrami ray of $\widehat\omega$ at zero. 

Let $B$ be the constant in the complex bounds (Corollary VII.3) and fix a
constant $M$ so large that $M>B^{-1}\log{(1+2e^B)}$. Take $L=2e^{MB}$, $\ell
={1\over{2}}\log{L}$ and $0<\varepsilon <1$ to be specified later, and then
choose $\varepsilon_1 =\delta (\varepsilon ,L)$ using the almost
geodesic principle for $\widehat\omega$. We get $K(1+\varepsilon)\ge L$, where
$K$ is the smallest dilatation of all qc-morphisms leafwise isotopic to the
identity in ${\cal L}_\phi$ between the standard structure and
$c(\widehat\mu(\ell))$. Therefore
$$
K\;=\;
{{1+|\mu(\ell)|_G }\over {1-|\mu(
\ell)|_G}}\;\ge\; {{2e^{MB}}\over {1+\varepsilon }}\;\ge\; e^{MB}\ .
\eqno{(20)}
$$
But if $n$ is sufficiently large, then by Corollary VII.3 we have
$$
d_{JT}(\<\Gamma_n\>,\<\Gamma_n^{\widetilde\mu_n(\ell)}\>)
\le B
$$
and therefore
$$
|\mu_n(\ell)|_G\;\le\;{{e^B-1}\over {e^B+1}}\ .
\eqno{(21)}
$$
Combining (20) and (21), we get
$$
|\mu(\ell)|_G\;\ge\; {{e^{MB}-1}\over {e^{MB}+1}}\;\ge\; {{e^B}\over
{e^B-1}}\, |\mu_n(\ell)|_G \ ,
\eqno{(22)}
$$
by our choice of $M$. Now let $k>1$ be such that $k(1-e^{-B})\le 1-e^{-2B}$;
we can choose $k$ as close to 1 as we like. Then
{\it either} $|\mu |_\infty >k|\mu_n|_\infty$, in which case
$$
|\mu_n|_G\;\le\; {{1+\varepsilon_1}\over k}\, |\mu |_G \ ,
\eqno{(23)}
$$
{\it or} $|\mu |_\infty \le k|\mu_n|_\infty$, in which case (22) gives us
$|\mu_n(\ell)|_G\le (1-e^{-2B})|\mu_n(\ell)|_\infty$. In this last
case, applying the coiling lemma to the admissible Beltrami vector
$\mu_n(\ell)$ with $\varepsilon'=e^{-2B}$ and
$d=1$, we 
see that there exists $0<\theta <1$ depending only on $B$ such that, for 
$0< \tanh{t} < 1$, the admissible Beltrami coefficients
$\mu_n(t)$ cannot be $\theta$-efficient. In particular, taking
$t={\rm arctanh}\,{|\mu|_\infty}$, we have
$$
|\mu_n|_G\;\le\; (1-\theta)|\mu_n|_\infty
\;\le\;k(1-\theta)(1+\varepsilon_1)\, |\mu |_G \ .
\eqno{(24)}
$$
Now choose $k$ first so that $k(1-\theta)<1$ and then $\varepsilon$ so small
that
$$
\lambda_1= {\rm
max}\,\{k^{-1}(1+\varepsilon_1),\,k(1-\theta)(1+\varepsilon_1)\}<1 \ .
$$
Then in both (23) and (24) we have $|\mu_n|_G\le \lambda_1|\mu |_G$. Using (19),
we deduce that
$$
d_{G}(\<\phi_n\>,\<\phi_n^{\mu}\>)\;\le\;\lambda_2\,d_G(\<\phi\>,\<\phi^\mu\>)\ ,
$$
for some $0<\lambda_2<1$, and this is the desired contraction.\hfill\endproof

This contraction of the Julia-Teichm\"uller distance results in stronger forms
of renormalization convergence.
Following [MS, Ch.~VI, \S 8], we define strong convergence as follows.

\ni{\it Definition 6.\enspace} A sequence $g_n: W\to {\bf C}$, where
$W\subseteq {\bf C}$ is compact, {\it converges strongly} to $g: W\to {\bf
C}$ if there exist an open neighborhood ${\cal O}$ of $W$ and  
holomorphic extensions $G: {\cal O}\to {\bf C}$ of $g$ and $G_n: {\cal O}\to
{\bf C}$ of $g_n$, for all but finitely many $n$, such that $G_n$ converges to
$G$ uniformly in ${\cal O}$. 

Notice that if $W$ is an interval on the line, say, then strong convergence of
$g_n$ to $g$ in $W$ implies $C^k$-convergence for all $k<\infty$.
Now let ${\cal B}^\omega (N)$ (resp. ${\cal B}^3 (N)$) be the class of {\it
normalized} real-analytic (resp. $C^3$-smooth) critical commuting pairs with
irrational rotation number of combinatorial type bounded by $N$.
We say that a sequence $\zeta_n$ in ${\cal B}^\omega (N)$ converges
strongly to $\zeta$ in ${\cal B}^\omega (N)$ if both $\eta_n-\eta$ and
$\xi_n-\xi$ converge strongly to zero. This last condition makes sense
because the first implies that $\eta_n(0)\to \eta (0)$ and therefore any fixed
neighborhood of $[\eta (0),0]$ contains $[\eta_n(0), 0]$ for all sufficiently
large $n$. Strong convergence of a sequence of holomorphic
commuting pairs, or of a sequence of Cantor repellers, can be similarly
defined. We need the following statement, which is Lemma 8.4 of [MS, Ch.~VI, \S 8].

\proclaim{Lemma IX.4}.{Given $\varepsilon >0$ and $R_1>1$, there exist $\delta
>0$ and $R_2>R_1$ with the following property. If $h$ is a $(1+\delta)$-qc homeo
that fixes $0$ and $1$ and whose domain and range contain the disk of radius
$R_2$ about zero, then $|h(z)-z|<\varepsilon$ for all $|z|<R_1$.\hfill\endproof}

Now, given $\zeta\in {\cal B}^\omega (N)$, let
$\zeta_n\,=\,(\xi_n,\eta_n)\,=\,{\cal R}^n\zeta$ be the {\it normalized}
renormalizations of $\zeta$. 

\proclaim{Theorem IX.5}.{Let $\zeta,\zeta'\in {\cal B}^\omega (N)$ be
critical commuting pairs of the same combinatorial type which either
belong to some Epstein class or extend to holomorphic commuting pairs. Then
$\xi_n-\xi'_n$ and $\eta_n-\eta'_n$ converge strongly to zero.}

\proof By Corollary VII.3, if  $n$ is sufficiently large then $\zeta_n$ and
$\zeta'_n$ extend to normalized holomorphic commuting pairs $\Gamma_n$ and
$\Gamma'_n$ with conformal types bounded from below.
Let ${\cal U}_n$ and ${\cal U}'_n$ be the inner domains of the bow-ties of
$\Gamma_n$ and $\Gamma'_n$. Also, for each $k>0$, let ${\cal U}_{n,k}\subseteq
{\cal U}_n$ be the linear copy of ${\cal U}_{n+k}$ corresponding to the $k$-th
renormalization of $\Gamma_n$ without rescaling, and let ${\cal
U}'_{n,k}\subseteq {\cal U}'_n$ be
similarly defined. We have ${\cal U}_{n,k+1}\subseteq {\cal U}_{n,k}$ for all
$k$, and ${\rm mod}\,({\cal U}_n\setminus {\cal U}_{n,k})\to\infty$ as
$k\to\infty$, by the complex bounds. Likewise, ${\rm mod}\,({\cal
U}'_n\setminus {\cal U}'_{n,k})\to\infty$ as $k\to\infty$.
Given $\varepsilon>0$ and $R_1>1$ so large that the disk of radius $R_1$ about
the origin contains the small dynamical intervals of $\Gamma_n$ and
$\Gamma'_n$ for all $n$, take $\delta$ and $R_2$ as in
Lemma IX.4. By Theorem IX.3, there exist $m>0$ and a
$(1+\delta)$-quasiconformal conjugacy $h:{\cal U}_m\to {\cal U}'_m$
between $\Gamma_m$ and $\Gamma'_m$. Note that the restriction of $h$
to ${\cal U}_{m,k}$ is a conjugacy between the $k$-th
renormalizations without rescaling of $\Gamma_m$ and $\Gamma'_m$. Let
$\Lambda_k$ and $\Lambda'_k$ be the linear maps that perform such
rescaling, so that $\Lambda_k({\cal U}_{m,k})\,=\,{\cal U}_{m+k}$.
Writing $H_k=\Lambda'_k\circ h\circ\Lambda_k^{-1}$, we have $H_k(0)=0$ and
$H_k(1)=1$, each $H_k$ is $(1+\delta)$-qc, and also
$$
\cases{\xi_{m+k}=H_k^{-1}\circ \xi'_{m+k}\circ H_k & {}\cr
{} & {}\cr
\eta_{m+k}=H_k^{-1}\circ \eta'_{m+k}\circ H_k & {}\cr}
\ .\eqno{(25)}
$$ 
Moreover, if ${\cal W}_k=\Lambda_k({\cal U}_m)$ and ${\cal
W}'_k=\Lambda'_k({\cal U}'_m)$, then 
$$
{\rm mod}\,({\cal W}_k\setminus [0,1])\,>\,{\rm mod}\,({\cal
W}_k\setminus {\cal U}_{m+k})\,=\,{\rm mod}\,({\cal
U}_m\setminus {\cal U}_{m,k})\to\,\infty
$$
as $k\to\infty$, and similarly for ${\cal W}'_k$. Therefore there
exists $k_0$ such that ${\cal W}_k$ and ${\cal W}'_k$ contain the
disk of radius $R_2$ about the origin for all $k\ge k_0$. By Lemma
IX.4, we have  $|H_k(z)-z|<\varepsilon$ and
$|H_k^{-1}(z)-z|<\varepsilon$ for all $|z|<R_1$, for all $k\ge
k_0$. Going back to (25), we get strong convergence as
claimed.\hfill\endproof  

\proclaim{Lemma IX.6}.{Let $(U_n,\phi_n,V_n)$ and $(U'_n,\phi'_n,V'_n)$, $n\ge
0$, be two sequences of Cantor repellers of the same topological type, and
suppose they converge strongly to $(U,\phi,V)$ and $(U',\phi',V')$,
respectively. If $d_G(\<\phi_n\>,\<\phi'_n\>)\,=\,0$ for all $n$, then 
$d_G(\<\phi\>,\<\phi'\>)\,=\,0$ also.}

\proof If $d_G(\<\phi_n\>,\<\phi'_n\>)\,=\,0$ then, by Theorem VIII.2,
there exist neighborhoods ${\cal O}_n\supseteq K_{\phi_n}$ and ${\cal
O}'_n\supseteq K_{\phi'_n}$, and an analytic homeo $h_n:{\cal O}_n\to {\cal
O}'_n$ conjugating $\phi_n$ to $\phi'_n$. Take $k>0$ large enough (depending
on $n$) so that
$$
\Delta_n\;=\;\phi_n^{-k}(V_n\setminus U_n)\subseteq {\cal O}_n
\ .
$$
Then $\Delta_n$ and $\Delta'_n\,=\,h_n(\Delta_n)\subseteq {\cal O}'_n$ are
fundamental domains for $\phi_n$ and $\phi'_n$, respectively. Now, let
$\varepsilon>0$. Since $\phi_n$ converges strongly to $\phi$, we can find a
fundamental domain $D_n$ for $\phi$ and a homeomorphism $\psi_n: D_n\to
\Delta_n$ very close to the identity which conjugates $\phi$ to $\phi_n$ on
corresponding boundaries and is $(1+\varepsilon )$-quasiconformal,
provided $n$ is sufficiently large. Similarly, we can find a fundamental
domain $D'_n$ for $\phi'$ and a $(1+\varepsilon )$-quasiconformal map
$\psi'_n: D'_n\to \Delta'_n$ conjugating $\phi'$ to $\phi'_n$, possibly by
making $n$ larger still. This gives us a conjugacy $(\psi'_n)^{-1}\circ
h_n\circ \psi_n : D_n\to D'_n$ between the fundamental domains of $\phi$ and
$\phi'$ which is $(1+\varepsilon )^2$-quasiconformal. By a simple pull-back
argument, this map extends to a qc-conjugacy with the same dilatation between
$\phi$ and $\phi'$ in full-neighborhoods of their limit sets. Therefore
$d_G(\<\phi\>, \<\phi'\>)\,=\,0$ as claimed.\hfill\endproof

\proclaim{Lemma IX.7}.{Let $\Gamma_n$ and $\Gamma'_n$, $n\ge
0$, be two sequences of holomorphic commuting pairs, and
suppose they converge strongly to $\Gamma$ and $\Gamma'$, respectively. If
$d_{JT}(\<\Gamma_n\>,\<\Gamma'_n\>)\,=\,0$ for all $n$, then there exists
$k\ge 0$ such that $d_{JT}(\<{\cal R}^k\Gamma\>,\<{\cal R}^k\Gamma'\>)\,=\,0$
also.} 

\proof This follows at once from Lemma IX.6 and Theorems VIII.1 and
VIII.2.\hfill\endproof 

\proclaim{Theorem IX.8}.{For each $n\in\bz$, let $\zeta_n$ and $\zeta'_n$ be
normalized critical commuting pairs with the same bounded 
combinatorial type and suppose they extend to holomorphic commuting
pairs $\Gamma(\zeta_n)$ and $\Gamma(\zeta'_n)$,
respectively, whose conformal types are uniformly bounded from below. If
$\Gamma(\zeta_{n+1})\,=\,{\cal R}\Gamma(\zeta_n)$ and
$\Gamma(\zeta'_{n+1})\,=\,{\cal R}\Gamma(\zeta'_n)$ for all $n$, then
$\zeta_0=\zeta'_0$.}

\proof The proof of [MS, Ch.~VI, Lemma 8.3] can be reproduced here almost
verbatim. Lemma IX.7 replaces the argument on continuity of the
Douady-Hubbard external class used in that proof.\hfill\endproof

Finally, we give a characterization of the {\it attractor} of the
renormalization operator for critical commuting pairs. Theorems A and B of the
Introduction are straightforward consequences of this last theorem, which is
the exact analogue of [MS, Ch.~VI, Theorem 1.1]. We denote
by $W^s(\zeta)$ the {\it stable set} of $\zeta\in {\cal B}^3(N)$, i.e. the set
of all $\zeta'\in {\cal B}^3(N)$ whose successive renormalizations are
$C^0$-asymptotic to those of $\zeta$.

\sproclaim{Theorem IX.9}.{Let $N$ be a positive integer. There exists a
renormalization-invariant, strongly compact set ${\cal A}\subseteq {\cal
B}^\omega (N)$ such that
\itemitem{(a)} If $\zeta\in {\cal B}^3(N)$, then the $C^3$-distance
between ${\cal R}^n(\zeta )$ and ${\cal A}$ converges to zero as
$n\to\infty$; 
\itemitem{(b)} There exist $a>0$ and $\tau>0$ such that ${\cal A}\subseteq
{\cal E}_a$ and each element of ${\cal A}$ extends to a holomorphic commuting
pair with conformal type bounded by $\tau$;
\itemitem{(c)} The restriction of ${\cal R}$ to ${\cal A}$ is a homeomorphism
topologically conjugate to the two-sided full-shift on $N$ symbols;
\itemitem{(d)} If $\zeta\in {\cal A}$ then $W^s(\zeta )$ is the set of
critical commuting pairs $\zeta'$ such that ${\cal R}^m(\zeta' )$ and ${\cal
R}^m(\zeta )$ have the same bounded combinatorial type for some $m>0$.}.{\ni
Moreover, there exists a strongly compact set ${\cal C}\supseteq {\cal 
A}$ such that (i) for any real-analytic $\zeta$ of combinatorial type bounded
by $N$ in some Epstein class, there exists $n_0(\zeta )>0$ such that ${\cal
R}^n(\zeta )\in {\cal C}$ for all $n\ge n_0(\zeta )$, and (ii) if $\zeta ,
\zeta' \in {\cal C}$ have the same bounded combinatorial type then
$\xi_n-\xi'_n$ and $\eta_n-\eta'_n$ converge strongly to zero.}

\proof Define ${\cal A}$ as the set of all $C^0$-limits of successive
renormalizations of critical commuting pairs in ${\cal B}^3(N)$. Then (a)
follows from Theorem I.4 which, combined with Theorem VII.2 and Corollary
VII.3, proves (b) also. Proceeding as in [MS, Ch.~VI, Thm.~8.4], one shows
that for each bi-infinite sequence $(\ldots ,r_{-1}, r_0,r_1,\ldots
,r_n,\ldots )$ with $r_n\in\{0,1,\ldots ,N\}$, there exists a bi-infinite
sequence $(\ldots ,\zeta_{-1}, \zeta_0,\zeta_1,\ldots ,\zeta_n,\ldots )$ of
critical commuting pairs $\zeta_n\in {\cal A}$ such that
$\rho(\zeta_n)=[r_n+1,r_{n+1},\ldots ]$ and $\zeta_{n+1}={\cal R}\zeta_n$ for
all $n$, and such sequence is unique by Theorem IX.8. Therefore
the map
$$
(\ldots ,r_{-1}, r_0,r_1,\ldots ,r_n,\ldots )\mapsto \zeta_0 \in {\cal A}
$$
is one-to-one and onto and conjugates the full-shift on $N$ symbols to the
restriction of renormalization to ${\cal A}$. If ${\cal A}$ is given the
strong topology of Definition 6, this map is continuous, hence a
homeomorphism, and this proves (c). In particular, ${\cal A}$ is strongly
compact. Finally, let ${\cal C}\subseteq {\cal B}^\omega(N)$ be the set of
critical commuting pairs that can be extended to holomorphic commuting pairs
with conformal type bounded from below by $\tau$. Then ${\cal C}$ is strongly
compact,  (i) is Theorem VII.1, and (ii) is Theorem IX.5. \hfill\endproof 

\bigskip
\removelastskip
\relax
\bigskip
\centre{\smc References}
\bigskip
\refsize
\medtype
\ref{A$_1$}{L.~Ahlfors, {\it{Lectures on Quasiconformal Mappings}}, 
Van-Nostrand (1966).}
\ref{A$_2$}{L.~Ahlfors, {\it{Conformal Invariants: Topics in 
Geometric Function Theory}}, McGraw-Hill, New-York (1973).}
\ref{B}{L.~Bers, On moduli of Kleinian groups, {\it Russian Math.
Surveys {\bf29}}, 88-102, (1974).}
\ref{Be}{A.~Beardon, {\it{The Geometry of Discrete Groups}},
Springer-Verlag, Berlin and New York, (1983).}
\ref{Bo}{R.~Bowen, {\it{On Axiom A Diffeomorphisms}}, in Regional 
Conference Series in Mathematics {\bf 35}, Providence AMS, 1977.}
\ref{C}{C.~Carath\'eodory, {\it{Theory of Functions of a 
Complex Variable}}, Chelsea, New York, vol. 2 (1964).}
\ref{Crl}{L.~Carleson, On mappings conformal at the boundary, 
{\it J. Anal. Math. {\bf 19}}, 1-13, (1967).}
\ref{CS}{P.~Cvitanovic \& B.~S\"oderberg, Scaling laws for mode lockings in
circle maps, {\it Physica Scripta {\bf 32}}, 263, (1988).}
\ref{Cv}{P.~Cvitanovic, {\it{Universality in Chaos}}, Adam Hilger, Bristol,
(1989).}
\ref{dF$_1$}{E.~de~Faria, {\it{Proof of universality for critical circle
mappings}}. Ph.~D. Thesis, CUNY, (1992).}
\ref{dF$_2$}{E.~de~Faria, On conformal distortion and Sullivan's
sector theorem. Preprint IHES/M/96/20, (1996).}
\ref{dFM$_1$}{E.~de~Faria \& W.~de~Melo, Renormalization of
circle mappings: compactness and convergence to the Epstein class.
Preprint (1996).} 
\ref{dFM$_2$}{E.~de~Faria \& W.~de~Melo, Towers of holomorphic
pairs and rigidity of circle mappings. In preparation.} 
\ref{DH}{A.~Douady \& J.~Hubbard, On the dynamics of 
polynomial-like mappings, {\it Ann. Scient. Ec. Norm. Sup. {\bf 18}},
287-343, (1985).} 
\ref{EE}{J.-P.~Eckmann \& H.~Epstein, On the existence of 
fixed points of the composition operator for circle maps, {\it Commun. 
Math. Phys. {\bf 107}}, 213-231, (1986). }
\ref{EL}{H.~Epstein \& J.~Lascoux, Analiticity properties of 
the Feigenbaum function, {\it Commun. Math. Phys. {\bf 81}}, 437-453,
(1981).} 
\ref{FKS}{M.~Feigenbaum, L.~Kadanoff \& S.~Shenker, Quasi-Periodicity
in dissipative systems. A renormalization group analysis, {\it
Physica {\bf 5D}}, 370-386, (1982).}
\ref{G$_1$}{F.~Gardiner, {\it{Teichm\"uller Theory and Quadratic 
Differentials}}, John Wiley and Sons, New York, (1987).}
\ref{G$_2$}{F.~Gardiner, On Teichm\"uller contraction, {\it Proc. of the
Am. Math. Soc. {\bf 118}}, 865-875, (1993).} 
\ref{GS}{F.~Gardiner \& D.~Sullivan, Symmetric and 
quasi-symmetric structures on a closed curve, {\it Am. J. of Math.
{\bf114}}, 683-736, (1992).}
\ref{GrS}{J.~Graczyk \& G.~\'Swi\c{a}tek, Critical circle maps near
bifurcation, {\it Commun.~Math.~Phys. {\bf 176}}, 227-260, (1991).}
\ref{H$_1$}{M.~Herman, Sur la conjugaison differentiable des 
diff\'eomorphismes du cercle a des rotations, {\it Publ.Math.I.H.E.S. 
{\bf49}}, 5-234, (1979).}
\ref{H$_2$}{M.~Herman, Conjugaison quasi-sym\'etrique des
diff\'eomorphismes du cercle et applications aux disques singuliers
de Siegel. Manuscript, (1986).}
\ref{Ha}{C.~R.~Hall, A $C^\infty$ Denjoy counterexample, {\it Ergod. Th. \&
Dynam. Sys. {\bf 1}}, 261-272, (1981).}
\ref{K}{L.~Keen, Dynamics of holomorphic self-maps of 
{\bf{C}}$^*$. In {\it Holomorphic Functions and Moduli I} (ed. D. Drasin
et al.) Springer-Verlag, New York, (1988). } 
\ref{KO}{Y.~Katznelson \& D.~Ornstein, The differentiability of
conjugation of certain diffeomorphisms of the circle, {\it Ergod. Th.
\& Dynam. Sys. {\bf 9}}, 643-680, (1989).}
\ref{Kh}{K.~Khanin, Thermodynamic formalism for critical circle
mappings. In {\it Chaos}, ed. D.~K~Campbell, American Institute of
Physics, New York, 71-90, (1990).}
\ref{KS}{K.~Khanin \& Y.~Sinai, A new proof of M.~Herman's theorem,
{\it Comm.~Math.~Phys. {\bf 112}}, 89-101, (1987).}
\ref{L$_1$}{O.E.~Lanford, Renormalization group methods for
circle mappings, {\it Nonlinear evolution and chaotic phenomena},
NATO Adv.Sci.Inst.Ser.B:Phys., {\bf 176}, Plenum, New York,
25-36, (1988).}
\ref{L$_2$}{O.E.~Lanford, Renormalization group methods for critical circle
mappings with general rotation number, {\it VIIIth International Congress on
Mathematical Physics (Marseille, 1986)}, World Sci. Publishing, Singapore,
532-536, (1987).}
\ref{Le}{O.~Lehto, {\it{Univalent functions and Teichm\"uller spaces}},
Springer-Verlag, Berlin and New York, (1986).}
\ref{LV}{O.~Lehto \& K.~Virtannen, {\it{Quasiconformal mappings}},
Springer-Verlag, Berlin and New York, (1965). }
\ref{Ly}{M.Y.~Lyubich, On typical behavior of the trajectories of a
rational mapping of the sphere, {\it Soviet Math. Dokl. {\bf 27}},
22-24, (1983).}
\ref{McM}{C.~McMullen, Self-similarity of Siegel disks and Hausdorff
dimension of Julia sets. Preprint, (1995).}
\ref{Me}{B.~Mestel, Ph.D. Dissertation. Department of Mathematics, 
Warwick University, (1985).}
\ref{Mi}{J.~Milnor, {\it Dynamics in one-complex variable: introductory
lectures}. IMS Stony Brook preprint 90/5, (1990).}
\ref{MS}{W.~de Melo \& S.~van Strien, {\it{One-dimensional Dynamics}}, 
Springer-Verlag, Berlin and New York, (1993). }
\ref{ORSS}{R.~Ostlund, D.~Rand, J.~Sethna \& E.~Siggia, Universal 
Properties of the Transition from Quasi-Periodicity to Chaos in Dissipative 
Systems, {\it Physica {\bf8D}}, 303, (1983).}
\ref{Ra$_1$}{D.~Rand, Existence, non-existence and universal 
breakdown of dissipative golden invariant tori. I. Golden critical circle 
mappings, Preprint, I.H.E.S., (1989).}
\ref{Ra$_2$}{D.~Rand, Universality and renormalisation in dynamical 
systems. In {\it New Directions in Dynamical Systems} (ed. T. Bedford and
J. W. Swift) , Cambridge University Press, Cambridge, (1987).} 
\ref{Ra$_3$}{D.~Rand, Global phase-space universality, smooth
conjugacies and renormalization, $C^{1+\alpha}$ case, {\it
Nonlinearity {\bf 1}}, 181-202, (1988).}
\ref{Ric}{S.~Rickmann, Removability theorems for quasiconformal 
mappings, {\it Ann. Ac. Scient. Fenn., {\bf49}}, 1-8, (1969).}
\ref{S$_1$}{D.~Sullivan, Bounds, quadratic differentials, and 
renormalization conjectures. In {\it Mathematics into the Twenty-first
Century}, Amer. Math. Soc. Centennial Publication, vol. 2, Amer. Math.
Soc., Providence, RI., (1991).} 
\ref{S$_2$}{D.~Sullivan, Bounded structure of infinitely 
renormalizable mappings. In {\it Universality in Chaos}, $2^{nd}$
edition, Adam Hilger, Bristol, (1989).}
\ref{S$_3$}{D.~Sullivan, Differentiable structures on fractal-like 
sets, determined by intrinsic scaling functions on dual Cantor sets,
{\it Proceedings of Symposia in Pure Mathematics {\bf48}}, (1988).}
\ref{S$_4$}{D.~Sullivan, Quasiconformal homeomorphisms and 
Dynamics I. Fatou-Julia Problem on Wandering Domains, {\it Ann. of
Math. {\bf122}}, 401-418, (1985).}
\ref{S$_5$}{D.~Sullivan, Linking the universalities of Milnor-Thurston,
Feigenbaum and Ahlfors-Bers. In {Topological Methods in Modern
Mathematics}, Proceedings of the Symposium held in honor of John
Milnor's 60th birthday, SUNY at Stony Brook, 543-564,
(1991).}
\ref{Sh}{S.~Shenker, Scaling behavior of a map of a circle onto 
itself: empirical results, {\it Physica {\bf5D}}, 405, (1982).}
\ref{St}{J.~Stark, Smooth conjugacy and renormalization for
diffeomorphisms of the circle, {\it Nonlinearity {\bf 1}}, 541-575,
(1988).}
\ref{ST}{D.~Sullivan \& W.~Thurston, Extending holomorphic 
motions, {\it Acta Math. {\bf 157}} n.3/4, 243-257, (1986).}
\ref{Sw$_1$}{G.~\'Swi\c{a}tek, Rational rotation numbers for maps of 
the circle, {\it Commun. Math. Phys. {\bf 119}}, 109-128, (1988).}
\ref{Sw$_2$}{G.~\'Swi\c{a}tek, One-dimensional maps and Poincar\'e 
metric, {\it Nonlinearity {\bf5}}, 81-108, (1992).}
\ref{Ya}{M.~Yampolsky, Complex bounds for critical circle maps.
IMS Stony Brook preprint 95/12, (1995).}
\ref{Yo$_1$}{J.-C.~Yoccoz, Conjugaison differentiable des 
diffeomorphismes du cercle dont le nombre de rotation verifie une
condition Diophantienne, {\it Ann. Sci. de l'Ec. Norm. Sup. {\bf17}},
333-361, (1984).} 
\ref{Yo$_2$}{J.-C.~Yoccoz, Il n'y a pas de contre-example de Denjoy 
analytique, {\it C. R. Acad. Sci. Paris {\bf298}} s\'erie I, 141-144,
(1984).}

\bye